\magnification=1200

\loadmsam
\loadmsbm
\loadeufm
\loadeusm
\UseAMSsymbols

\font\BIGtitle=cmr10 scaled\magstep3
\font\bigtitle=cmr10 scaled\magstep1
\font\boldsectionfont=cmb10 scaled\magstep1
\font\section=cmsy10 scaled\magstep1

\def\scr#1{{\fam\eusmfam\relax#1}}

\def\scrA{{\scr A}}
\def\scrB{{\scr B}}
\def\scrC{{\scr C}}

\def\scrF{{\scr F}}

\def\scrJ{{\scr J}}
\def\scrM{{\scr M}}
\def\scrN{{\scr N}}

\def\scrP{{\scr P}}

\def\scrS{{\scr S}}
\def\scrU{{\scr U}}

\def\scrT{{\scr T}}

\def\gr#1{{\fam\eufmfam\relax#1}}

	\def\gra{{\gr a}}
\def\grB{{\gr B}}

\def\db#1{{\fam\msbfam\relax#1}}

\def\dbA{{\db A}} 
\def\dbC{{\db C}} 
 \def\dbF{{\db F}}
\def\dbG{{\db G}} \def\dbH{{\db H}}

 \def\dbN{{\db N}}
 
\def\dbQ{{\db Q}} \def\dbR{{\db R}}
\def\dbS{{\db S}}

 \def\dbZ{{\db Z}}

\def\bigg{\text{big}}
\def\Ker{\text{Ker}}
\def\der{\text{der}}
\def\Sh{\hbox{\rm Sh}}

\def\Res{\text{Res}}
\def\ab{\text{ab}}
\def\an{\text{an}}
\def\ad{\text{ad}}

\def\sh{\text{sh}}
\def\Gal{\text{Gal}}
\def\Hom{\text{Hom}}
\def\End{\text{End}}
\def\Spec{\text{Spec}}
\def\Spf{\text{Spf}}
\def\stc{\text{stc}}

\def\sc{\text{sc}}

\def\Lie{\text{Lie}}

\def\leaderfill{\leaders\hbox to 1em
     {\hss.\hss}\hfill}
\def\nspace{\lineskip=1pt\baselineskip=12pt\lineskiplimit=0pt}

\def\finishproclaim{\par\rm
     \ifdim\lastskip<\medskipamount\removelastskip
     \penalty55\medskip\fi}
\def\proof{\par\noindent {\it Proof:}\enspace}
\def\references#1{\par
  \centerline{\boldsectionfont References}\medskip
     \parindent=#1pt\nspace}
\def\Ref[#1]{\par\hang\indent\llap{\hbox to\parindent
     {[#1]\hfil\enspace}}\ignorespaces}
\def\Item#1{\par\smallskip\hang\indent\llap{\hbox to\parindent
     {#1\hfill$\,\,$}}\ignorespaces}
\def\ItemItem#1{\par\indent\hangindent2\parindent
     \hbox to \parindent{#1\hfill\enspace}\ignorespaces}

\def\Le{{\mathchoice{\,{\scriptstyle\le}\,}
  {\,{\scriptstyle\le}\,}
  {\,{\scriptscriptstyle\le}\,}{\,{\scriptscriptstyle\le}\,}}}
\def\Ge{{\mathchoice{\,{\scriptstyle\ge}\,}
  {\,{\scriptstyle\ge}\,}
  {\,{\scriptscriptstyle\ge}\,}{\,{\scriptscriptstyle\ge}\,}}}

\def\arrowsim{\,\smash{\mathop{\to}\limits^{\lower1.5pt
  \hbox{$\scriptstyle\sim$}}}\,}

\def\doublemaprights#1#2#3#4{\raise3pt\hbox{$\mathop{\,\,\hbox to     
#1pt{\rightarrowfill}\kern-30pt\lower3.95pt\hbox to
     #2pt{\rightarrowfill}\,\,}\limits_{#3}^{#4}$}}

\def\rightcapdownarrow{\raise9pt\hbox{$\ssize\cap$}\kern-7.75pt
     \Big\downarrow}

\def\rcapmapdown#1{\rightcapdownarrow\kern-1.0pt\vcenter{
     \hbox{$\scriptstyle#1$}}}

\def\rmapdown#1{\Big\downarrow\kern-1.0pt\vcenter{
     \hbox{$\scriptstyle#1$}}}
\def\rightsubsetarrow#1{{\ssize\subset}\kern-4.5pt\lower2.85pt
     \hbox to #1pt{\rightarrowfill}}
\def\longtwoheadedrightarrow#1{\raise2.2pt\hbox to #1pt{\hrulefill}
     \!\!\!\twoheadrightarrow}

\def\Gal{\operatorname{\hbox{Gal}}}
\def\Hom{\operatorname{\hbox{Hom}}}

\def\im{\hbox{Im}}

\NoBlackBoxes
\parindent=25pt
\document
\footline={\hfil}

\null
\noindent 
\centerline{\BIGtitle Good Reductions of Shimura Varieties of Preabelian}
\medskip
\centerline{\BIGtitle  Type in Arbitrary Unramified Mixed Characteristic, I}
\vskip 0.3 cm
\centerline{\bigtitle Adrian Vasiu, UA, 11/4/03}
\footline={\hfill}
\vskip 0.3 cm
\noindent
{\bf ABSTRACT.} We prove the existence of weak integral canonical models of Shimura varieties of Hodge type in arbitrary unramified mixed characteristic $(0,p)$. As a first application we solve a conjecture of Langlands for Shimura varieties of Hodge type. As a second application we prove the existence of integral canonical models of Shimura varieties of preabelian (resp. of abelian) type in mixed characteristic $(0,p)$ with $p\Ge 3$ (resp. with $p=2$) and with respect to hyperspecial subgroups; if $p=3$ (resp. if $p=2$) we restrict in this part I either to the $A_n$, $C_n$, $D_n^{\dbH}$ (resp. $A_n$ and $C_n$) types or to the $B_n$ and $D_n^{\dbR}$ (resp. $B_n$, $D_n^{\dbH}$ and $D_n^{\dbR}$) types which have compact factors over $\dbR$ (resp. which have compact factors over $\dbR$ in some $p$-compact sense). Though the second application is new just for $p\Le 3$, a great part of its proof is new even for $p\Ge 5$ and corrects [Va1, 6.4.11] in most of the cases. The second application forms progress towards the proof of a conjecture of Milne. It also provides in arbitrary mixed characteristic the very first examples of general nature of projective varieties over number fields which are not embeddable into abelian varieties and which have N\'eron models over certain local rings of rings of integers of number fields.
\medskip\noindent
{\bf Key words}: Shimura varieties, reductive group schemes, abelian schemes, integral models, p-divisible groups and $F$-crystals.
\medskip\noindent
{\bf MSC 2000}: Primary 11G10, 11G18, 14F30, 14G35, 14G40, 14K10 and 14J10.

\bigskip
\centerline{\bigtitle Contents}

{\nspace{

\bigskip
\line{\item{\bf\S1.}{Introduction}\leaderfill 1}

\smallskip
\line{\item{\bf \S2.}{Preliminaries}\leaderfill 8}

\smallskip
\line{\item{\bf \S3.}{On $\dbZ_{(p)}$ embeddings}\leaderfill 13}

\smallskip
\line{\item{\bf \S4.}{Crystalline applications} \leaderfill 22}

\smallskip
\line{\item{\bf \S5.}{Proof of 1.5}\leaderfill 33}

\smallskip
\line{\item{\bf \S6.}{Proof of 1.6 and complements}\leaderfill 41}

\smallskip
\line{\item{}{References}\leaderfill 54}

}}

\footline={\hss\tenrm \folio\hss}
\pageno=1

\bigskip
\noindent
{\boldsectionfont \S1. Introduction}
\bigskip

For a module $M$ over a commutative ring with unit $R$, let $M^*:=\Hom_R(M,R)$ and let $GL(M)$ be the group scheme over $R$ of linear automorphisms of $M$. A bilinear form on $M$ is called perfect if it induces an isomorphism $M\arrowsim M^*$. If $M$ is a free $R$-module of even rank and if $\psi_M$ is a perfect alternating form on $M$, then $Sp(M,\psi_M)$ and $GSp(M,\psi_M)$ are viewed as reductive group schemes. If $*_R$ or $*$ is an object or a morphism of the category of $R$-schemes, let $*_U$ be its pull back via an affine morphism $m_R:\Spec(U)\to\Spec(R)$. If there are several ways to take $m_R$, we mention the homomorphism $R\to U$ defining $m_R$ but often we do not add it as part of the notation $*_U$. 

If $F$ is a reductive group scheme over a scheme $Y$, then its fibres are connected and $Z(F)$, $F^{\der}$, $F^{\ad}$ and  $F^{\ab}$ denote the center, the derived group, the adjoint group and respectively the abelianization of $F$. So $Z^{\ab}(F)=F/F^{\der}$ and $F^{\ad}=F/Z(F)$. Let $Z^0(F)$ be the maximal subtorus of $Z(F)$. If $Y=\Spec(R)$ is affine, let $\Lie(F)$ be the $R$-Lie algebra of $F$ and let $\dbG_{mR}$ be the rank 1 split torus over $R$. Similarly, the group schemes $\dbG_{aR}$, $SL_{nR}$, etc., will be understood to be over $R$. If $F_1\hookrightarrow F$ is a monomorphism of group schemes over $R$ which is a closed embedding, then we identify $F_1$ with its image in $F$ and we consider intersections of subgroups of $F_1(R)$ with subgroups of $F(R)$. 

If $Y\to Z$ is an \'etale cover, then $\Res_{Y/Z} F$ is the reductive group scheme over $Z$ obtained from $F$ through the Weil restriction of scalars (see [BLR, 7.6] and [BrT, 1.5]). Let $\dbS:=\Res_{\dbC/\dbR} \dbG_{m\dbC}$. We identify $\dbS(\dbR)=\dbG_{m\dbC}(\dbC)$ and $\dbS(\dbC)=\dbG_{m\dbC}(\dbC)\times\dbG_{m\dbC}(\dbC)$ in such a way that the monomorphism $\dbR\hookrightarrow\dbC$ induces the map $z\to (z,\bar z)$, where $z\in\dbG_{m\dbC}(\dbC)$. 

Let $p$ be a rational prime. Let $\dbZ_{(p)}$ be the localization of $\dbZ$ with respect to its ideal $(p)$. If $E$ is a number field, let $E_{(p)}$ be the normalization of $\dbZ_{(p)}$ in $E$. For $v$ a finite prime of $E$, let $E_v$ be the completion of $E$ with respect to $v$ and let $e(v)$ be the index of ramification of $v$. Let $W(k)$ be the Witt ring of a perfect field $k$ of characteristic $p$. Always $n\in\dbN$. Let $\dbA_f:=\dbQ\otimes_{\dbZ}\widehat{\dbZ}$ be the topological ring of finite ad\`eles of $\dbQ$. Warning: let $\dbA_f^{(p)}$ be the ring of finite ad\`eles of $\dbQ$ with the $p$-component omitted; so $\dbA_f=\dbQ_p\times \dbA_f^{(p)}$. 

If $R\in\{\dbA_f,\dbA_f^{(p)},\dbQ_p\}$, then the group $F(R)$ is endowed with the coarsest topology making all maps $R=\dbG_{aR}(R)\to F(R)$ associated to morphisms $\dbG_{aR}\to F$ of $R$-schemes, to be continuous; so $F(R)$ is a totally discontinuous locally compact group. 

A continuous action of a totally discontinuous locally compact group on a scheme is always in the sense of [De2, 2.7.1] and is a right action. 

\bigskip\noindent
{\bf 1.1. A review of Shimura varieties.} A Shimura pair $(G,X)$ comprises from a reductive group $G$ over $\dbQ$ and a $G(\dbR)$-conjugacy class $X$ of homomorphisms $\dbS\to G_{\dbR}$ satisfying Deligne's axioms of [De2, 2.1]: the Hodge $\dbQ$--structure on $\Lie(G)$ defined by any $h\in X$ is of type $\{(-1,1),(0,0),(1,-1)\}$, $Ad\circ h(i)$ is a Cartan involution of $\Lie(G^{\ad}_{\dbR})$, and no simple factor of $G^{\ad}$ becomes compact over $\dbR$. Here $Ad:G_{\dbR}\to GL(\Lie(G^{\ad}_{\dbR}))$ is the adjoint representation. These axioms imply that $X$ has a natural structure of a hermitian symmetric domain, cf. [De2, 1.1.17]. For $h\in X$ we consider the Hodge cocharacter
$$\mu_h:\dbG_{m\dbC}\to G_{\dbC}$$ 
defined on complex points by the rule: $z\in\dbG_{m\dbC}(\dbC)$ is mapped into $h_{\dbC}(z,1)$. 

For generalities on Shimura pairs and varieties and their types we refer to [De2], [Mi2, \S 1], [Mi3] and [Va1, 2.2 to 2.8]. The most studied Shimura pairs are constructed as follows. Let $\psi$ be a non-degenerate symplectic form on an even dimensional vector space $W$ over $\dbQ$. Let $S$ be the set of all monomorphisms $\dbS\hookrightarrow GSp(W\otimes_{\dbQ} {\dbR},\psi)$ defining Hodge $\dbQ$--structures on $W$ of type $\{(-1,0),(0,-1)\}$ and having either $2\pi i\psi$ or $-2\pi i\psi$ as polarizations. The pair $(GSp(W,\psi),S)$ is a Shimura pair defining a Siegel modular variety. Let $L$ be a $\dbZ$-lattice of $W$ such that $\psi$ induces a perfect form $\psi:L\otimes_{\dbZ} L\to\dbZ$. Let $K_p:=GSp(L,\psi)(\dbZ_p)$.

The adjoint (resp. the toric part) of a Shimura pair $(G,X)$ is denoted as $(G^{\ad},X^{\ad})$ (resp. $(G^{\ab},X^{\ab})$), cf. [Va1, 2.4.1]. We say that $(G,X)$ is of Hodge type if there is an injective map $f:(G,X)\hookrightarrow (GSp(W,\psi),S)$ for some symplectic space $(W,\psi)$ over $\dbQ$ as above. So $f:G\hookrightarrow GSp(W,\psi)$ is a monomorphism such that $f_{\dbR}\circ h\in S$, $\forall h\in X$. 

We say that $(G,X)$ is of preabelian type if its adjoint $(G^{\ad},X^{\ad})$ is isomorphic to the adjoint of a Shimura pair $(G_1,X_1)$ of Hodge type. If moreover, we can choose $(G_1,X_1)$ such that $G_1^{\der}$ is an isogeny cover of $G^{\der}$, then we say that $(G,X)$ is of abelian type. 

The adjoint Shimura pair $(G^{\ad},X^{\ad})$ is called simple if $G^{\ad}$ is a simple $\dbQ$--group. In this case all simple factors of $G^{\ad}_{\dbC}$ have the same Lie type $LT$. If $LT$ is not $D_n$ for some $n\Ge 4$, then $(G^{\ad},X^{\ad})$ is said to be of $LT$ type. Let now $LT$ be $D_n$ with $n\Ge 4$. If $n\Ge 5$, then $(G^{\ad},X^{\ad})$ is said to be of  $D_n^{\dbH}$ (resp. $D_n^{\dbR}$) type iff all simple, non-compact factors of $G^{\ad}_{\dbR}$ are isomorphic to $SO^*(2n)_{\dbR}^{\ad}$ (resp. $SO(2,2n-2)^{\ad}_{\dbR}$); see [He, p. 445] for these classical groups. If $n\Ge 5$ and if $(G^{\ad},X^{\ad})$ is not of $D_n^{\dbH}$ or $D_n^{\dbR}$ type, then $(G^{\ad},X^{\ad})$ is said to be of $D_n^{\text{mixed}}$ mixed type. See [De2, p. 272] or [Va4, 2.2.1] for the difference between the $D_4^{\text{mixed}}$, $D_4^{\dbH}$ and $D_4^{\dbR}$ types. If $(G^{\ad},X^{\ad})$ is also of abelian type, then a classical result of Satake says that $LT$ is a classical Lie type (see [Sa]). Deligne proved that a simple, adjoint Shimura pair is of abelian type iff it is of $A_n$, $B_n$, $C_n$, $D_n^{\dbH}$ or $D_n^{\dbR}$ type (cf. [De2, 2.3.10]).

Let $E(G,X)\hookrightarrow\dbC$ be the number subfield of $\dbC$ which is the field of definition of the $G(\dbC)$-conjugacy class of the cocharacters $\mu_h$'s of $G_{\dbC}$, cf. [Mi2, p. 163--164]. We recall that $E(G,X)$ is called the reflex field of $(G,X)$. Let $O_{(v)}$ be the localization of the ring of integers $O(G,X)$ of $E(G,X)$ with respect to a prime $v$ of $E(G,X)$ dividing $p$. Let $k(v)$ be the residue field of $v$. 

The Shimura variety $\Sh(G,X)$ is identified with the canonical model over $E(G,X)$ of the complex Shimura variety 
$$\Sh(G,X)_{\dbC}:={\text{proj}.}{\text{lim}.}_{K\in CO(G)} G(\dbQ)\backslash X\times G(\dbA_f)/K,$$ 
where $CO(G)$ is the set of compact, open subgroups of $G(\dbA_f)$ endowed with the inclusion relation (see [De1], [De2] and [Mi1] to [Mi4]). So $\Sh(G,X)$ is an $E(G,X)$-scheme together with a continuous action of $G(\dbA_f)$. For $K$ a compact subgroup of $G(\dbA_f)$ let 
$$\Sh_K(G,X):=\Sh(G,X)/K.$$ 
If $K\in CO(G)$, then a classical result of Baily and Borel allows us to view $\Sh_K(G,X)_{\dbC}=G(\dbQ)\backslash X\times G(\dbA_f)/K$ as a complex variety and not just as a complex space (see [BB]). 

Let $H$ be a compact, open subgroup of $G_{\dbQ_p}(\dbQ_p)$. We refer to the quadruple $(G,X,H,v)$ as a Shimura quadruple. The group $G_{\dbQ_p}$ is called unramified iff it has a Borel subgroup and splits over an unramified, finite field extension of $\dbQ_p$. See [Ti2] or 2.2 for hyperspecial subgroups. Here we will add just two things. First $G_{\dbQ_p}(\dbQ_p)$ has hyperspecial subgroups iff $G_{\dbQ_p}$ is unramified. Second $H$ is a hyperspecial subgroup of $G_{\dbQ_p}(\dbQ_p)$ iff it is the group of $\dbZ_p$-valued points of a reductive group scheme over $\dbZ_p$ having $G_{\dbQ_p}$ as its generic fibre.

In this paper we study integral models of $\Sh_H(G,X)$ over different localizations of $E_{(p)}$ (like over $O_{(v)}$). The subject has a long history, the first fundamental result being the construction of (coarse) moduli schemes over $\dbZ$ of principally polarized abelian schemes endowed with level structures obtained by Mumford in 1965 (see [MFK, Ths. 7.9 and 7.10]). In 1976 Langlands conjectured the existence of a good integral model of $\Sh_H(G,X)$ over $O_{(v)}$, provided $H$ is a hyperspecial subgroup of $G_{\dbQ_p}(\dbQ_p)$ (see [La, p. 411]). But only in 1992, an idea of Milne made it significantly clearer how to characterize and identify the ``good" integral models. We now recall some definitions from [Va1] which are just slight modifications of the original ideas of Milne expressed in definitions [Mi2, 2.1, 2.2, 2.5 and 2.9].

\medskip\noindent
{\bf 1.1.1. Definitions.} Let $*$ be a faithfully flat $\dbZ_{(p)}$-algebra which is a localization of $E_{(p)}$.

\medskip
{\bf 1)} A smooth (resp. normal) integral model of $\Sh_H(G,X)$ over $*$ is a faithfully flat $*$-scheme $\scrN$ together with a continuous action of $G(\dbA_f^{(p)})$ on it such that the following two properties hold:

\medskip
{\bf (i)} {\it its generic fibre $\scrN_{E(G,X)}$ with its induced $G(\dbA_f^{(p)})$-action is $\Sh_H(G,X)$ with its canonical $G(\dbA_f^{(p)})$-action;}

\smallskip
{\bf (ii)} {\it there is a compact, open subgroup $H_0$ of $G(\dbA_f^{(p)})$ with the property that for any inclusion $H_1\leqslant H_2$ of open subgroups of $H_0$, the canonical morphism $\scrN/H_1\to\scrN/H_2$ induced by the action of $G(\dbA_f^{(p)})$ on $\scrN$, is an \'etale morphism between smooth (resp. normal), separated schemes of finite type over $*$.}

\medskip
The faithfully flat $*$-scheme $\scrN$ is separated (cf. (ii)) and so from (i) we get that the action of $G(\dbA_f^{(p)})$ on $\scrN$ is uniquely determined; so often we will not mention explicitly this action and we will refer just to $\scrN$ itself as a smooth (resp. normal) integral model of $\Sh_H(G,X)$ over $*$. 
Also if one (any) such $H_0$ is such that $\scrN/H_0$ is a smooth (resp. normal), (quasi-) projective $*$-scheme, then we say $\scrN$ is a pro-\'etale cover of a smooth (resp. normal), (quasi-) projective $*$-scheme.

\smallskip
{\bf 2)} Let $Y$ be a faithfully flat $\dbZ_{(p)}$-scheme which is regular. We say $Y$ is healthy regular if for any open subscheme $U$ of $Y$ containing $Y_{\dbQ}$ and whose complement in $Y$ is of codimension at least $2$, every abelian scheme over $U$ extends to an abelian scheme over $Y$. A flat $*$-scheme $Z$ is said to have the extension property if for any healthy regular scheme $\tilde Y$ over $*$, every morphism $\tilde Y_{E(G,X)}\to Z_{E(G,X)}$ extends uniquely to a morphism $\tilde Y\to Z$.

\smallskip
{\bf 3)} We say that a smooth integral model $\scrN$ of $\Sh_H(G,X)$ over $*$ is an integral canonical model of $\Sh_H(G,X)$ over $*$ if $\scrN$ (viewed just as a scheme) has the extension property and is a healthy regular scheme.

\medskip
For practical reasons we add the following definitions. 

\bigskip\noindent
{\bf 1.2. Definitions.} {\bf 1)} A smooth integral model $\scrN$ of $\Sh_H(G,X)$ over $O_{(v)}$ is called a weak integral canonical model of $\Sh_H(G,X)$ over $O_{(v)}$ if it has the following smooth extension property: for any regular, formally smooth  scheme $Y$ over a DVR $O$ faithfully flat over $O_{(v)}$ and of the same index of ramification as $O_{(v)}$, every morphism $Y_{E(G,X)}\to \scrN_{E(G,X)}=\Sh_H(G,X)$ of $E(G,X)$-schemes extends uniquely to a morphism $Y\to\scrN$ of $O_{(v)}$-schemes. 

\smallskip
{\bf 2)} We say $(G,X)$ has compact factors if the group $G^{\ad}_{\dbR}$ is non-trivial and has simple factors which are compact.

\smallskip
{\bf 3)} We assume that $e(v)=1$. An affine, flat group scheme $G_{\dbZ_{(p)}}$ over $\dbZ_{(p)}$ having $G$ as its generic fibre is called a quasi-reductive group scheme for $(G,X,v)$ if there is a reductive, normal subgroup $G^1_{\dbZ_{p}}$ of $G_{\dbZ_p}$ such that there is a cocharacter $\mu_v:\dbG_{mW(k(v))}\to G^1_{W(k(v))}$ whose extension to $\dbC$ via an (any) $O_{(v)}$-monomorphism $W(k(v))\hookrightarrow\dbC$ is $G(\dbC)$-conjugate with the cocharacters $\mu_h$ ($h\in X$) defined in the beginning of 1.1.

\smallskip
{\bf 4)} We say the prime $v$ is  $p$-compact for $(G,X)$ if either $G^{\ad}$ is trivial or the product decomposition $G^{\ad}_{\dbQ_p}=\prod_{i\in I_p} G_i$ in simple, adjoint $\dbQ_p$-groups is such that $\forall i\in I_p$, the pull back of $G_i$ to $\Spec(\dbC)$ via an (any) $E(G,X)$-embedding $E(G,X)_v\hookrightarrow\dbC$ has simple factors which are pull backs of simple, compact factors of $G^{\ad}_{\dbR}$.

\medskip
It is easy to see that any (weak) integral canonical model is uniquely determined up to unique isomorphism. Moreover, if $e(v)\Le\max\{1,p-2\}$, then any 
integral canonical model of $\Sh_H(G,X)$ over $O_{(v)}$ is also a weak integral canonical model. This is so due to the following theorem (see [Va1, 3.2.2 1)] for $p>2$ and see [Va2, 1.3] for the case when $e(v)=1$).

\medskip\noindent
{\bf 1.2.1. Theorem.} {\it If $e(v)\Le\max\{1,p-2\}$, then any regular, formally smooth scheme over a DVR which is a faithfully flat $O_{(v)}$-algebra and which is of index of ramification at most $\max\{1,p-2\}$, is healthy regular.}

\medskip
From 1.2.1 and [Va1, rm. 3.2.3.1 2)] we get that for $e(v)\Le\max\{1,p-2\}$ and for $*=O_{(v)}$, the definition 1.1.1 3) coincides with the definition [Va1, 3.2.3 6)]. So as in this paper we deal in essence just with the case $e(v)=1$ and with normal models which are separated, we will not make in what follows any distinction between the two possible definitions of integral canonical models.

\bigskip\noindent
{\bf 1.3. Basic notations.} 
Let $(GSp(W,\psi),S)$, $L$ and $K_p$ be as in 1.1. Let $L_{(p)}:=L\otimes_{\dbZ} \dbZ_{(p)}$. For the rest of this introduction we will assume that we have an injective map 
$$f\colon (G,X)\hookrightarrow (GSp(W,\psi),S)$$ 
and that 
$H:=K_p\cap G(\dbQ_p).$ 
Let $G_{\dbZ_{(p)}}$ be the Zariski closure of $G$ in $GL(L_{(p)})$. So $H=G_{\dbZ_{(p)}}(\dbZ_p)$. Let $E(G,X)$, $v$, $k(v)$ and $O_{(v)}$ be as in 1.1. Let 
$$\scrM$$ 
be the integral canonical model of $\Sh_{K_p}(GSp(W,\psi),S)$ over $\dbZ_{(p)}$. It is known that $\scrM$ is the moduli scheme over $\dbZ_{(p)}$ parameterizing isomorphism classes of principally polarized abelian schemes of relative dimension $\dim_{\dbQ}(W)/2$ and having level $N$ symplectic similitude structure, $\forall N\in\dbN$ prime to $p$ (cf. [Va1, 3.2.9 and 4.1] and [MFK, Ths. 7.9 and 7.10]). See [Va1, 3.2.9 and 4.1] and [De1, 4.21] for the natural action of $GSp(W,\psi)(\dbA_f^{(p)})$ on $\scrM$. It is known that $\Sh_H(G,X)$ is a finite scheme over $\scrM_{E(G,X)}=\Sh_{K_p}(GSp(W,\psi),S)_{E(G,X)}$ (for instance, cf. 3.1.1 3)). If $G_{\dbZ_{(p)}}$ is a reductive group scheme, then in fact $\Sh_H(G,X)$ is a closed subscheme of $\scrM_{E(G,X)}$ (cf. [Va1, 3.2.14]). 

Let 
$$\scrN^\prime$$ 
be the normalization of $\scrM_{O_{(v)}}$ in the ring of fractions of $\Sh_H(G,X)$; so if $G_{\dbZ_{(p)}}$ is a reductive group scheme, then $\scrN^\prime$ is the normalization of the Zariski closure of $\Sh_H(G,X)$ in $\scrM_{O_{(v)}}$. As $G(\dbA_f^{(p)})$ acts on $\Sh_H(G,X)$ and $\scrM$, we get a natural action of $G(\dbA_f^{(p)})$ on $\scrN^\prime$. We consider the formally smooth locus 
$$\scrN$$ 
of $\scrN^\prime$ over $O_{(v)}$. It is a $G(\dbA_f^{(p)})$-invariant, open subscheme of $\scrN^\prime$ such that $\scrN_{E(G,X)}=\scrN^\prime_{E(G,X)}$, cf. 3.1.2. Let 
$$(\scrA^\prime,\scrP_{\scrA^\prime})$$ 
be the pull back to $\scrN^\prime$ of the universal abelian scheme over $\scrM$. 

The main goal of this paper is to prove the following three basic Theorems. 

\bigskip\noindent
{\bf 1.4. Basic Theorem.}
{\it We assume that $e(v)=1$ (i.e. that $v$ is unramified over $p$) and that the $k(v)$-scheme $\scrN_{k(v)}$ is non-empty. We have:

\medskip
{\bf 1)} The $O_{(v)}$-scheme $\scrN$ together with the natural action of $G(\dbA_f^{(p)})$ on it is a weak integral canonical model of $\Sh_H(G,X)$ over $O_{(v)}$. 

\smallskip
{\bf 2)} For any algebraically closed field $k$ of characteristic $p$, the natural morphism $\scrN_{W(k)}\to\scrM_{W(k)}$ is a formally closed embedding at every $k$-valued point of $\scrN_{W(k)}$.} 

\bigskip\noindent
{\bf 1.5. Basic Theorem.} {\it We assume that $e(v)=1$ and that $G_{\dbZ_{(p)}}$ is a quasi-reductive group scheme for $(G,X,v)$. If $p=2$ we also assume that any pull back of $\scrA^\prime$ via a geometric point of $\scrN_{k(v)}$ has $p$-rank $0$. Then $\scrN_{k(v)}$ is a non-empty, open closed subscheme of $\scrN^{\prime}_{k(v)}$.}

\bigskip\noindent
{\bf 1.6. Main Theorem.} {\it {\bf 1)} Let $(G^{\ad},X^{\ad})=\prod_{i\in I_0} (G_1^i,X_1^i)$ be the product decomposition 
in simple, adjoint  Shimura pairs.
We assume that $G_{\dbZ_{(p)}}$ is a reductive group scheme and that one of the following three conditions holds:

\medskip
{\bf (i)} $p\Ge 5$;

\smallskip
{\bf (ii)} $p=3$ and $\forall i\in I_0$, the simple, adjoint Shimura pair $(G_1^i,X_1^i)$ either is of $A_n$, $C_n$ or $D_n^{\dbH}$ type or is of $B_n$ or $D_n^{\dbR}$ type and moreover has compact factors;

\smallskip
{\bf (iii)} $p=2$ and $\forall i\in I_0$, the simple, adjoint Shimura pair $(G_1^i,X_1^i)$ either is of $A_n$ or $C_n$ type or is of $B_n$, $D_n^{\dbH}$ or $D_n^{\dbR}$ type and moreover the prime $v_1^i$ of $E(G_1^i,X_1^i)$ divided by $v$ is $p$-compact for $(G_1^i,X_1^i)$. 

\medskip
Then $e(v)=1$ and $\scrN=\scrN^\prime$ is an integral canonical model of $\Sh_H(G,X)$ over $O_{(v)}$. 

\medskip
{\bf 2)} Let $(G_1,X_1)$ be a Shimura pair of preabelian type. Let $p$ be a prime such that $G_{1\dbQ_p}$ is unramified; so $E(G_1,X_1)$ is unramified over $p$ (cf. [Mi3, Cor. 4.7 (a)]). Let $H_1$ be a hyperspecial subgroup of $G_{1\dbQ_p}(\dbQ_p)$. Let $(G_1^{\ad},X_1^{\ad})=\prod_{i\in I_0} (G_1^i,X_1^i)$ be the product decomposition 
in simple, adjoint  Shimura pairs. We consider the set $S(G_1,X_1,2)$ of primes $v_1$ of $E(G_1,X_1)$ dividing $2$ and such that $\forall i\in I_0$ with the property that $(G_1^i,X_1^i)$ is of $B_n$, $D_n^{\dbH}$ or $D_n^{\dbR}$ type, the prime $v_1^i$ of $E(G_1^i,X_1^i)$ divided by $v$ is $2$-compact for $(G_1^i,X_1^i)$. Let $O(G_1,X_1,2)$ be the localization of $E(G_1,X_1)_{(2)}$ whose finite primes are precisely the primes in $S(G_1,X_1,2)$. If $p\Ge 3$ let $O(G_1,X_1,p):=E(G_1,X_1)_{(p)}$.

If $p=3$ we assume that condition (ii) of 1) holds. If $p=2$ we assume that $(G_1,X_1)$ is of abelian type and that the set $S(G_1,X_1,2)$ is non-empty.  Then there is an integral canonical model $\scrN_1$ of $\Sh_{H_1}(G_1,X_1)$ over $O(G_1,X_1,p)$ which is a pro-\'etale cover of a smooth, quasi-projective $O(G_1,X_1,p)$-scheme. If moreover $(G_1^i,X_1^i)$ has compact factors for all $i\in I_0$, then $\scrN_1$ is a pro-\'etale cover of a smooth, projective $O(G_1,X_1,p)$-scheme.}

\bigskip\noindent
{\bf 1.7. On contents.} \rm
We detail on the contents of the paper. In 2.1 we list conventions and notations. In 2.2 we include some results on extending reductive group schemes; the main new result (see 2.2.3) is a purity type of theorem for the classical Lie types. In 2.3 we review properties of dilatations of affine group schemes over Witt rings. In \S3 we list some properties of embeddings between Shimura pairs  in a $\dbZ_{(p)}$ context. 

In \S4 we list different crystalline applications. In 4.1 and 4.2 we introduce basic notations and review two recent results pertaining to $p$-divisible groups and which play a central role in the rest of the paper. The results are: the extension theorem of de Jong (see [dJ2]) and a conjecture of Milne proved in [Va6, 1.2]. In 4.3 we prove 1.4. In 4.4 we present a simple criterion on when the $k(v)$-scheme $\scrN_{k(v)}$ is non-empty. In 4.5 we apply 1.4 1) and 4.4 to get a solution of the mentioned conjecture of Langlands for the case of Shimura varieties of Hodge type. In 4.6 we refine  the deformation theories of [Fa, \S7] and [Va1, 5.4] as allowed by the recent tools developed in [Va6].

In \S5 we use 2.2 and 4.1 to 4.3 to prove 1.5. The proof of 1.6 is carried on in 6.1 and 6.2. In 6.3 we use the projectiveness part of 1.6 2) to provide new examples of general nature of N\'eron models in the sense of [BLR, p. 12]. In 6.4 we include some complements; in particular in 6.4 5) and 6) we explain why we can view 1.5 as a generalization of [Mo] for the unramified context with $p>2$. 

\bigskip\noindent
{\bf 1.8. More on literature.} 
Referring to 1.4 1), in [No] it is shown that all ordinary points of $\scrN^\prime_{k(v)}$ belong to $\scrN_{k(v)}$. Except the part involving projectiveness, Theorem 1.6 was proved for $p\Ge 5$ in [Va1, 3.2.12, 5.1 and 6.4.1] (see also [Va5, 7.9.3] for the passage from the abelian type case to the preabelian type case) and for unitary Shimura varieties in [Va4, \S5]. The works [MFK], [Dr], [Mo], [Zi], [LR], [Ko], [Va1], [Va3] to [Va5] are the most relevant literature on the existence of integral canonical models. See also [HT, \S5] for a more recent translation of part of [Dr] in terms of unitary Shimura varieties. Theorems 1.4 and 1.5 can be also viewed as steps in proving Conjectures B 3.7 and B 3.12 of [Re]. 

In part II of this paper we will use [Va1, 6.6.2] and the formalism of smooth toroidal compactifications and elementary inductions to show that 1.6 1) always holds for $p\Le 3$ provided $G_{\dbZ_{(p)}}$ is a reductive group scheme and that 1.6 2) always holds with $O(G_1,H_1,p)$ being replacement by $E(G_1,X_1)_{(p)}$ even for $p=2$ and under no extra restriction for $p\Le 3$. This will complete the proof of Milne's conjecture of [Mi2, 2.7] and [Va1, 3.2.5] for Shimura varieties of preabelian type in mixed characteristics $(0,p)$ with $p\Ge 3$ and of abelian type in mixed characteristic $(0,2)$. 

This paper brings several completely new ideas in order to essentialize, shorten and simplify [Va1] and to extend many parts of loc. cit. worked out just for $p\Ge 5$ to the case of small primes $p\in\{2,3\}$. Moreover, the projectiveness part of 1.6 corrects for $p\Ge 5$ in most of the cases an error in the proof of [Va1, 3.2.3.2 ii)] which invalidated [Va1, 6.4.1.1 i) and most of 6.4.11]. Such a correction was started in [Va4, 5.2.1] and is significantly continued here; it will be implicitly completed in part II of the paper. We would like to thank Universities of Utah and Arizona for providing us with good conditions for the writing of this paper. This research was partially supported by the NSF grant DMF 97-05376.

\bigskip
\noindent
{\boldsectionfont \S2. Preliminaries}
\bigskip 

In 2.1 we include some conventions and notations. In 2.2 we first prove two extension results on reductive group schemes whose adjoints have classical Lie type and then we get a more general version of [Va1, 3.1.6]. In 2.3 we recall few basic properties of dilatations of affine group schemes over Witt rings.

\bigskip\noindent
{\bf 2.1. Conventions and notations.} 
Let $\dbF$ be the algebraic closure of the field with $p$ elements. Let $k$ be a perfect field of characteristic $p$. Let $\sigma:=\sigma_k$ be the Frobenius automorphism of $k$, $W(k)$ and $B(k):=W(k)\fracwithdelims[]1p$. Let $F^{\sc}$ be the simply connected semisimple group cover of the derived group $F^{\der}$ of a reductive group scheme $F$. Warning: we will use the notations of 1.3 only when this is explicitly stated so; however, the standard conventions and notations of the paragraphs before 1.1 as well as of 1.1 will be used everywhere. 

If $M$ and $R$ are as in the beginning of \S1, then let
$$\scrT(M):=\oplus_{s,t\in\dbN\cup\{0\}} M^{\otimes s}\otimes_R M^{*\otimes t}.$$
We identify naturally $\End(M)=M\otimes_R M^*$ and $\End(\End(M))=M^{\otimes 2}\otimes_R M^{*\otimes 2}$. Let $x_R\in R$ be a non-divisor of $0$. A family of tensors of $\scrT(M[{1\over x_R}])=\scrT(M)[{1\over x_R}]$ is denoted $(w_{\alpha})_{\alpha\in\scrJ}$, with $\scrJ$ as the set of indices. Let $M_1$ be another $R$-module. Let $(w_{1\alpha})_{\alpha\in\scrJ}$ be a family of tensors of $\scrT(M_1[{1\over x_R}])$ indexed by the same set $\scrJ$. By an isomorphism $f:(M,(w_{\alpha})_{\alpha\in\scrJ})\arrowsim (M_1,(w_{1\alpha})_{\alpha\in\scrJ})$ we mean an isomorphism $f:M\arrowsim M_1$ giving birth to an isomorphism $\scrT(M[{1\over x_R}])\arrowsim\scrT(M_1[{1\over x_R}])$ taking $w_{\alpha}$ into $w_{1\alpha}$, $\forall\alpha\in\scrJ$. If $R$ is local, let $R^{\sh}$ be its strict henselization. 

We denote two tensors or two bilinear forms in the same way provided they are obtained one from another via a reduction or a scalar extension.

\bigskip\noindent
{\bf 2.2. Extending reductive group schemes.} 
Let $Y$ be a reduced scheme such that the ring of fractions $\eta$ of $Y$ is a finite product of fields. Let $\scrF$ be a locally free sheaf on $Y$ of finite rank. Let $GL(\scrF)$ be the reductive group scheme over $Y$ of linear automorphism of $\scrF$. Let $U$ be an open subscheme of $Y$. Let $C:=Y\setminus U$ be endowed with the reduced scheme structure. Let $E_U$ be a reductive group scheme over $U$. We next recall the following result [SGA3, Vol. III, Prop. 6.1 of p. 32]: 

\medskip\noindent
{\bf 2.2.1. Proposition.} {\it Locally in the \'etale topology of $U$, $E_U$ has maximal split tori.} 

\medskip
We refer to [Va1, 3.1.2.1 c) and 3.1.2.2 3)] for the proof of the next Proposition.

\medskip\noindent
{\bf 2.2.2. Proposition.} {\it We assume that $R$ is a DVR. Let $q:E_1\to E_2$ be a homomorphism between flat, affine groups over $R$ such that $E_1$ is a reductive group scheme and $q_{\eta}$ is a closed embedding. Then $q$ is a closed embedding.}

\medskip\noindent
{\bf 2.2.3. Theorem.} {\it We assume $Y$ is a regular scheme whose local rings have dimension at most 2. We also assume that $C$ is of codimension 2 in $Y$ and that all simple factors of pull backs of $E_U^{\ad}$ via geometric points of $U$ are of classical Lie type. Then $E_U$ extends uniquely to a reductive group scheme $E$ over $Y$.}

\medskip
\proof
It suffices to prove the Theorem under the extra assumption that $Y$ is a local scheme $\Spec(R)$ and $C=\Spec(k(R))$, where $k(R)$ is the residue field of $R$. So $\eta$ is a field. Let $R_E$ be the $\dbZ$-algebra of global functions of $E_U$. If $E$ exists, then $E$ as a scheme is $\Spec(R_E)$ and moreover the identity section, the multiplication operation and the inverse operation of $E$ are uniquely determined by the corresponding analogues of $E_U$. This takes care of the uniqueness part. We are left to show that the morphism $\Spec(R_E)\to Y$ is smooth and defines naturally a reductive group scheme over $Y$ extending $E_U$. It suffices to prove this under the extra assumption that $R$ is strictly henselian.

We show that $E^{\ab}_U$ extends to a torus over $Y$. Let $\eta_1$ be the smallest Galois extension of $\eta$ such that $E^{\ab}_{\eta_1}$ is a split torus. Let $R_1$ be the normalization of $R$ in $\eta_1$. From 2.2.1 we get that $\Spec(R_1)$ has \'etale points above any point of $Y$ of codimension 1. So as $\eta_1$ is a Galois extension of $\eta$, the morphism $\Spec(R_1)\to Y$ is an \'etale cover above points of $Y$ of codimension at most 1 and so from the classical purity theorem (see [SGA1, p. 275]) we get that it is an \'etale cover. So as $R$ is strictly henselian, we have $R_1=R$. So $E^{\ab}_{\eta}$ is a split torus and so also $E^{\ab}_U$ is a split torus. So $E^{\ab}_U$ extends to a split torus $E^{\ab}$ over $Y$.

We show that $E^{\ad}_U$ extends to an adjoint group scheme over $U$; the argument for this extends until the first paragraph after Case 4 below. We write $E^{\ad}_{\eta}=\prod_{i\in I} E_{i\eta}$ as a product of simple, adjoint groups over $\eta$ (cf. [Ti1, 3.1.2]). Let $\eta_i$ be the finite field extension of $\eta$ such that $E_{i\eta}=\Res_{\eta_i/\eta} E^i_{\eta_i}$, with $E^i_{\eta_i}$ as an absolutely simple, adjoint group over $\eta_i$ (cf. loc. cit.). Let $R_i$ be the normalization of $R$ in $\eta_i$. Let $V$ be a local ring of $R$ which is a DVR. Let $V_0$ be an integral, faithfully flat, \'etale $V$-algebra such that the group $E^{\ad}_{V_0}$ is split (cf. 2.2.1). Let $\eta_0$ be the field of fractions of $V_0$. As the group $E^{\ad}_{\eta_0}=\prod_{i\in I} \Res_{\eta_i\otimes_{\eta} \eta_0/\eta_0} E^i_{\eta_i\otimes_{\eta} \eta_0}$ is split and adjoint, we get that for each $i\in I$ the $\eta_0$-algebra $\eta_i\otimes_{\eta} \eta_0$ is isomorphic to a product of a finite number of copies of $\eta_0$. This implies that $R_i\otimes_R V_0$ is an \'etale $V_0$-algebra, $\forall i\in I$. So $R_i\otimes_R V$ is an \'etale $V$-algebra, $\forall i\in I$. So as in the previous paragraph we argue that $R_i=R$, $\forall i\in I$. So $E_{i\eta}$ is an absolutely simple, adjoint group, $\forall i\in I$. Let $E_{iU}$ be the Zariski closure of $E_{i\eta}$ in $E^{\ad}_U$. We have a direct product decomposition $E_U=\prod_{i\in I} E_{iU}$ as one can easily check this locally in the \'etale topology of $U$.

Let $S_{iU}$ be the split, simple, adjoint group scheme over $U$ of the same Lie type $LT_i$ as any geometric fibre of $E_{iU}$. We have a short exact sequence $0\to S_{iU}\to \text{Aut}(S_{iU})\to M_{iU}\to 0$, where $M_{iU}$ is a finite, \'etale  group scheme over $U$ (cf. [SGA3, Vol. III, Th. 1.3 of p. 328]). It is well known that $M_{iU}$ is $\dbZ/2\dbZ$ if $LT_i\in\{A_n|n\Ge 2\}\cup\{D_n|n\Ge 5\}$, is the symmetric group $S_3$ if $LT_i=D_4$ and is trivial otherwise. Let $\gamma_i\in H^1(\eta,\text{Aut}(S_{i\eta}))$ be the class defining $E_{i\eta}$ and let $\delta_i$ be its image in $H^1(\eta,M_{i\eta})$. The non-trivial torsors of $\dbZ/2\dbZ$ and $\dbZ/3\dbZ$ over $\eta$ are in one-to-one correspondence to quadratic and respectively cubic Galois extensions of $\eta$. Let $\eta^i$ be the smallest Galois extension of $\eta$ over which $\delta_i$ becomes the trivial class. As in the previous two paragraphs we easily get that the normalization of $R$ in $\eta^i$ is $R$ itself. So $\delta_i$ is the trivial class. So $E_{i\eta}$ is an inner form of $S_{i\eta}$. We now show that $E_{iU}$ extends to an adjoint group scheme over $Y$. We consider four Cases.

\smallskip
{\bf Case 1: $LT_i=A_n$}. As $E_{i\eta}$ is an inner form of $S_{i\eta}$, $Z(E_{i\eta}^{\sc})$ is the multiplicative group $\mu_{n+1\eta}$ over $\eta$. So $Z(E_{iU}^{\sc})$ is also the multiplicative group scheme $\mu_{n+1U}$ over $U$. Let $G_{iU}$ be the quotient of $\dbG_{mU}\times_U E_{iU}^{\sc}$ by the $\mu_{n+1U}$ subgroup embedded via the natural embedding $\mu_{n+1U}\hookrightarrow\dbG_{mU}$ and via the composite of the inverse automorphism $\mu_{n+1U}\arrowsim\mu_{n+1U}$ with the natural embedding $\mu_{n+1U}\hookrightarrow E_{iU}^{\sc}$. Let $\gra_i$ be the locally free $R$-module of finite rank of global sections of the Lie sheaf $\Lie(G_{iU})$ on $U$, cf. [FC, Lemma 6.2 of p. 181]. The group $G_{i\eta}$ is an inner form of $GL_{n+1\eta}$ and so is the group of invertible elements of a semisimple $\eta$-algebra. So $\gra_i\otimes_R \eta$ has a canonical structure of a semisimple $\eta$-algebra. Let $Tr_i$ be the canonical trace form on $\gra_i\otimes_R \eta$. As $G_{iU}$ locally in the \'etale topology of $U$ is isomorphic to $GL_{n+1U}$, the $\eta$-algebra structure and the trace form on $\gra_i\otimes_R \eta$ extend uniquely to an $R$-algebra structure and respectively to a trace form $Tr_i$ on $\gra_i$ which is perfect in codimension at most 1. So the trace form $Tr_i$ on $\gra_i$ is as well perfect. 

Let $G_i$ be the group scheme over $Y$ of invertible elements of the $R$-algebra $\gra_i$. As a scheme, $G_i$ is an affine, open subscheme of the vector group scheme over $Y$ defined by the $R$-module $\gra_i$. So $G_i$ is smooth, has connected fibres and it is easy to see that $\Lie(G_i)$ is $\gra_i$ endowed with the Lie bracket defined by the $R$-algebra structure of $\gra_i$. Moreover, the notations match, i.e. the restriction of $G_i$ to $U$ is the group scheme $G_{iU}$ considered in the previous paragraph. Argument: the canonical identification of $G_{i\eta}$ with the fibre of $G_i$ over $\eta$ extends to a canonical identification of $G_{iU}$ with the restriction of $G_i$ to $U$, as one can easily check this locally in the \'etale topology of $U$. 

We now check that $G_i$ is a reductive group scheme over $Y$. To check this it suffices to show that the fibre of $G_i$ over the maximal point of $Y$ is a  reductive group. To check this last thing it suffices to show that $\gra_i\otimes_R k(R)$ is a semisimple $k(R)$-algebra. Let $\tilde V$ be an $R$-algebra which is a complete DVR having an algebraically closed residue field and such that $G_{i\tilde\eta}$ splits over the field of fractions $\tilde\eta$ of $\tilde V$ and the natural morphism $\Spec(\tilde V)\to Y$ has as image two distinct points, one of them being the maximal point of $Y$. We consider a homomorphism $G_{i\tilde V}\to GL(\tilde V^{n+1})$ over $\tilde V$ which over $\tilde\eta$ is an isomorphism, cf. [Ja, 10.4 of Part 1]. As the trace form $Tr_i$ is perfect, the Lie homomorphism $\gra_i\otimes_R \tilde V=\Lie(G_{i\tilde V})\to \Lie(GL(\tilde V^{n+1}))$ is an isomorphism. So the $\tilde V$-homomorphism $\gra_i\otimes_R\tilde V\to \End(\tilde V^{n+1})$ is as well an isomorphism. So $\gra_i\otimes_R\tilde V$ is a semisimple $\tilde V$-algebra. So $\gra_i\otimes_R k(R)$ is a semisimple $k(R)$-algebra. So $G_i$ is a reductive group scheme. 

So $G_i^{\ad}$ is an adjoint group scheme over $Y$ whose restriction to $U$ is $G_{iU}^{\ad}=E_{iU}$. 

\smallskip
{\bf Case 2: $LT_i=B_n$.} We consider the faithful representation $W_i$ of $E_{i\eta}$ of dimension $2n+1$. We check that there is a free $R$-submodule $L_i$ of $W_i$ of rank $2n+1$ such that the natural monomorphism $E_{i\eta}\hookrightarrow GL(W_i)$ extends to a homomorphism $E_{iU}\to GL(L_i)_U$. If $V$ is a local ring of $R$ which is a DVR, then from [Ja, 10.4 of Part 1] we get the existence of a $V$-lattice $L_i(V)$ of $W_i$ such that the monomorphism $E_{i\eta}\hookrightarrow GL(W_i)$ extends to a homomorphism $E_{iV}\to GL(L_i(V))$. Let $\tilde V$ be a faithfully flat $V$-algebra which is a complete DVR having an algebraically closed residue field. So $E_{i\tilde V}$ is split and so there is a $\tilde V$-lattice $\tilde L_i(V)$ of $W_i\otimes_V \tilde V$ normalized by $E_{i\tilde V}$ and such that the special fibre of the homomorphism $E_{i\tilde V}\to GL(\tilde L_i(V))$ is an irreducible, faithful representation. So from [Ja, 10.9 of Part I] we get that the special fibre of the representation of $E_{iV}$ on $L_i(V)$ is absolutely irreducible. So $L_i(V)$ is uniquely determined up to a $\dbG_{mR}(\eta)$-multiple. 

Let $(\Spec(R[{1\over {f_j}}]))_{j\in J}$ be a finite, affine, open cover of $U$ such that there is a free $R[{1\over {f_j}}]$-submodule $L_i(f_j)$ of $W_i$ with the property that $\forall j\in J$ the monomorphism $E_{i\eta}\hookrightarrow GL(W_i)$ extends to a homomorphism $E_{iR[{1\over {f_j}}]}\to GL(L_i(f_j))$ over $\Spec(R[{1\over {f_j}}])$. If $j_1$, $j_2\in J$, then $L_i(f_{j_1})[{1\over {f_{j_2}}}]$ and $L_i(f_{j_2})[{1\over {f_{j_1}}}]$ differ just by a $\dbG_{mR}(R[{1\over {f_{j_1}f_{j_2}}}])$-multiple. This is so as each $L_i(V)$ is uniquely determined up to a $\dbG_{mR}(\eta)$-multiple and as $R[{1\over {f_{j_1}f_{j_2}}}]$ is a unique factorization domain. So as $R$ is itself a unique factorization domain, there are non-zero elements $s_j\in\eta$ ($j\in J$) such that $s_{j_1}L_i(f_{j_1})[{1\over {f_{j_2}}}]=s_{j_2}L_i(f_{j_2})[{1\over {f_{j_1}}}]$, $\forall j_1$, $j_2\in J$. So the free sheaves on $\Spec(R[{1\over {f_j}}])$ defined by $s_{j}L_i(f_{j})$'s, glue together to define a locally free sheaf on $U$. It extends to a free sheaf on $\Spec(R)$ (cf. [FC, Lemma 6.2 of p. 181]) and the global sections of it are the searched for free $R$-module $L_i$. So the monomorphism $E_{i\eta}\hookrightarrow GL(W_i)$ extends to a homomorphism $E_{iU}\to GL(L_i)_U$ which is a closed embedding, cf. 2.2.2.

Let $\eta_{\text{sep}}$ be the separable closure of $\eta$. As $E_{i\eta_{\text{sep}}}$ is split, there is a non-degenerate quadratic form $F_{i\eta_{\text{sep}}}$ on $W_i\otimes_{\eta} \eta_{\text{sep}}$ such that $E_{i\eta_{\text{sep}}}$ is the subgroup of $SL(W_i\otimes_{\eta} \eta_{\text{sep}})$ preserving $F_{i\eta_{\text{sep}}}$. Argument: if $\eta$ is of characteristic $2$, then this follows from [Bo, 23.5 and 23.6]; the case when $\eta$ is not of characteristic $2$ is well known. The quadratic form $F_{i\eta_{\text{sep}}}$ is uniquely determined up to a $\dbG_{mR}(\eta_{\text{sep}})$-multiple. Argument: any non-degenerate quadratic form on $W_i\otimes_{\eta} \eta_{\text{sep}}$ is isomorphic to the quadratic form $x_1x_2+x_3x_4+...+x_{2n-1}x_{2n}+x_{2n+1}^2$ on $W_i\otimes_{\eta} \eta_{\text{sep}}$ corresponding to a fixed ordered basis of $W_i\otimes_{\eta} \eta_{\text{sep}}$ (see loc. cit. for the case when $\eta$ is of characteristic $2$); so the statement follows from the fact that the normalizer of $E_{i\eta_{\text{sep}}}$ in $GL(W_i\otimes_{\eta} \eta_{\text{sep}})$ is generated by $E_{i\eta_{\text{sep}}}$ and by $Z(GL(W_i\otimes_{\eta} \eta_{\text{sep}}))$. 

So from Hilbert's Theorem 90 we get that a non-zero scalar multiple of $F_{i\eta_{\text{sep}}}$ is definable over $\eta$, i.e. is the scalar extension of a quadratic form $F_{i\eta}$ on $W_i$. For any DVR $V$ as above, a $\dbG_{mR}(\eta)$-multiple of $F_{i\eta}$ extends to a quadratic form $F_{iV}$ on $L_i\otimes_R V$. Argument: to check this one can work locally in the \'etale topology of $\Spec(V)$ and so it suffices to consider the case when $E_{iV}$ is split; but this case is trivial.

As $R$ is a unique factorization domain, we deduce the existence of a quadratic form $F_i$ on $L_i$ which over each DVR $V$ as above is a $\dbG_{mR}(V)$-multiple of $F_{iV}$. The quadratic form $F_i$ on $L_i$ is perfect as it is perfect in codimension at most 1. So the subgroup of $SL(L_i)$ fixing  $F_i$ is an adjoint group scheme over $Y$ extending $E_{iU}$.

\smallskip
{\bf Case 3: $LT_i=C_n$.} This case is entirely the same as Case 2. We just have to replace quadratic form by alternating form and the reference to [Bo, 23.5 and 23.6] by [Bo, 23.3]. 

\smallskip 
{\bf Case 4: $LT_i=D_n$.} We denote also by $\gamma_i$ the class of $H^1(\eta,E_{i\eta})$ whose image in $H^1(\eta,\text{Aut}(E_{i\eta}))$ is $\gamma_i$. Let $\tilde S_{iU}$ be the central, isogeny cover of $S_{iU}$ which is a split $SO$ group. We have a central short exact sequence $0\to\mu_2\to \tilde S_{i\eta}\to S_{i\eta}\to 0$. Let $\gamma_i(2)\in H^2(\eta,\mu_2)$ be the coboundary of $\gamma_i\in H^1(\eta,S_{i\eta})$. Let $G_{i\eta}$ be the group scheme of invertible elements of the central division algebra $D_{i\eta}$ over $\eta$ defining $\gamma_i(2)$. As $E_{iU}$ is an adjoint group scheme, we easily get that $D_{i\eta}$ extends naturally to a semisimple algebra over any local ring of $U$ which is a DVR. So $G_{i\eta}$ extends naturally to a reductive group scheme $G_{iU}$ over $U$. As in Case 1 we argue that $G_{iU}$ extends to a reductive group scheme $G_i$ over $Y$. As $Y$ is strictly henselian, $G_i$ is split (cf. 2.2.1). So the $\eta$-algebra $D_{i\eta}$ is a matrix $\eta$-algebra and so it is $\eta$ itself. So $\gamma_i(2)$ is the identity element of $H^2(\eta,\mu_2)$. So as in Case 2 we can consider the faithful representation $W_i$ of $E_{i\eta}$ of dimension $2n$. The rest is as in Case 2 (with $2n+1$ and the quadratic form $x_1x_2+x_3x_4+...+x_{2n-1}x_{2n}+x_{2n+1}^2$ being replaced by $2n$ and respectively by the quadratic form $x_1x_2+x_3x_4+...+x_{2n-1}x_{2n}$). So $E_{iU}$ extends to an adjoint group scheme over $Y$.

\medskip
So each $E_{iU}$ extends to an adjoint group scheme over $Y$. So $E_U^{\ad}=\prod_{i\in I} E_{iU}$ extends to an adjoint group scheme $E^{\ad}$ over $Y$. 

As $Y$ is strictly henselian, $E^{\ad}$ is split (cf. 2.2.1). As $E_U$ is isogenous to $E^{\ad}_U\times_U E^{\ab}_U$, the group scheme $E_U$ is split. So the fact that $E_U$ extends to a split, reductive group scheme $E$ over $Y$ follows from the uniqueness of a split reductive group scheme associated to a root datum, see [SGA3, Vol. III, Cor. 5.2 of p. 314]. We now include a second (more direct and elementary) way to check the existence of $E$. 

Let $E^1_{U}:=E_U\times_{E^{\ad}_U\times_U E^{\ab}_U} E^{\sc}_U\times_U E^{\ab}_U$. So $E^1_{U}$ is an affine group scheme over $U$ and $E_U^{\sc}$ is naturally a normal, closed subgroup of it. The quotient group $E_U^1/E_U^{\sc}$ is a group scheme of multiplicative type and so the extension of a finite, flat group scheme of finite type by a torus $E_U^T$. Let $E_U^2$ be the normal, flat, closed subgroup of $E_U^1$ containing $E_U^{\sc}$ and such that the quotient group $E_U^2/E_U^{\sc}$ is naturally identified with $E_U^T$. The group scheme $E_U^2$ is a reductive group scheme isomorphic to $E^{\sc}_U\times_U E^{\ab}_U$ and so it extends to a reductive group scheme $E^2$ over $Y$ isomorphic to $E^{\sc}\times_Y E^{\ab}$. The kernel of the natural central epimorphism $E^2_{U}\to E_U$ extends to a finite, flat, closed subgroup $Z^2$ of $Z(E^2)$ of multiplicative type. The quotient group scheme $E^2/Z^2$ is a reductive group scheme $E$ over $Y$ extending $E_U$. This ends the proof.  

\medskip\noindent
{\bf 2.2.4. Theorem.} {\it Under the hypotheses of 2.2.3, if $E_U$ is a closed subgroup of $GL(\scrF)_U$, then $E$ is a closed subgroup of $GL(\scrF)$.}

\medskip
\proof
As in the proof of 2.2.3 we argue that to prove the Theorem it suffices to consider the case when $Y=\Spec(R)$ is local and strictly henselian. As $E$ is an affine $R$-scheme, we have a natural affine morphism $q:E\to GL(\scrF)$ of $R$-schemes of finite type. It is a homomorphism as this is so generically. But from 2.2.2 we get that each fibre of $q$ is a closed embedding and that $q$ is a proper morphism. So $q$ is a finite morphism. So as $R$ is noetherian, from Nakayama's Lemma we get that $q$ is a closed embedding. This ends the proof.

\medskip\noindent
{\bf 2.2.5. Proposition.} {\it Let $m\in\dbN$. For $j\in\{1,...,m\}$ let $G_{j\eta}$ be a reductive subgroup of $GL(\scrF)_{\eta}$ such that the groups $G_{j\eta}$'s commute among themselves and the direct sum $\oplus_{j=1}^m \Lie(G_{j\eta})$ is a Lie subalgebra of $\Lie(GL(\scrF)_{\eta})$. Let $G_{0\eta}$ be the reductive subgroup of $GL(\scrF)_{\eta}$ generated by $G_{j\eta}$'s. 

\medskip
{\bf 1)} We have $\Lie(G_{0\eta})=\oplus_{j=1}^m \Lie(G_{j\eta})$. 

\smallskip
{\bf 2)} We assume that the Zariski closure $G_j$ of $G_{j\eta}$ in $GL(\scrF)$ is a reductive group scheme over $Y$, $\forall j\in\{1,...,m\}$. Then the Zariski closure $G_0$ of $G_{0\eta}$ in $GL(\scrF)$ is a reductive group scheme over $Y$.}

\medskip
\proof
Using induction on $m\in\dbN$, it suffices to prove the Proposition for $m=2$. As 1) and 2) are local statements for the \'etale topology of $Y$, from 2.2.1 we get that it suffices to prove the Proposition under the extra assumption that $G_1$ and $G_2$ are split and $Y$ is affine and connected. Let $C_{\eta}:=G_{1\eta}\cap G_{2\eta}$. It is a subgroup of $G_{j\eta}$ commuting with $G_{j\eta}$, $j\in\{1,2\}$. The Lie algebra $\Lie(C_{\eta})$ is included in $\Lie(G_{1\eta})\cap\Lie(G_{2\eta})=\{0\}$ and so is trivial. So $C_{\eta}$ is a finite, \'etale subgroup of $Z(G_{j\eta})$. Let $C$ be the Zariski closure of $C_{\eta}$ in $GL(\scrF)$. It is contained in any split torus $T_j$ of $G_j$. 

We consider the product homomorphism $T_1\times_Y T_2\to GL(\scrF)$. Its kernel $KER$ is a group scheme of multiplicative type isomorphic to $T_1\cap T_2$. The fibre over $\Spec(\eta)$ of $KER$ is isomorphic to $C_{\eta}$ and so is a finite, flat group scheme. So $KER$ is a finite, flat group scheme of multiplicative type. So $T_1\cap T_2$ is a finite, flat group scheme of multiplicative type and so equal to $C$. We embed $C$ in $G_1\times_Y G_2$ via the natural embedding $C\hookrightarrow G_1$ and via the composite of the inverse isomorphism $C\arrowsim C$ with the natural embedding $C\hookrightarrow G_2$. We have a natural product homomorphism $q:(G_1\times_Y G_2)/C\to GL(\scrF)$ which over $\eta$ can be identified naturally with the monomorphism $G_{0\eta}\hookrightarrow GL(\scrF)_{\eta}$. But again from 2.2.2 we get that each fibre of $q$ is a closed embedding and that $q$ is a proper morphism. We check that $q$ is a closed embedding. Argument: we can assume that $Y$ is finitely generated over $\Spec(\dbZ)$ and so the spectrum of a noetherian $\dbZ$-algebra; so as in the proof of 2.2.4 we get that $q$ is a closed embedding. So $G_0=\im(q)\arrowsim (G_1\times_Y G_2)/C$ is a reductive group scheme over $Y$. So 2) holds. Also as $\Lie(C_{\eta})=\{0\}$, we get that $\Lie(G_{0\eta})=\Lie(G_{1\eta})\oplus\Lie(G_{2\eta})$. So 1) also holds. This ends the proof. 

\bigskip\noindent
{\bf 2.3. Canonical dilatations.}  
We assume $k=\overline{k}$. Let $\tilde G$ be an affine, flat group scheme over $W(k)$. Let $a\in\tilde G(W(k))$. The N\'eron measure of the defect of smoothness $\delta(a)\in\dbN\cup\{0\}$ of $\tilde G$ at $a$ is the length of the torsion part of $a^*(\Omega_{\tilde G/\Spec(W(k))})$. As $\tilde G$ is a group scheme, the value of $\delta(a)$ does not depend on $a$ and so we denote it by $\delta(\tilde G)$. We have $\delta(\tilde G)\in\dbN$ iff $\tilde G$ is not smooth, cf. [BLR, Lemma of p. 65]. Let $\tilde S_k$ be the Zariski closure in $\tilde G_k$ of all special fibres of $W(k)$-valued points of $\tilde G$. It is a reduced subgroup of $\tilde G_k$. We write $\tilde S_k=\Spec(R_{\tilde G}/J_{\tilde G})$, where $\tilde G=\Spec(R_{\tilde G})$ and $J_{\tilde G}$ is the ideal of $R_{\tilde G}$ defining $\tilde S_k$. 

By the canonical dilatation of $\tilde G$ we mean the affine $\tilde G$-scheme $\tilde G_1=\Spec(R_{\tilde G_1})$, where $R_{\tilde G_1}$ is the $R_{\tilde G}$-subalgebra of $R_{\tilde G}[{1\over p}]$ generated by ${x\over p}$ with $x\in J_{\tilde G}$. The $W(k)$-scheme $\tilde G_1$ is flat and has a canonical group structure such that the morphism $\tilde G_1\to \tilde G$ is a homomorphism of group schemes over $W(k)$, cf. [BLR, p. 63 and (d) of p. 64]. The morphism $\tilde G_1\to \tilde G$ has the following universal property: any morphism $\tilde Y\to \tilde G$ of flat $W(k)$-schemes whose special fibre factors through the closed embedding $\tilde S_k\hookrightarrow \tilde G_k$, factors through $\tilde G_1\to \tilde G$ (cf. [BLR, (b) of p. 63]). If $\tilde G$ is smooth, then $\tilde S_k=\tilde G_k$ and so $\tilde G_1=\tilde G$.

Either $\tilde G_1$ is smooth or we have $0<\delta(\tilde G_1)<\delta(\tilde G)$, cf. [BLR, Prop. 5 of p. 68]. So by using at most $\delta(\tilde G)$ canonical dilatations (the first one of $\tilde G$, the second one of $\tilde G_1$, etc.), we get the existence of a unique smooth, affine group scheme $\tilde G^\prime$ over $W(k)$ endowed with a homomorphism $\tilde G^\prime\to \tilde G$ whose fibre over $B(k)$ is an isomorphism and which has the following universal property: any morphism $\tilde Y\to \tilde G$ of $W(k)$-schemes, with $\tilde Y$ a smooth $W(k)$-scheme, factors through $\tilde G^\prime\to \tilde G$. We refer to $\tilde G^\prime$ as the universal smoothening of $\tilde G$.

\bigskip
\noindent
{\boldsectionfont \S3. On $\dbZ_{(p)}$ embeddings}
\bigskip

Let $(G_1,X_1)$ be a Shimura pair. Let $p$ be a prime such that $G_{1\dbQ_p}$ is unramified. The fact that $G_{1\dbQ_p}$ is unramified is equivalent to the fact that $G_{1\dbQ_p}$ extends to a reductive group scheme $G_{1\dbZ_p}$ over $\dbZ_p$, cf. [Ti2, 1.10.2 and 3.8.1]. Any subgroup of $G_{1\dbQ_p}(\dbQ_p)$ of the form $G_{1\dbZ_p}(\dbZ_p)$ is called hyperspecial. The hyperspecial subgroups of $G_{1\dbQ_p}(\dbQ_p)$ are precisely the compact subgroups of $G_{1\dbQ_p}(\dbQ_p)$ of maximum volume with respect to the natural Haar measure on $G_{1\dbQ_p}(\dbQ_p)$, cf. [Ti2, 3.2]. Let $H_1:=G_{1\dbZ_p}(\dbZ_p)$. It is known that $H_1$ determines $G_{1\dbZ_{p}}$ (see [Ti2, 3.8.1] and [Va1, 3.1.2.1 e)]). Let $G_{1\dbZ_{(p)}}$ be the unique extension of $G_1$ to a reductive group scheme over $\dbZ_{(p)}$ whose pull back to $\dbZ_p$ is $G_{1\dbZ_p}$, cf. [Va1, 3.1.3]. Let $v_1$ be a  prime of $E(G_1,X_1)$ dividing $p$. The prime $v_1$ is unramified over $p$, cf. [Mi3, Cor. 4.7 (a)]. So $e(v_1)=1$. We recall from \S1 that we call $(G_1,X_1,H_1,v_1)$ a Shimura quadruple. In 3.1 to 3.3 we impose different assumptions on $(G_1,X_1,H_1,v_1)$ and we consider different properties of $(G_1,X_1,H_1,v_1)$ and of other similar quadruples related to it.

\bigskip\noindent
{\bf 3.1. General properties.} We consider an injective map 
$$f_1\colon (G,X)\hookrightarrow (G_1,X_1)$$ 
of Shimura pairs. Let $G_{\dbZ_{(p)}}$ be the Zariski closure of $G$ in $G_{1\dbZ_{(p)}}$. Let $H:=G(\dbQ_p)\cap H_1=G_{\dbZ_{(p)}}(\dbZ_p)$. Let $v$ be a prime of $E(G,X)$ dividing $v_1$. So $v$ divides also $p$. We also view $f_1:(G,X,H,v)\hookrightarrow (G_1,X_1,H_1,v_1)$ as an injective map between Shimura quadruples. 

The functorial morphism $\Sh(G,X)\to \Sh(G_1,X_1)_{E(G,X)}$ (see [De1, 5.4]) is a closed embedding as it is so over $\dbC$ (cf. [De1, 1.15]). So the natural morphism $\Sh_H(G,X)\to \Sh_{H_1}(G_1,X_1)_{E(G,X)}$ is a pro-finite morphism. We recall that any $O_{(v)}$-scheme of finite type is excellent (for instance, cf. [Ma, (34.A) and (34.B)]).

\medskip\noindent
{\bf 3.1.1. Proposition.} {\it We assume that $\Sh_{H_1}(G_1,X_1)$ has an integral canonical model $\scrN_1$ over $O_{(v_1)}$. We also assume that either $Z(G_1)$ is a subgroup of $G$ or $Z(G_1)(\dbQ)$ is discrete in $Z(G_1)(\dbA_f)$. Then we have the following four properties:

\medskip
{\bf 1)} The normalization $\scrN^\prime$ of $\scrN_{1O_{(v)}}$ in the ring of fractions of $\Sh_H(G,X)$ is a normal integral model of $(G,X,H,v)$ having the extension property. 

\smallskip
{\bf 2)} The $O_{(v)}$-scheme $\scrN^\prime$ is a pro-\'etale cover of a normal, faithfully flat $O_{(v)}$-scheme $\scrP_1$ of finite type and of relative dimension equal to $\dim_{\dbC}(X)$.

\smallskip
{\bf 3)} The morphism $\scrN^\prime\to\scrN_{1O_{(v)}}$ is finite.

\smallskip
{\bf 4)} If $(G_1,X_1)$ is a Siegel modular variety, then $\scrP_1$ is a quasi-projective $O_{(v)}$-scheme.}

\medskip
\proof
For $i\in\{1,2\}$ let $H_{0i1}$ be a compact, open subgroup of $G_1(\dbA_f^{(p)})$ such that $\scrN_1$ is a pro-\'etale cover of $\scrN_1/H_{0i1}$. Let $H_{0i}$ be a compact, open subgroup of $G(\dbA_f^{(p)})\cap H_{0i1}$ such that $\Sh(G,X)$ is a pro-\'etale cover of $\Sh_{H_{0i}\times H}(G,X)$. The morphism $\Sh_{H_{0i}\times H}(G,X)_{\dbC}\to\Sh_{H_{0i1}\times H_1}(G_1,X_1)_{\dbC}$ is of finite type as well as a formally closed embedding at each $\dbC$-valued point of $\Sh_{H_{0i}\times H}(G,X)_{\dbC}$. Let $\scrP_i$ be the normalization of $\scrN_{1O_{(v)}}/H_{0i1}$ in $\Sh_{H_{0i}\times H}(G,X)$. The $O_{(v)}$-scheme $\scrN_{1O_{(v)}}/H_{0i1}$ is faithfully flat, of finite type and so also excellent. So $\scrP_i$ is a normal scheme finite over $\scrN_{1O_{(v)}}/H_{0i1}$ and so also a faithfully flat $O_{(v)}$-scheme of finite type. So the relative dimension of $\scrP_i$ over $O_{(v)}$ is $\dim_{\dbC}(\Sh_{H_{0i}\times H}(G,X))=\dim_{\dbC}(X)$. The scheme $\scrP_i$ is naturally an $\scrN_{1O_{(v)}}/H_{0i}$-scheme. 

We now take $H_{01}$ such that it is also a normal subgroup of $H_{02}$. The natural morphism 
$$q_{12}:\scrP_1\to\scrP_2\times_{\scrN_{1O_{(v)}}/H_{02}} \scrN_{1O_{(v)}}/H_{01}$$ 
of normal schemes is pro-finite. We now check that $q_{12E(G,X)}$ is an open closed embedding. As $q_{12E(G,X)}$ is a pro-\'etale cover between normal $E(G,X)$-schemes of finite type, it is enough to show that $q_{12}(\dbC)$ is injective. We have 
$$\Sh_{H_{0i}\times H_1}(G_1,X_1)(\dbC)=G_{1\dbZ_{(p)}}(\dbZ_{(p)})\backslash X_1\times G_1(\dbA_f^{(p)})/H_{0i}\overline{Z(G_{1\dbZ_{(p)}})(\dbZ_{(p)})},$$ 
where $\overline{Z(G_{1\dbZ_{(p)}})(\dbZ_{(p)})}$ is the Zariski closure of $Z(G_{1\dbZ_{(p)}})(\dbZ_{(p)})$ in $Z(G_{1\dbZ_{(p)}})(\dbA_f^{(p)})$ (cf. [Mi3, Prop. 4.11]). Also we have a natural disjoint union decomposition
$$\Sh_{H_{0i}\times H}(G,X)(\dbC)=\cup_{[g_j]\in G(\dbQ)\backslash G(\dbQ_p)/H} C_j\backslash X\times G(\dbA_f^{(p)})/H_{0i},$$
where $g_j\in G(\dbQ_p)$ is a representative of the class $[g_j]\in G(\dbQ)\backslash G(\dbQ_p)/H$ and $C_j:=G(\dbQ)\cap g_jHg_j^{-1}$ is independent of $i\in\{1,2\}$. As $G_1(\dbQ_p)=G_1(\dbQ)H_1$ (cf. [Mi3, Lemma 4.9]), we can write $g_j=a_jh_j$, with $a_j\in G_1(\dbQ)$ and $h_j\in H_1$. So $C_j\leqslant G_1(\dbQ)\cap g_jH_1g_j^{-1}=G_1(\dbQ)\cap a_jH_1a_j^{-1}=a_jG_{1\dbZ_{(p)}}(\dbZ_{(p)}))a_j^{-1}=:C_{1j}$. Even more, we have $C_j=G(\dbQ)\cap C_{1j}$. This is so as $g_jHg_j^{-1}$ is the group of $\dbZ_p$-valued points of the Zariski closure of $G$ in $a_jG_{1\dbZ_{(p)}}a_j^{-1}$.

We first show that $q_{12}(\dbC)$ is injective if $Z(G_1)(\dbQ)$ is discrete in $Z(G_1)(\dbA_f^{(p)})$. So $\overline{Z(G_{1\dbZ_{(p)}})(\dbZ_{(p)})}=Z(G_{1\dbZ_{(p)}})(\dbZ_{(p)})$ and so we have $\Sh_{H_{0i}\times H_1}(G_1,X_1)(\dbC)=G_{1\dbZ_{(p)}}(\dbZ_{(p)})\backslash X_1\times G_1(\dbA_f^{(p)})/H_{0i}$. So to show that $q_{12}(\dbC)$ is injective, it suffices to show that each one of the following commutative diagrams 
$$
\spreadmatrixlines{1\jot}
\CD 
C_j\backslash X\times G(\dbA_f^{(p)})/H_{01} @>{s_1}>> G_{1\dbZ_{(p)}}(\dbZ_{(p)})\backslash X_1\times G_1(\dbA_f^{(p)})/H_{01}\\
@V{q}VV @VV{q_1}V\\
C_j\backslash X\times G(\dbA_f^{(p)})/H_{02}@>{s_2}>> G_{1\dbZ_{(p)}}(\dbZ_{(p)})\backslash X_1\times G_1(\dbA_f^{(p)})/H_{02},
\endCD
$$
is such that the maps $q$ and $s_1$ define an injective map of $C_j\backslash X\times G(\dbA_f^{(p)})/H_{01}$ into the fibre product of $s_2$ and $q_1$. The vertical maps $q$ and $q_1$ are the natural projections. The horizontal maps $s_1$ and $s_2$ are defined by the rule: the equivalence class $[h,g]$ is mapped into the equivalence class $[a_j^{-1}h,a_j^{-1}g]$, where $h\in X$ and $g\in G(\dbA_f^{(p)})$. So the fact that $q$ and $s_1$ define an injective map of $C_j\backslash X\times G(\dbA_f^{(p)})/H_{01}$ into the fibre product of $s_2$ and $q_1$ is an immediate consequence of the fact that $C_j=G(\dbQ)\cap C_{1j}$. So $q_{12}(\dbC)$ is injective. 

We now show that $q_{12}(\dbC)$ is injective if $Z(G_1)$ is a subgroup of $G$. So $H_{0i}\overline{Z(G_{1\dbZ_{(p)}})(\dbZ_{(p)})}$ is a subgroup of $G(\dbA_f^{(p)})$. So entirely as in the previous paragraph we argue that $q_{12}(\dbC)$ is injective (we just need to replace ``$/H_{0i}$" by ``$/H_{0i}\overline{Z(G_{1\dbZ_{(p)}})(\dbZ_{(p)})}$"). 

So $q_{12E(G,X)}$ is an open closed embedding. So as $q_{12}$ is also a pro-finite morphism of normal $O_{(v)}$-schemes of finite type, $q_{12}$ itself is an open closed embedding. So $\scrP_1$ is a pro-\'etale cover of $\scrP_2$ which in characteristic 0 is an \'etale cover inducing Galois covers between connected components. So $\scrP_1$ is an \'etale cover of $\scrP_2$ inducing Galois covers between connected components. By allowing $H_{01}$ to vary among the normal, open subgroups of $H_{02}$ and by natural passage to limits, we get that $\scrN^\prime$ is a pro-\'etale cover of $\scrP_2$ and so also that $\scrP_2=\scrN^\prime/H_{02}$. So the $O_{(v)}$-scheme $\scrN^\prime$ together with the natural action of $G(\dbA_f^{(p)})$ on it is a normal integral model of $\Sh_H(G,X)$ over $O_{(v)}$. It has the extension property as $\scrN_1$ has it. So 1) holds. As $\scrN^\prime$ is a pro-\'etale cover of $\scrP_2$, 2) holds too.  

As each $q_{12}$ is an open closed embedding, by allowing $H_{01}$ to vary through all normal, open subgroups of $H_{02}$ we also get that $\scrN^\prime$ is an open closed subscheme of $\scrP_2\times_{\scrN_{1O_{(v)}}} \scrN_{1O_{(v)}}$ and so is a finite $\scrN_{1O_{(v)}}$-scheme. So 3) holds.

To prove 4) we choose $H_{0i2}$ such that $\scrN_1/H_{0i2}$ is a moduli scheme of principally polarized abelian schemes having level $N$ symplectic similitude structure, where $N\in\dbN\setminus\{1,2\}$ is prime to $p$. Based on 2), to prove 4) it suffices to check that $\scrN_1/H_{0i2}$ is a quasi-projective scheme over $O_{(v_1)}=\dbZ_{(p)}$. But this is implied by [MFK, Ths. 7.9 and 7.10 of p. 139]. This ends the proof.

\medskip\noindent
{\bf 3.1.2. Lemma.} {\it We assume the hypotheses of 3.1.1 hold. Let $\scrN$ be the formally smooth locus of $\scrN^\prime$ over $O_{(v)}$. Then $\scrN$ is an open subcheme of $\scrN^\prime$ and we have $\scrN_{E(G,X)}=\scrN^\prime_{E(G,X)}$. Moreover, if $\scrN_{k(v)}$ is a non-empty scheme, then $\scrN$ together with the resulting action of $G(\dbA_f^{(p)})$ on it is a smooth integral model of $\Sh_H(G,X)$ over $O_{(v)}$.}

\medskip
\proof
Always $\scrN^\prime$ is a normal integral model of $\Sh_H(G,X)$ having the extension property, cf. 3.1.1 1). As $\scrN^\prime$ is a pro-\'etale cover of an excellent scheme $\scrP_1$ which is a normal $O_{(v)}$-scheme of finite type (see 3.1.1 2)), $\scrN$ is a pro-\'etale cover of the open subscheme $\scrU_1$ of $\scrP_1$ which is the smooth locus of $\scrP_1$ over $O_{(v)}$ and so $\scrN$ is an open subscheme of $\scrN^\prime$. It is well known that $\scrP_{1E(G,X)}$ is a smooth $E(G,X)$-scheme and so we have $\scrN_{E(G,X)}=\scrN^\prime_{E(G,X)}$. Obviously the open subscheme $\scrN$ of $\scrN^\prime$ is $G(\dbA_f^{(p)})$-invariant. So if $\scrN_{k(v)}$ is a non-empty scheme, then $\scrN$ together with the resulting action of $G(\dbA_f^{(p)})$ on it is a smooth integral model of $\Sh_H(G,X)$ over $O_{(v)}$. This ends the proof.

\medskip\noindent
{\bf 3.1.3. Lemma.} {\it We assume that the hypotheses of 3.1.1 hold and that $e(v)\Le\max\{1,p-2\}$. Let $\scrN$ be as in 3.1.2. Then the following three statements are equivalent:

\medskip
{\bf (i)} the $O_{(v)}$-scheme $\scrN$ has the extension property;

\smallskip
{\bf (ii)} $\scrN$ is the integral canonical model of $\Sh_H(G,X)$ over $O_{(v)}$;

\smallskip
{\bf (iii)} $\scrN=\scrN^\prime$.}

\medskip
\proof
As $\scrN$ is a healthy regular scheme (cf. 1.2.1), the equivalence $(i)\Leftrightarrow (ii)$ follows directly from 3.1.2 and def. 1.1.1 3). The implication $(iii)\Rightarrow (i)$ follows from 3.1.1 1). The spectrum of any faithfully flat $O_{(v)}$-algebra $V$ which is a DVR is a healthy regular scheme. So if (i) holds, then any morphism $\Spec(V)\to\scrN^\prime$ of $O_{(v)}$-schemes factors through $\scrN$. This implies that any morphism $\Spec(V)\to\scrP_1$ of $O_{(v)}$-schemes factors through $\scrU_1$. So we have $\scrU_1=\scrP_1$ and so $\scrN=\scrN^\prime$. So $(i)\Rightarrow (iii)$. So the above three statements are equivalent. This ends the proof.

\medskip\noindent
{\bf 3.1.4. Proposition.} {\it We assume that the hypotheses of 3.1.1 hold, that $e(v)\Le\max\{1,p-2\}$ and that $\Sh_H(G,X)$ has an integral canonical model $\scrN$ over $O_{(v)}$. Then $\scrN=\scrN^\prime$.}

\medskip
\proof
As $\scrN$ is a healthy regular scheme (cf. def. 1.1.1 3)) and as $\scrN^\prime$ has the extension property, we have a natural morphism $\scrN\to\scrN^\prime$ respecting the $G(\dbA_f^{(p)})$-actions and extending the natural identifications $\scrN_{E(G,X)}=\scrN^\prime_{E(G,X)}=\Sh_H(G,X)$. Let $H_0$ be a compact, open subgroup of $G(\dbA_f^{(p)})$ such that for any inclusion $H_2\leqslant H_1$ of open subgroups of $H_0$, the morphisms $\scrN/H_2\to \scrN/H_1$ and $\scrN^\prime/H_2\to \scrN^\prime/H_1$ are \'etale covers. To end the proof it suffices to show that the induced morphism $q:\scrN/H_0\to\scrN^\prime/H_0$ is an isomorphism. 

Let $q_V:\Spec(V)\to\scrN^\prime/H_0$ be a morphism with $V$ a strictly henselian DVR  of mixed characteristic. It lifts to a morphism $q_V^{\infty}:\Spec(V)\to\scrN^\prime$. As $\scrN$ also has the extension property and as $\Spec(V)$ is a healthy regular scheme, $q_V^{\infty}$ factors through $\scrN$. So $q_V$ factors through $\scrN/H_0$. So $q$ satisfies the valuative criterion of properness with respect to discrete valuation rings of mixed characteristic $(0,p)$. So as $q_{E(G,X)}$ is an isomorphism and as $\scrN/H_0$ is a flat $O_{(v)}$-scheme, we get that $q$ is a separated morphism of $O_{(v)}$-schemes. From Nagata's embedding theorem (see [Na] and [Vo]) we get that $\scrN/H_0$ is an open, Zariski dense subscheme of a proper, normal $\scrN^\prime/H_0$-scheme $(\scrN/H_0)^c$ which is flat over $O_{(v)}$. So as $q$ satisfies the valuative criterion of properness with respect to discrete valuation rings of mixed characteristic $(0,p)$, we get that $\scrN/H_0=(\scrN/H_0)^c$. So $q$ is a proper morphism of $O_{(v)}$-schemes. 

From the smoothening process of [BLR, Th. 3 of p. 61] we get the existence of an $\scrN^\prime/H_0$-scheme $Z$, smooth over $O_{(v)}$, quasi-projective over $\scrN^\prime/H_0$, whose fibre over $E(G,X)$ is $\scrN^\prime_{E(G,X)}/H_0=\Sh_{H_0\times H}(G,X)$ and such that the induced map $Z(O_{(v)}^{\sh})\to \scrN^\prime/H_0(O_{(v)}^{\sh})$ is a bijection and the smooth locus $\scrS$ of $\scrN^\prime/H_0$ over $O_{(v)}$ is naturally an open subscheme of $Z$. The scheme $Z\times_{\scrN^\prime/H_0} \scrN^\prime$ is regular, formally smooth over $O_{(v)}$ and so (cf. 1.2.1) healthy regular. So as $\scrN$ has the extension property, $Z\times_{\scrN^\prime/H_0} \scrN^\prime$ is naturally an $\scrN$-scheme. So $Z$ itself is naturally an $\scrN/H_0$-scheme. 

To show that $q$ is an isomorphism it suffices to show that $q$ is an isomorphism in codimension 1. We show that the assumption that $q$ is not an isomorphism in codimension 1 leads to a contradiction. This assumption implies that there is a connected component $\scrC_{k(v)}$ of $\scrN_{k(v)}/H_0$ dominating a reduced closed subscheme $Z_0$ of $\scrN^\prime_{k(v)}/H_0$ of dimension $d<\dim(\scrC_{k(v)})$. As $\scrS$ is an open subscheme of $Z$ and so an $\scrN/H_0$-scheme, the subscheme $Z_0$ of $\scrN^\prime/H_0$ is a closed subscheme of the complement of $\scrS$ in $\scrN^\prime/H_0$. Let $\scrC$ be the open subscheme of $\scrN/H_0$ defined by $\scrC_{k(v)}$ and the generic fibre of $\scrN/H_0$. The morphism $\scrC\to\scrN^\prime/H_0$ lifts to a morphism $\tilde q:\scrC\to\tilde Z$, where $\tilde Z$ is obtained from $\scrN^\prime/H_0$ through the first dilatation needed to get $Z$. Argument: this first dilatation is obtained by blowing up the maximal reduced closed subscheme $Z_{00}$ of $\scrN^\prime_{k(v)}/H_0$ having the property that it is included in $\scrN^\prime/H_0\setminus\scrS$ and the points of it with values in the residue field $\dbF$ of $O_{(v)}^{\text{\sh}}$ and which admit lifts (in $\scrN^\prime/H_0$) to $O^{\text{\sh}}$-valued points, are Zariski dense in it (cf. [BLR, p. 72]); so as $Z_0$ is a closed subscheme of $Z_{00}$, the existence of $\tilde q$ is implied by [BLR, Prop. 1 of p. 63]. We know that $\tilde Z$ is an affine $\scrN^\prime/H_0$-scheme, cf. the properties of dilatations [BLR, p. 62]. So as $\scrC$ is smooth and as its fibres over points of $\scrN^\prime/H_0$ are proper schemes, $\tilde q$ dominates a closed subscheme $\tilde Z_0$ of $\tilde Z_{k(v)}$ of dimension $d$. As above we argue that $\tilde Z_0$ is included in the non-smooth locus of $\tilde Z$. By induction on the number of dilatations needed to construct $Z$ from $\scrN^\prime/H_0$, we get that $\scrC$ is a $Z$-scheme and that the image of $\scrC_{k(v)}$ in $Z$ is of dimension $d$. So the composite morphism $\scrC\to Z\to\scrN/H_0$ is not an open embedding. Contradiction. So $q$ is an isomorphism. This ends the proof of the Proposition.

\medskip\noindent
{\bf 3.1.5. Example.} With the notations of 1.3, we take $(G_1,X_1)$ and $\scrN_1$ to be $(GSp(W,\psi),S)$ and respectively $\scrM$. So the definitions of $\scrN^\prime$ in 1.3 and 3.1.1 1) coincide. Always, the rank 1 split torus $Z(GSp(W,\psi))$ is a subtorus of $G$. So the hypotheses of  3.1.1 hold. So 3.1.1 and 3.1.2 apply entirely. In particular, $\scrN$ is always a $G(\dbA_f^{(p)})$-invariant, open subscheme of $\scrN^\prime$. If moreover we have $e(v)\Le\max\{1,p-2\}$, then 3.1.3 and 3.1.4 also apply.

\bigskip\noindent
{\bf 3.2. Standard PEL situations.} 
In this section let $L$, $L_{(p)}$, $f:(G,X)\hookrightarrow (GSp(W,\psi),S)$, $G_{\dbZ_{(p)}}$, $E(G,X)$, $v$, $k(v)$, $O_{(v)}$, $\scrN^\prime$ and $\scrN$ be as in 1.3. In this section we assume that $G_{\dbZ_{(p)}}$ is a reductive group scheme and that there is a $\dbZ_{(p)}$-subalgebra $\scrB_{(p)}$ of $\End(L_{(p)})$ which over $W(\dbF)$ is a product of matrix algebras, which is self dual with respect to $\psi$ and is such that $G$ is the identity component of the subgroup of $GSp(W,\psi)$ fixing the elements of $\scrB_{(p)}$. We recall that in this case $\Sh(G,X)$ is called a Shimura variety of PEL type (see [Ko], etc.). We also refer to the triple $(f,L,v)$ as a standard PEL situation.

\medskip\noindent
{\bf 3.2.1. Theorem.} {\it Then $\scrN=\scrN^{\prime}$ is an integral canonical model of $\Sh_H(G,X)$ over $O_{(v)}$.}

\medskip
\proof 
As $G_{\dbZ_{(p)}}$ is a reductive group scheme, the group $G_{\dbQ_p}$ is unramified and so from [Mi3, Cor. 4.7 (a)] we get $e(v)=1$. So based on 3.1.3 and 3.1.5, it suffices to show that $\scrN^\prime=\scrN$ and so that $\scrN^\prime$ is formally smooth over $O_{(v)}$. The case $p>2$ is well known (see [Zi], [LR],  [Ko, \S5] and [Va1, 4.3.11 and 5.6.3]). We refer to [Va3, 1.3] for the case $p=2$. This ends the proof.

\bigskip\noindent
{\bf 3.3. Proposition.} {\it We assume $(G_1,X_1)$ is a simple, adjoint Shimura pair of abelian type. Then there is an injective map $f:(G,X)\hookrightarrow (GSp(W,\psi),S)$ and a $\dbZ_{(p)}$-lattice $L_{(p)}$ of $W$ on which $\psi$ induces a perfect alternating form, such that the following four properties hold:

\medskip
{\bf (i)} we have $(G^{\ad},X^{\ad})=(G_1,X_1)$;

\smallskip
{\bf (ii)} the Zariski closure $G_{\dbZ_{(p)}}$ of $G$ in $GL(L_{(p)})$ is a reductive group scheme whose adjoint is the adjoint group $G_{1\dbZ_{(p)}}$ introduced in the beginning of \S3;

\smallskip
{\bf (iii)} if $(G_1^\prime,X_1^\prime)$ is a Shimura pair of abelian type having $(G_1,X_1)$ as its adjoint, then $G^{\der}$ is an isogeny cover of $G_1^{\prime\der}$ (i.e. $G^{\der}$ is the maximal isogeny cover of $G_1$ allowed by the abelian types);

\smallskip
{\bf (iv)} if either $(G_1,X_1)$ is of $A_n$ type or $(G_1,X_1)$ is of $C_n$ or $D_n^{\dbH}$ type and moreover $(G_1,X_1)$ has no compact factors, then there is a $\dbZ_{(p)}$-subalgebra $\scrB_{(p)}$ of $\End(L_{(p)})$ as in 3.2 (i.e. $(f,L,v)$ is a standard PEL situation, for any $\dbZ$-lattice $L$ of $W$ on which $\psi$ induces a perfect alternating form and which has the property that $L_{(p)}=L\otimes_{\dbZ} \dbZ_{(p)}$).}

\medskip
\proof
If $(G_1,X_1)$ is of $A_n$ type, then this is proved in [Va4, Prop. 4.1]. Let now $(G_1,X_1)$ be of $B_n$ or $C_n$ type with $n\Ge 2$ or of $D_n^{\dbH}$ or $D_n^{\dbR}$ type with $n\Ge 4$. If $p\Ge 3$, then the Proposition is just a weaker form of [Va4, Th. 4.8 and Rm. 4.8.1 4)]. We now recall a slight variant of loc. cit. which is closer in spirit to [Va1, 6.5 and 6.6] and which has two extra features which will play key roles in \S6. We work with an arbitrary prime $p\Ge 2$.

Let $F_1$ be a totally real number field such that $G_1$ is $\Res_{F_1/\dbQ} \tilde G_1$, with $\tilde G_1$ as an absolutely simple, adjoint group over $F_1$ (cf. [De2, 2.3.4 (a)]). We write 
$$F_1\otimes_{\dbQ} \dbQ_p=\prod_{i\in I_p} F_{1i}$$ 
as a product of $p$-adic fields. As $G_{1\dbZ_p}$ is a reductive group scheme, it splits over an unramified extension of $\dbZ_p$ (see 2.2.1). So $G_{1\dbQ_p}$ splits over an unramified, finite field extension $F_{10}$ of $\dbQ_p$. So each $F_{1i}$ is a subfield of $F_{10}$. So $F_1$ is unramified over $p$. Let $F_2$ be a totally real number field containing $F_1$, unramified over $p$ and such that we can write 
$$F_2\otimes_{\dbQ} \dbQ_p=\prod_{i\in I_p} F_{2i}$$ 
as a product of $p$-adic fields indexed by the set $I_p$ and with each $F_{1i}$ as a subfield of $F_{2i}$. 

An element $i\in I_p$ will be called compact if the following property holds:

\medskip
{\bf (v)} {\it the extension of the normal subgroup $\Res_{F_{1i}/\dbQ_p} \tilde G_{F_{1i}}$ of $G_{1\dbQ_p}$ to $\dbC$ via an (any) $O_{(v_1)}$-embedding $W(k(v_1))\hookrightarrow\dbC$, is such that it has simple factors which are pull backs to $\Spec(\dbC)$ of compact factors of $G_{1\dbR}$.}

\medskip
The set $I_p$ has compact elements iff $(G_1,X_1)$ has compact factors. Let $E_2$ be a totally imaginary quadratic extension of $F_2$ unramified over $p$. If $(G_1,X_1)$ has no compact factors, then we take $F_2=F_1$. If $(G_1,X_1)$ has compact factors, then we take $F_2$ and $E_2$ such that the following two additional properties hold:

\medskip
{\bf (vi)} {\it there is a compact element $i_0\in I_p$ such that we have $F_{1i}=F_{2i}$ for all $i\in I_p\setminus\{i_0\}$; so we have $[F_2:F_1]=[F_{2i_0}:F_{1i_0}]$;}

\smallskip
{\bf (vii)} {\it for any compact element $i\in I_p$ the field $E_2\otimes_{F_2} F_{2i}$ is the unramified quadratic field extension of $F_{2i}$.}

\medskip
Let $\tilde G_{1F_{1(p)}}$ be the adjoint group scheme over $F_{1(p)}$ extending $\tilde G_1$ and such that we have $G_{1\dbZ_{(p)}}=\Res_{F_{1(p)}/\dbZ_{(p)}} \tilde G_{1F_{1(p)}}$. One constructs $\tilde G_{1F_{1(p)}}$ as a direct factor of $G_{1F_{1(p)}}$. 
Let $F_3$ be a number field containing $F_1$, unramified over $p$ and such that a maximal torus $T_{1F_{1(p)}}$ of $\tilde G_{1F_{1(p)}}$ splits over $F_{3(p)}$. Let $F_4$ be the composite field of $F_3$ and $F_2$. Warning: the number fields $F_3$ and $F_4$ are not in general totally real number fields. 

Let $G_{2\dbZ_{(p)}}:=\Res_{F_{2(p)}/\dbZ_{(p)}} \tilde G_{1F_{2(p)}}$. Let $G_2:=G_{2\dbQ}=\Res_{F_2/\dbQ} \tilde G_{1F_2}$. Let $(G_2,X_2)$ be the adjoint Shimura pair such that the natural monomorphism $G_1=\Res_{F_1/\dbQ} \tilde G_1\hookrightarrow G_2=\Res_{F_2/\dbQ} \tilde G_{1F_2}$ extends to an injective map $(G_1,X_1)\hookrightarrow (G_2,X_2)$ of Shimura pairs (it exists as $F_2$ is totally real).  

See [Bou, planche II to IV] for the standard notations of weights associated to the $B_n$, $C_n$ and $D_n$ Lie types. We consider the representation  $\rho_0:\tilde G^{\sc}_{1F_{4(p)}}\to GL(L_{1(p)})$ over $F_{4(p)}$ such that the following properties hold:

\medskip
{\bf (viii)} {\it $L_{1(p)}$ is a free $F_{4(p)}$-module of finite rank;}

\smallskip
{\bf (ix)} {\it the representation $\rho_0$ is a direct sum of representations associated to the following minuscule weights: $\varpi_n$ if $(G_1,X_1)$ is of $B_n$ type, $\varpi_1$ if $(G_1,X_1)$ is of $C_n$ type, $\varpi_1$ if $(G_1,X_1)$ is of $D_n^{\dbH}$ type, and $\varpi_{n-1}$ and $\varpi_n$ if $(G_1,X_1)$ is of $D_n^{\dbR}$ type.}

\medskip
Here the weights are with respect to the inverse image in $\tilde G^{\sc}_{1F_{4(p)}}$ of the split torus $T_{1F_{4(p)}}$ of $\tilde G_{1F_{4(p)}}$. Warning: if $\tilde G_1$ is of $D_4$ Dynkin type, then we choose the numbering of weights as in [De2, p. 272]. The representation $\rho_0$ is a closed embedding iff $(G_1,X_1)$ is of $B_n$, $C_n$ or $D_n^{\dbR}$ type. If $(G_1,X_1)$ is of $D_n^{\dbH}$ type, then $\rho_0$ is the composite of an isogeny of degree 2 with a monomorphism which is a closed embedding.

Let 
$$L_{(p)}:=E_{2(p)}\otimes_{F_{2(p)}} L_{1(p)}$$ 
but viewed as a $\dbZ_{(p)}$-module and not as an $E_{2(p)}$-module; it is naturally a $G^{\sc}_{2\dbZ_{(p)}}$-module and so also a $G^{\sc}_{1\dbZ_{(p)}}$-module. Let $W:=L_{(p)}[{1\over p}]$; it is a $\dbQ$--vector space on which $G_2^{\sc}$ and $G_1^{\sc}$ act naturally. Let $G^{\der}$ (resp. $\tilde G^{\der}$) be the semisimple subgroup of $GL(W)$ which is the image of $G_1^{\sc}$ (resp. of $G^{\sc}_2$) in $GL(W)$. The Zariski closure of $G^{\der}$ (resp. $\tilde G^{\der}$) in $GL(L_{(p)})$ is a semisimple subgroup $G^{\der}_{\dbZ_{(p)}}$ (resp. $\tilde G^{\der}_{\dbZ_{(p)}}$) of $GL(L_{(p)})$ having $G_{1\dbZ_{(p)}}$ (resp. $G_{2\dbZ_{(p)}}$) as its adjoint. As $L_{(p)}$ has also a natural structure of an $E_{2(p)}$-module, the torus $\Res_{E_{2(p)}/\dbZ_{(p)}} \dbG_{mE_{2(p)}}$ acts naturally and faithfully on $L_{(p)}$. 

Let $T_{\dbZ_{(p)}}$ be the torus of $GL(L_{(p)})$ generated by $Z(GL(L_{(p)}))$ and by the maximal subtorus $T^c_{\dbZ_{(p)}}$ of $\Res_{E_{2(p)}/\dbZ_{(p)}} \dbG_{mE_{2(p)}}$ which over $\dbR$ is compact. The torus $T^c_{\dbZ_{(p)}}$ is isogenous to the quotient torus $\Res_{E_{2(p)}/\dbZ_{(p)}} \dbG_{mE_{2(p)}}/\Res_{F_{2(p)}/\dbZ_{(p)}} \dbG_{mF_{2(p)}}$ and so has rank $[F_2:\dbQ]$. Let $G$ (resp. $\tilde G)$ be generated by a subtorus $T^\prime_{\dbQ}$ of $T_{\dbQ}$ and by $G^{\der}$ (resp. and by $\tilde G^{\der}$). The notations match, i.e. the derived group of $G$ is indeed the group $G^{\der}$ introduced in the previous paragraph. Moreover, the Zariski closure of $T_{\dbQ}^\prime=Z^0(G)$ in $GL(L_{(p)})$ is a subtorus of $T_{\dbZ_{(p)}}$. So from 2.2.5 2) we get that the Zariski closures $G_{\dbZ_{(p)}}$ and $\tilde G_{\dbZ_{(p)}}$ of $G$ and respectively $\tilde G$ in $GL(L_{(p)})$ are reductive group schemes. Their adjoints are $G_{1\dbZ_{(p)}}$ and respectively $G_{2\dbZ_{(p)}}$. So (ii) holds. Using a standard embedding $SL_{m\dbZ_{(p)}}\hookrightarrow Sp_{2m\dbZ_{(p)}}$ with $m:=\dim_{\dbQ}(W)$ and replacing if needed $L_{1(p)}$ by $L_{1(p)}\oplus L_{1(p)}^*$, we can assume that there is a perfect alternating form $\tilde\psi$ on $L_{(p)}$ normalized by $\tilde G_{\dbZ_{(p)}}$.

In [De2, proof of 2.3.10] it is shown that we can choose the torus $T^\prime_{\dbQ}$ (for instance, we can take $T^\prime_{\dbQ}$ to be $T_{\dbQ}$ itself) such that there is a monomorphism $x:\dbS\hookrightarrow\tilde G_{\dbR}$ having the following three properties:

\medskip
{\bf (x)} {\it the Hodge $\dbQ$--structure on $W$ defined by $x$ is of type $\{(-1,0),(0,-1)\}$;}

\smallskip
{\bf (xi)} {\it there is a non-degenerate alternating form $\psi$ on $W$ normalized by $\tilde G$ and such that $2\pi i\psi$ is a polarization of the Hodge $\dbQ$--structure on $W$ mentioned in property (x);}

\smallskip
{\bf (xii)} {\it if $\tilde X$ is the $\tilde G(\dbR)$-conjugacy class of $x$, then we get a Shimura pair $(\tilde G,\tilde X)$ whose adjoint is $(G_2,X_2)$.}

\medskip
The group $G_{2\dbZ_{(p)}}(\dbZ_{(p)})$ permutes transitively the connected components of $X_2$, cf. [Va1, 3.3.3]. So by replacing $x$ by its inner conjugate under an element of $\tilde G(\dbR)G_{2\dbZ_{(p)}}(\dbZ_{(p)})$ and by replacing if needed $\psi$ by another non-degenerate alternating form on $W$ such that (xi) holds (cf. [De2, 2.3.3]), we can choose $x$ such that its image in $X_2$ (i.e. the composite of $x$ with the natural epimorphism $\tilde G_{\dbR}\twoheadrightarrow \tilde G^{\ad}_{\dbR}=G_{2\dbR}$) belongs to $X_1$. So $x$ factors through $G_{\dbR}$ and so we can view $x$ also as a monomorphism $x:\dbS\hookrightarrow G_{\dbR}$. 

Let $ALT$ be the free $\dbZ_{(p)}$-module of alternating forms on $L_{(p)}$ normalized by $\tilde G_{\dbZ_{(p)}}$. So $\tilde\psi\in ALT$ and $\psi\in ALT[{1\over p}]$. Based on [De2, 1.1.18 b)] and on standard arguments of inequivalent valuations, we get that we can choose $\psi$ such that we also have $\psi\in ALT$  and moreover $\psi$ and $\tilde\psi$ are congruent mod $p$. So $\psi$ induces a perfect alternating form on $L_{(p)}$.

Let $X$ be the $G(\dbR)$-conjugacy class of $x:\dbS\hookrightarrow G_{\dbR}$. From (xii) we get that the pair $(G,X)$ is a Shimura pair. As $x$ maps into an element of $X_1$, the adjoint of $(G,X)$ is $(G_1,X_1)$. So (i) holds. We already checked that (ii) holds. The fact that (iii) holds is obvious if $(G_1,X_1)$ is of $B_n$, $C_n$ or $D_n^{\dbR}$ (as $G^{\der}$ is simply connected) and  is implied by [De2, 2.3.8 and 2.3.10 (ii)] if $(G_1,X_1)$ is of $D_n^{\dbH}$ type. 

Let now $(G_1,X_1)$ be of $C_n$ or $D_n^{\dbH}$ (resp. of $B_n$ or $D_n^{\dbR}$) type. Let $C_{\dbZ_{(p)}}$ be the centralizer of $G_{\dbZ_{(p)}}$ in $GL(L_{(p)})$. The representation of $G_{W(\dbF)}$ on $L_{(p)}\otimes_{\dbZ_{(p)}} W(\dbF)$ is a direct sum of irreducible representations of rank $2n$ (resp. of rank $2^{n}$ or $2^{n-1}$ and which are either spin or half spin representations and) on which a product factor of $G^{\der}_{W(\dbF)}$ isomorphic to $Sp_{2nW(\dbF)}$ or $SO_{2nW(\dbF)}$ (resp. isomorphic to $\text{Spin}_{2n+1W(\dbF)}$ or to a quotient of $\text{Spin}_{2nW(\dbF)}$ by $\mu_{2W(\dbF)}$) acts non-trivially. So the special fibres of any two such irreducible representations are also irreducible (see [Va5, Lemma of the proof of 7.8.3.1]) and moreover are isomorphic iff the two representations over $W(\dbF)$ are isomorphic. This implies that $C_{W(\dbF)}$ is a reductive group scheme. So $C_{\dbZ_{(p)}}$ is a reductive group scheme. Let $\scrB_{(p)}$ be the semisimple $\dbZ_{(p)}$-subalgebra of $\End(L_{(p)})$ such that as $\dbZ_{(p)}$-modules we have $\scrB_{(p)}=\Lie(C_{\dbZ_{(p)}})$. 

If $(G_1,X_1)$ is of $C_n$ or $D_n^{\dbH}$ type and has no compact factors, then we can choose $T^\prime_{\dbQ}$ to be $Z(GL(W))$ (cf. [De2, Rm. 2.3.13]) and so the centralizer $D$ of $\scrB_{(p)}[{1\over p}]$ in $GL(W)$ is a reductive group  such that $D_{\dbC}$ is isomorphic to $GL_{2n\dbC}^{[F_1:\dbQ]}$. Moreover the monomorphism $G^{\der}_{\dbC}\to D^{\der}_{\dbC}$ is a product of $[F_1:\dbQ]$-monomorphism of the form $Sp_{2n\dbC}\hookrightarrow SL_{2n\dbC}$ or of the form $SO_{2n\dbC}\hookrightarrow SL_{2n\dbC}$ and the representation on $W\otimes_{\dbQ} \dbC$ of each direct factor $SL_{2n\dbC}$ of $D^{\der}_{\dbC}$ is isotypic. So the Lie algebra of each such direct factor is normalized without being fixed by the involution of $\End(W\otimes_{\dbQ} \dbC)$ defined by $\psi$. We easily get that $\Lie(G)$ is the Lie algebra of the identity component $\tilde G^\prime$ of the centralizer of $\scrB_{(p)}[{1\over p}]$ in $GSp(W,\psi,S)$. So the natural monomorphism $G\hookrightarrow\tilde G^\prime$ is an isomorphism. So (iv) holds. This ends the proof. 

\medskip\noindent
{\bf 3.3.1. Remark.} We refer to the proof of 3.3. Always $\scrB_{(p)}$ is an $F_{1(p)}$-algebra and the action of $F_{1(p)}$ on $L_{(p)}$ has no fixed non-trivial subspaces. The centralizers of $T_{\dbZ_{(p)}}$ and $\Res_{E_{2(p)}/\dbZ_{(p)}} \dbG_{mE_{2(p)}}$ in $GL(L_{(p)})$ are the same. So if we choose $T_{\dbQ}^\prime$ to be $T_{\dbQ}$, then $\scrB_{(p)}$ is an $E_{2(p)}$-algebra and so also an $F_{2(p)}$-algebra.

\medskip\noindent
{\bf 3.3.2. Lemma.} {\it We refer to the proof of 3.3, with $(G_1,X_1)$ of $B_n$ or $C_n$ type with $n\Ge 1$ or of $D_n^{\dbH}$ type with $n\Ge 4$. Let $m$ be $2n$ if $(G_1,X_1)$ is of $C_n$ or $D_n^{\dbH}$ type and be $2^n$ if $(G_1,X_1)$ is of $B_n$ type. We take $T_{\dbQ}^\prime$ to be $T_{\dbQ}$. Then $Z(G)=Z^0(G)$ and the centralizer of $\scrB_{(p)}$ in $GSp(L_{(p)},\psi)$ is a reductive group scheme $\tilde G^\prime_{\dbZ_{(p)}}$ having $Z(G_{\dbZ_{(p)}})$ as its center and having a derived group which over $W(\dbF)$ is isomorphic to $SL_{mW(\dbF)}^{[F_2:\dbQ]}$.}

\medskip
\proof
The torus $T^c_{W(\dbF)}$ is a product $\prod_{j=1}^{[F_2:\dbQ]} T_j$ of rank 1 tori such that the representation of $T^c_{W(\dbF)}$ on $L_{(p)}\otimes_{\dbZ_{(p)}} W(\dbF)$ is a direct sum $\oplus_{j=1}^{2[F_2:\dbQ]} L_j$ having the following property:

\medskip
(*) {\it for any $j\in\{1,...,[F_2:\dbQ]\}$, the torus $T_j$ acts via the identical character on $L_{j}$, via the inverse of the identical character on $L_{j+[F_2:\dbQ]}$ and trivially on $L_{j_1}$, $\forall j_1\in\{1,...,2[F_2:\dbQ]\}\setminus\{j,j+[F_2:\dbQ]\}$.}

\medskip
As $T^c_{W(\dbF)}$ fixes $\psi$, for $j_1$, $j_2\in\{1,...,2[F_2:\dbQ]\}$ we have 

\medskip
(**) {\it $\psi(L_{j_1},L_{j_2})\neq \{0\}$ iff $|j_2-j_1|=[F_2:\dbQ]$.} 

\medskip
The derived group $G^{\der}_{W(\dbF)}$ acts on each $L_j$ via a direct factor $DF_j$ of it isomorphic to an $Sp_{2nW(\dbF)}$, an $SO_{2nW(\dbF)}$ or an $\text{Spin}_{2n+1W(\dbF)}$ group. Moreover, the representation of $DF_j$ on $L_j$ is a direct sum of isomorphic representations associated to a unique minuscule weight and which shows up in 3.3 (ix). This implies that $Z(G^{\der}_{W(\dbF)})\arrowsim\mu_2^{[F_1:\dbQ]}$ is a subgroup of $T^c_{W(\dbF)}$. So $Z(G_{W(\dbF)})=Z^0(G_{W(\dbF)})$ and so $Z(G)=Z^0(G)$. 

The group $C_{W(\dbF)}$ is a product $\prod_{j=1}^{2[F_2:\dbQ]} C_j$, with $C_j$ a $GL_{mW(\dbF)}$ group acting trivially on each $L_{j_1}$ with $j_1\in\{1,...,2[F_2:\dbQ]\}\setminus\{j\}$. The representation of $C_j$ on $L_j$ is a direct sum of standard rank $m$ representations. Due to (**), the intersection $C_{W(\dbF)}\cap Sp(L_{(p)}\otimes_{\dbZ_{(p)}} W(\dbF),\psi)$ is a product $\prod_{j=1}^{[F_2:\dbQ]} D_j$, with $D_j$ as a $GL_{mW(\dbF)}$ group having $T_j$ as its center. So the centralizer of $\scrB_{(p)}\otimes_{\dbZ_{(p)}} W(\dbF)$ in $GSp(L_{(p)}\otimes_{\dbZ_{(p)}} W(\dbF),\psi)$ is a reductive group scheme having $Z(G_{W(\dbF)})$ as its center and having a derived group isomorphic to $SL^{[F_2:\dbQ]}_{mW(\dbF)}$. So the centralizer $\tilde G^\prime_{\dbZ_{(p)}}$ of $\scrB_{(p)}$ in $GSp(L_{(p)},\psi)$ is a reductive group scheme such that its center is $Z(G_{\dbZ_{(p)}})$ and $\tilde G^{\prime\der}_{W(\dbF)}$ is isomorphic to $SL^{[F_2:\dbQ]}_{mW(\dbF)}$. This ends the proof. 

\medskip\noindent
{\bf 3.3.3. Remark.} We refer to the proof of 3.3, with $(G_1,X_1)$ of $D_n^{\dbR}$ type with $n\Ge 4$. If we choose $T_{\dbQ}^\prime$ to be $T_{\dbQ}$, then entirely as in the proof of 3.3.2 we argue that the centralizer of $\scrB_{(p)}$ in $GSp(L_{(p)},\psi)$ is a reductive group scheme $\tilde G^\prime_{\dbZ_{(p)}}$ whose derived group is such that its extension to $W(\dbF)$ is isomorphic to $SL_{2^{n-1}W(\dbF)}^{2[F_2:\dbQ]}$. The only difference is that this time we have two minuscule weights $\varpi_{n-1}$ and $\varpi_n$ (see 3.3 (ix)) and the representations associated to them are half spin representations and so are of rank $2^{n-1}$; this justifies the double number $2[F_2:\dbQ]$ of copies of $SL_{2^{n-1}W(\dbF)}$.  

\bigskip
\noindent
{\boldsectionfont \S4. Crystalline applications}
\bigskip

In 4.1 we recall the simplest variant of the main result of [dJ2]. In 4.2 we first introduce several notations needed to prove 1.4 to 1.6 and then we apply the main result of [Va6]. In 4.3 we prove 1.4. In 4.4 we present a simple criterion on when the special fibre $\scrN_{k(v)}$ of the $O_{(v)}$-scheme $\scrN$ of 1.3 is non-empty.  In 4.5 we apply 1.4 1) and 4.4 to get a solution to Langlands' conjecture of [La, p. 411] for the case of Shimura varieties of Hodge type. In 4.6 we refine the deformation theories of [Fa, \S7] and [Va1, 5.4] as allowed by [Va6]. 

For (crystalline or de Rham) Fontaine comparison theory we refer to [Fo]; see also [Fa, \S5] and [Va6]. Let $k$ be as in 2.1. Let $x$ be an independent variable. As the Verschiebung maps of $p$-divisible groups will not be mentioned at all in what follows, we will use the terminology $F$-crystals (resp. filtered) $F$-crystals associated to $p$-divisible groups over $k$, $k[[x]]$ or $k((x))$ (resp. over $W(k)$ or $W(k)[[x]]$) instead of the terminology Dieudonn\'e (resp. filtered Dieudonn\'e) $F$-crystals used in [BBM, Ch. 3], [BM, Ch. 2 and 3] or [dJ1]. 

The simplest form of [dJ2, 1.1] says:

\bigskip\noindent
{\bf 4.1. Theorem (de Jong).} {\it The natural functor from the category of non-degenerate $F$-crystals over $\Spec(k[[x]])$ to the category of non-degenerate $F$-crystals over $\Spec(k((x)))$ is fully faithful.}

\medskip
For the notion non-degenerate crystal [dJ2] refers to [Saa, 3.1.1 of p. 331]. For us the only important thing will be that any $F$-crystal of an abelian scheme over $\Spec(k[[x]])$ is non-degenerate and that the non-degenerate $F$-crystals are stable under tensor products. 

\bigskip\noindent
{\bf 4.2. A setting.} Warning: from now on until the end of the paper, the field $k$ will be assumed to be algebraically closed and the following notations $L$, $L_{(p)}$, $f:(G,X)\hookrightarrow (GSp(W,\psi),S)$, $E(G,X)$, $v$, $k(v)$, $O_{(v)}$, $K_p=GSp(L,\psi)(\dbZ_p)$, $G_{\dbZ_{(p)}}$, $H=G_{\dbZ_{(p)}}(\dbZ_p)$, $\scrM$, $\scrN^\prime$, $\scrN$, $(\scrA^\prime,\scrP_{\scrA^\prime})$ will be as in 1.3. 

We will also use the terminology of [De3] on Hodge cycles on an abelian scheme $A_Z$ over a reduced $\dbQ$--scheme $Z$. So we write each Hodge cycle $v$ on $A_Z$ as a pair $(v_{dR},v_{\acute et})$, where $v_{dR}$ and $v_{\acute et}$ are the de Rham and respectively the \'etale component of $v$. The \'etale component $v_{\acute et}$ as its turn has an $l$-component $v_{\acute et}^l$, for any rational prime $l$. For instance, if $Z$ is the spectrum of a field $E$, then $v_{\acute et}^p$ is a suitable $\Gal(\overline{E}/E)$-invariant tensor of the tensor algebra of $H^1_{\acute et}(A_{\overline{Z}},\dbQ_p)\oplus (H^1_{\acute et}(A_{\overline{Z}},\dbQ_p))^*\oplus \dbQ_p(1)$, where $\overline{Z}:=\Spec(\overline{E})$ and $\dbQ_p(1)$ is the usual Tate twist. If moreover $\overline{E}$ is a subfield of $\dbC$, then we also use the Betti realization of $v$: it corresponds to $v_{dR}$ (resp. to $v_{\acute et}^l$) via the canonical isomorphism relating the Betti cohomology with $\dbQ$--coefficients of $A_Z\times_Z \Spec(\dbC)$ with the de Rham (resp. the $\dbQ_l$ \'etale) cohomology of $A_{\overline{Z}}$ (see [De3, \S2]).

Let $TR$ be the trace form on $\End(W)$. If $K$ is a field of characteristic 0 and if $G_0$ is a reductive subgroup of $GL(W\otimes_{\dbQ} K)$, then the restriction $TR_0$ of $TR$ to $\Lie(G_0)$ is non-degenerate. Argument: it suffices to consider the cases when $K=\overline{K}$ and $G_0$ is either $\dbG_{mK}$ or a semisimple group whose adjoint is simple; the fact that $TR_0$ is non-degenerate is trivial in the first case and it follows from Cartan's solvability criterion in characteristic 0 in the second case. Let $\pi_{G_0}$ be the projector of $\End(W\otimes_{\dbQ} K)$ on $\Lie(G_0)$ along the perpendicular of $\Lie(G_0)$ with respect to $TR$. If $G_K$ normalizes $G_0$, then $G_K$ fixes $\pi_{G_0}$. 

As $Z(GL(W))$ is a subgroup of $G$, each tensor of $\scrT(W^*)$ fixed by $G$ belongs to $\oplus_{u\in\dbN\cup\{0\}} W^{*\otimes u}\otimes_{\dbQ} W^{\otimes u}$. So as $G$ is a reductive group, from [De3, 3.1 c)] we get the existence of a family $(v_{\alpha})_{\alpha\in\scrJ}$ of tensors in spaces of the form $W^{*\otimes u}\otimes_{\dbQ} W^{\otimes u}\subset\scrT(W^*)$ for some $u\in\dbN\cup\{0\}$, which contains $\pi_G$ and has the property that $G$ is the subgroup of $GL(W)$ fixing $v_{\alpha}$, $\forall\alpha\in\scrJ$. The choice of $L$ and $(v_\alpha)_{\alpha\in\scrJ}$ allows a moduli interpretation of $\Sh(G,X)$ (see [De1], [De2], [Mi3] and [Va1, 4.1 and 4.1.3]). For instance, $\Sh(G,X)(\dbC)=G(\dbQ)\backslash X\times G(\dbA_f)$ is the set of isomorphism classes of principally polarized abelian
varieties over $\dbC$ of dimension ${{\dim_{\dbQ}(W)}\over 2}$, carrying a family of Hodge
cycles indexed by $\scrJ$, having level $N$ symplectic similitude structure for any $N\in\dbN$ and satisfying some additional conditions. This interpretation endows the abelian scheme $\scrA^\prime_{E(G,X)}$ with a family $(w_{\alpha}^{\scrA^\prime})_{\alpha\in\scrJ}$ 
of Hodge cycles (the Betti realizations of pull backs of $w_{\alpha}^{\scrA^\prime}$ via $\dbC$-valued points of $\scrN_{E(G,X)}$ correspond naturally to $v_{\alpha}$).

Let $z\in\scrN^\prime(W(k))$. Let 
$$(A,p_A,(\tilde w_{\alpha})_{\alpha\in\scrJ}):=z^*(\scrA^\prime,\scrP_{\scrA^\prime},(w_{\alpha}^{\scrA^\prime})_{\alpha\in\scrJ}).$$ 
Let 
$$(M,F^1,\phi,p_M)$$
 be the principally quasi-polarized filtered $F$-crystal over $k$ of the principally quasi-polarized $p$-divisible group $(D,p_D)$ of $(A,p_A)$. So $p_M$ is a perfect alternating form on the free $W(k)$-module $M$ of rank $\dim_{\dbQ}(W)=2\dim_{W(k)}(A)$, $F^1$ is a maximal isotropic subspace of $M$ with respect to $p_M$ and the pair $(M,\phi)$ is a Dieudonn\'e module. Let $\tilde t_{\alpha}$ and $\tilde v_{\alpha}$ be the de Rham and respectively the $p$-component of the \'etale component of $\tilde w_{\alpha}$. We have $\tilde t_{\alpha}\in\scrT(M)[{1\over p}]$. Let $\tilde G$
be the Zariski closure in $GL(M)$ of the subgroup of $GL(M[{1\over p}])$ fixing $\tilde t_{\alpha}$, $\forall\alpha\in\scrJ$. 

It is known that $\tilde w_{\alpha}$ is a de Rham cycle; in other words, via de Rham and so also the crystalline Fontaine comparison theory, $\tilde t_{\alpha}$ and $\tilde v_{\alpha}$ correspond to each other. If $A_{B(k)}$ is definable over a number field contained in $B(k)$, then this was known previously (for instance, see [Bl, (0.3)]). The general case follows from loc. cit. and [Va1, 5.2.16] (in [Va1, 5.2] an odd prime is used; however the proof of [Va1, 5.2.16] applies for all primes). In particular we get that $\phi(\tilde t_{\alpha})=\tilde t_{\alpha}$, $\forall\alpha\in\scrJ$. Here $\phi$ acts on $M^*[{1\over p}]$ by mapping $e\in M^*[{1\over p}]$ into $\sigma\circ e\circ\phi^{-1}\in M^*[{1\over p}]$ and so it acts on $\scrT(M)[{1\over p}]$ in the natural tensor product way.

Let $\mu$ be the inverse of the canonical split cocharacter of $(M,F^1,\phi)$ defined in [Wi, p. 512]. It fixes $\tilde t_{\alpha}$, $\forall\alpha\in\scrJ$ (cf. the functorial aspects of [Wi, p. 513]). So $\mu$ factors through $\tilde G$; we denote also by 
$$\mu:\dbG_{mW(k)}\to\tilde G$$ 
the resulting factorization. Let 
$$M=F^1\oplus F^0$$ 
be the direct sum decomposition such that $\mu$ fixes $F^0$. Warning: often in what follows in connection to different Kodaira--Spencer maps, we will identity naturally $\Hom(F^1,F^0)$ with a direct summand of $\End(M)$ formed by endomorphisms of $M$ annihilating $F^0$. 

\medskip\noindent
{\bf 4.2.1. Fact.} {\it Let $w$ be a complex point of $\Sh_H(G,X)=\scrN_{E(G,X)}$. Let $(A_w,p_{A_w}):=w^*((\scrA^{\prime},\scrP_{\scrA^\prime})_{E(G,X)})$. We denote by $v^w_{\alpha}$ the $p$-component of the \'etale component of the Hodge cycle $w^*(w_{\alpha}^{\scrA^\prime})$ on $A_w$. Then there are isomorphisms $(H^1_{\acute et}(A_w,\dbZ_p),(v_{\alpha}^w)_{\alpha\in\scrJ})\arrowsim (L^*_{(p)}\otimes_{\dbZ_{(p)}} \dbZ_p,(v_{\alpha})_{\alpha\in\scrJ})$ taking the perfect bilinear form on $H^1_{\acute et}(A_w,\dbZ_p)$ defined by $p_{A_w}$ into a $\dbG_{m\dbZ_p}(\dbZ_p)$-multiple of the perfect bilinear form $\psi^*$ on $L^*_{(p)}\otimes_{\dbZ_{(p)}} \dbZ_p$ defined by $\psi$.}

\medskip
\proof
Let $w^{\infty}\in\Sh(G,X)(\dbC)$ lifting $w$. Let $\scrA^{\prime\infty}_{E(G,X)}$ be the pull back of $\scrA^{\prime}_{E(G,X)}$ to $\Sh(G,X)$. We identify canonically $H^1_{\acute et}(A_w,\dbZ_p)=H^1_{\acute et}(w^{\infty*}(\scrA^{\prime\infty}_{E(G,X)}),\dbZ_p)$. From the standard moduli interpretation of $\Sh(G,X)(\dbC)$ applied to $w^{\infty}\in\Sh(G,X)(\dbC)$ we get that $A_w$ as a complex manifold is $F^{0,-1}_w\backslash W\otimes_{\dbQ} \dbC/L_w$, where $L_w$ is a $\dbZ$-lattice of $W$ such that $L_w\otimes_{\dbZ} \widehat{\dbZ}$ is a $\widehat{\dbZ}$-lattice of $W\otimes_{\dbQ} \dbA_f=W\otimes_{\dbZ} \widehat{\dbZ}$ which is $G(\dbA_f)$-conjugate to $L\otimes_{\dbZ} \widehat{\dbZ}$ and where $W\otimes_{\dbQ} \dbC=F^{0,-1}_w\oplus F^{-1,0}_w$ is the Hodge decomposition of a Hodge $\dbQ$--structure on $W$ defined by some element $h_w\in X$ (see [Di1], [Mi2], [Mi3], [Va1, p. 454]). Moreover, the principal polarization $p_{A_w}$ of $A_w$ is defined naturally by a $\dbG_{m\dbQ}(\dbQ)$-multiple of $\psi$. So $(H^1_{\acute et}(A_w,\dbZ_p),(v_{\alpha}^w)_{\alpha\in\scrJ})$ is identified naturally with $(L_w^*\otimes_{\dbZ} \dbZ_p,(v_{\alpha})_{\alpha\in\scrJ})$ and so with a $G_{\dbQ_p}(\dbQ_p)$-conjugate of $(L^*_{(p)}\otimes_{\dbZ_{(p)}} \dbZ_p,(v_{\alpha})_{\alpha\in\scrJ})$. From this and the existence of the mentioned $\dbG_{m\dbQ}(\dbQ)$-multiple of $\psi$, the Fact follows. This ends the proof. 
 
\medskip\noindent
{\bf 4.2.2. Theorem.} {\it For the point $z\in\scrN^\prime(W(k))$ we consider the following condition 

\medskip
(*) there is a flat, closed subgroup $\tilde G^1$ of $\tilde G$, a set $\scrJ^1$ containing $\scrJ$ and a family of tensors $(\tilde t_{\alpha})_{\alpha\in\scrJ^1}$ of $\scrT(M)[{1\over p}]$ fixed by $\phi$ and $\mu_{B(k)}$, such that $\tilde G^1_{B(k)}$ is the subgroup of $GL(M[{1\over p}])$ fixing $\tilde t_{\alpha}$, $\forall\alpha\in\scrJ^1$, and for any $W(k)$-valued point $\tilde g$ of the universal smoothening $\tilde G^{1\prime}$ of $\tilde G^1$ lifting a $k$-valued point of the identity component of $\tilde G^{1\prime}_k$, the $F$-crystal $(M,\tilde g\phi)$ over $k$ has only positive slopes.

\medskip
If $p=2$ we assume that either $G_{\dbZ_{(p)}}$ is a torus or (*) holds. Then we have:

\medskip
{\bf 1)} There is an isomorphism $(M,(\tilde t_{\alpha})_{\alpha\in\scrJ})\arrowsim (H^1_{\acute et}(A_{B(k)},\dbZ_p)\otimes_{\dbZ_p} W(k),(\tilde v_{\alpha})_{\alpha\in\scrJ})$.

\smallskip
{\bf 2)} The group scheme $\tilde G$ is isomorphic to $G_{W(k)}$.}

\medskip
\proof
Part 1) follows from [Va6, 1.2 and 4.3 5)] applied in the context of the pair $(D,(\tilde t_{\alpha})_{\alpha\in\scrJ^1})$.${}^1$ $\vfootnote{1}{The assumption $G$ smooth used in [Va6, 3.4.5 2) and 4.3 5)] is not needed due to the following reason. Let $Q$, $k_0$, $\hat Q^0$ and $\tilde h$ be as in [Va6, 3.4.1]. Considering a Teichm\"uller lift $\Spec(W(k_0))\to\Spec(Q)$ factoring through $\Spec(\hat Q_0)$, we get that $(M\otimes_{W(k)} W(k_0),h(\phi\otimes\sigma_{k_0}))$ is isomorphic to $(M\otimes_{W(k)} W(k_0),\tilde h(\phi\otimes\sigma_{k_0}))$, for some $\tilde h\in G^0(W(\overline{k}))\cap G(W(k_0))$. So $(M\otimes_{W(k)} W(k_0),h(\phi\otimes\sigma_{k_0}))$ does not have both slopes $0$ and $1$ with positive multiplicities.}$ It suffices to prove 2) under the extra assumption that the transcendental degree of $k$ is countable. So there is an $E(G,X)$-monomorphism $B(k)\hookrightarrow\dbC$. Let $w\in\scrN_{E(G,X)}(\dbC)$ be the composite of the resulting morphism $\Spec(\dbC)\to\Spec(B(k))$ with the generic fibre $z_{B(k)}$ of $z$. From 4.2.1 and 1) we get the existence of an isomorphism $(M,(\tilde t_{\alpha})_{\alpha\in\scrJ})\arrowsim  (L^*_{(p)}\otimes_{\dbZ_{(p)}} W(k),(v_{\alpha})_{\alpha\in\scrJ})$. From this 2) follows. This ends the proof.

\medskip\noindent
{\bf 4.2.3. Remark.} We refer to 4.2.2 (*). If $\tilde G^1_{B(k)}$ is a reductive subgroup of $\tilde G_{B(k)}$ through which $\mu_{B(k)}$ factors and  if $\phi$ normalizes the Lie subalgebra $\Lie(\tilde G^1_{B(k)})$ of $\End(M[{1\over p}])=M[{1\over p}]\otimes_{B(k)} M^*[{1\over p}]$, then a family of tensors $(\tilde t_{\alpha})_{\alpha\in\scrJ^1}$ as in 4.2.2 (*) exists (cf. [Va6, 2.5.3]).

\bigskip\noindent
{\bf 4.3. Proof of 1.4.}
We start the proof of 1.4. So $e(v)=1$, i.e. the prime $v$ of $E(G,X)$ is unramified over $p$. Let $O$ be a faithfully flat $O_{(v)}$-algebra which is a DVR of index of ramification 1. Let $Y$ be a regular, formally smooth $O$-scheme such that there is a morphism $q_{Y_{E(G,X)}}:Y_{E(G,X)}\to\Sh_H(G,X)=\scrN_{E(G,X)}$. So $Y$ is a healthy regular scheme, cf. 1.2.1. So $q_{Y_{E(G,X)}}$ extends uniquely to a morphism $q_Y:Y\to\scrN^\prime$, cf. 3.1.1 1). To prove 1.4 1) we just need to show that $q_Y$ factors through $\scrN$, cf. 3.1.3, 3.1.5 and def. 1.2 1). It suffices to check this under the extra assumption that $Y=\Spec(R_Y)$, where $R_Y=W(k)[[T_1,...,T_m]])$ with $m\in\dbN\cup\{0\}$ and with $T_1$, ..., $T_m$ as independent variables. Let $z_Y\in Y(W(k))$ be defined by the $W(k)$-epimorphism $R_Y\twoheadrightarrow W(k)$ taking $T_i$ into $0$, $i\in\{1,...,m\}$. We will use the notations of 4.2 for $z:=q_Y\circ z_Y\in\scrN^\prime(W(k))$. As $\scrN$ is an open subscheme of $\scrN^\prime$ (cf. 3.1.5), to show that $q_Y$ factors through $\scrN$ it suffices to show that $z$ factors through $\scrN$. Let $\tilde G^\prime$ be the universal smoothening of $\tilde G$, cf. 2.3. 

From Fontaine comparison theory we get that $\tilde G_{B(k)}$ is isomorphic to $G_{B(k)}$ and so has dimension $l:=\dim_{\dbQ}(G)$. So the relative dimension of $\tilde G^\prime$ over $W(k)$ is also $l$. Let $R$ be the completion of the local ring of $\tilde G^\prime$ at the identity element of $\tilde G^\prime_k$. We choose an identification $R=W(k)[[x_1,...,x_l]]$ such that the identity section of $\tilde G^\prime$ is defined by $x_1=...=x_l=0$. 

Let $M_{R}:=M\otimes_{W(k)} R$ and $F^1_{R}:=F^1\otimes_{W(k)} R$. Let $\Phi_{R}$ be the Frobenius lift of $R$ compatible with $\sigma$ and taking $x_i$ into $x_i^p$. Let $\Phi:=g_{\text{univ}}(\phi\otimes\Phi_{R})$ be the $\Phi_{R}$-linear endomorphism of $M_{R}$, where $g_{\text{univ}}\in \tilde G^\prime(R)$ is the universal element. Let $\Omega^\wedge_{R/W(k)}$ be the $p$-adic completion of $\Omega_{R/W(k)}$; it is a free $R$-module having $\{dx_1,...,dx_l\}$ as an $R$-basis. Let $\nabla:M_{R}\to M_{R}\otimes_{R} \Omega^\wedge_{R/W(k)}$ be the unique connection such that we have $\nabla\circ\Phi=(\Phi\otimes d\Phi_{R})\circ\nabla$, where $d\Phi_{R}$ is the differential map of $\Phi_{R}$ (cf. [Fa, Th. 10]). 

Loc. cit. implies that $\nabla$ is integrable and nilpotent mod $p$. Moreover we have the following three properties (see [Va6, 3.3 1)] for the first two properties and see [Va6, proof of 5.1] for the third one):

\medskip
{\bf (i)} {\it the connection $\nabla$ is of the form $\delta+\eta$, where $\delta$ is the connection on $M_{R}$ annihilating $M\otimes 1$ and where $\eta\in(\Lie(\tilde G_{B(k)})\cap\End(M))\otimes_{W(k)} \Omega^\wedge_{R/W(k)}$;} 

\smallskip
{\bf (ii)} {\it the connection on $\scrT(M_{R})=\scrT(M)\otimes_{W(k)} R$ induced naturally by $\nabla$ annihilates the tensor $\tilde t_{\alpha}\in\scrT(M)\otimes_{W(k)} R[{1\over p}]$, $\forall\alpha\in\scrJ$;}

\smallskip
{\bf (iii)} {\it the connection $\nabla$ is versal and its Kodaira--Spencer map has as image the direct summand $(\Lie(\tilde G_{B(k)})\cap\Hom(F^1,F^0))\otimes_{W(k)} R$ of $\Hom(F^1,F^0)\otimes_{W(k)} R$.}

\medskip
Each element of $\Ker(\dbG_{mW(k)}(R)\to\dbG_{mW(k)}(R/(x_1,...,x_l)))$ is of the form $\beta\Phi_{R}(\beta^{-1})$ with $\beta\in\Ker(\dbG_{mW(k)}(R)\to\dbG_{mW(k)}(R/(x_1,...,x_l)))$. So there is a unique $\Ker(\dbG_{mW(k)}(R)\to\dbG_{mW(k)}(R/(x_1,...,x_l)))$-multiple $p_{M_{R}}$ of the perfect alternating form $p_M$ on $M_{R}$ such that we have $p_{M_{R}}(\Phi(x),\Phi(y))=p\Phi_{R}(p_{M_{R}}(x,y))$, $\forall x, y\in M_{R}$.

The categories of $p$-divisible groups over $\Spec(R)$ and respectively over $\Spf(R)$ are canonically isomorphic, cf. [dJ1, 2.4.4]; below we will use this fact without any extra comment. Let $(D_{R},p_{D_{R}})$ be the principally quasi-polarized $p$-divisible group over $\Spec(R)$ lifting $(D,p_D)$ and whose principally quasi-polarized filtered $F$-crystal is $(M_{R},F^1_{R},\Phi,\nabla,p_{M_{R}})$. 
The existence up to unique isomorphism of the fibre of $(D_{R},p_{D_{R}})$ over $\Spec(R/pR)$ is implied by [dJ1, Th. of introd.]. So as the the ideal $p(x_1,...,x_l)$ of $R/(x_1,...,x_l)^m$ has a natural nilpotent PD structure $\forall m\in\dbN$, the existence up to unique isomorphism of $(D_{R},p_{D_{R}})$ is implied by Grothendieck--Messing deformation theory. 

Let $(A_{R},p_{A_{R}})$ be the principally polarized abelian scheme over $R$ lifting $(A,p_A)$ and whose principally quasi-polarized $p$-divisible group is $(D_{R},p_{D_{R}})$, cf. Serre--Tate deformation theory and Grothendieck algebraization theorem. Let 
$$q_{R}:\Spec(R)\to\scrM$$ 
be the morphism corresponding to $(A_{R},p_{A_{R}})$ and the level $N$ symplectic similitude structures of $(A_{R},p_{A_{R}})$ lifting those of $(A,p_A)$ ($N\in\dbN$ being prime to $p$).

Let $d:=\dim_{\dbC}(X)$. Let $y:\Spec(k)\hookrightarrow\scrN^\prime_{W(k)}$ be the closed embedding defined naturally by the special fibre of $z$. Let $O_y^{\bigg}$ and $O_{y}$ be the completions of the local rings of $z$ viewed as $W(k)$-valued points of $\scrM_{W(k)}$ and respectively of $\scrN^\prime_{W(k)}$. The local ring $O_{y}$ is normal and of dimension $1+d$, cf. 3.1.1 2) and the fact that any normal $O_{(v)}$-scheme of finite type is excellent. The natural homomorphism $n_y:O_y^{\bigg}\to O_y$ is finite, cf. 3.1.1 3).

The rank $r$ of the direct summand $\Lie(\tilde G_{B(k)})\cap\Hom(F^1,F^0)$ of $\Hom(F^1,F^0)$ is the dimension of the unipotent radical of the maximal parabolic subgroup of $\tilde G_{B(k)}$ normalizing $F^1[{1\over p}]$ and so is the dimension of the compact dual of any connected component of $X$. So we have $r=d$. From (iii) and the equality $r=d$, we get that there is a natural $W(k)$-epimorphism $h_{R}:R\twoheadrightarrow \tilde R:=W(k)[[y_1,...,y_d]]$ such that the composite of the $W(k)$-homomorphism $s_R:O_y^{\bigg}\to R$ defining $q_{R}$ with $h_{R}$, is a $W(k)$-epimorphism $h_y^{\bigg}:O_y^{\bigg}\twoheadrightarrow \tilde R$. In order to show that there is a $W(k)$-homomorphism $H_y:O_y\to\tilde R$ making the following diagram commutative
$$
\spreadmatrixlines{1\jot}
\CD
O_y^{\bigg} @>{n_y}>> O_y\\
@V{s_R}VV @VV{h_y}V\\
R@>{h_R}>>\tilde R,
\endCD
$$
we will need to first recall a result of Faltings.

\medskip\noindent
{\bf 4.3.1. A result of Faltings.} As $\tilde t_{\alpha}\in\scrT(M)[{1\over p}]$ is the de Rham component of the Hodge cycle $\tilde w_{\alpha}$ on $A_{B(k)}$ and due to 4.3 (ii), a result of Faltings implies that $\tilde t_{\alpha}\in\scrT(M)\otimes_{W(k)} R[{1\over p}]$ is the de Rham component of a Hodge cycle on $A_{R[{1\over p}]}$. As the essence of this result is just outlined in [Va1, 4.1.5], we include a full proof of it here. 

As $\scrM$ is a pro-\'etale cover of a quasi-projective $\dbZ_{(p)}$-scheme (cf. 3.1.1 4)) and as $\scrJ$ is a countable set, it suffices to prove Faltings' result in the case when there is a morphism $e_k:\Spec(\dbC)\to\Spec(W(k))$. We view $\dbC$ as a $W(k)$-algebra via $e_k$. Let $\tilde R_{\dbC}:=\dbC[[y_1,...,y_d]]$ and let $R_{\dbC}:=\dbC[[x_1,...,x_l]]$. Let $I_{\dbC}$ be the maximal ideal of $R_{\dbC}$. 

Let $(A_{R_{\dbC}},p_{A_{R_{\dbC}}},(\tilde t_{\alpha})_{\alpha\in\scrJ})$ be the pull back of $(A_{R},p_{A_{R}},(\tilde t_{\alpha})_{\alpha\in\scrJ})$ via the natural $W(k)$-monomorphism $R=W(k)[[x_1,...,x_l]]\hookrightarrow \dbC[[x_1,...,x_l]]=R_{\dbC}$. It suffices to show that $\tilde t_{\alpha}\in \scrT(M)\otimes_{W(k)} R_{\dbC}$ is the de Rham component of a Hodge cycle on $A_{R_{\dbC}}$.

Let $(B_{\tilde R_{\dbC}},p_{B_{\tilde R_{\dbC}}},(w_{\alpha}^{\tilde R_{\dbC}})_{\alpha\in\scrJ})$ be the pull back of $(\scrA^\prime,\scrP_{\scrA^\prime},(w_{\alpha}^{\scrA^\prime})_{\alpha\in\scrJ})$ via a formally \'etale morphism $\Spec(\tilde R_{\dbC})\to\scrN$ whose composite with the natural embedding $\Spec(\dbC)\hookrightarrow\Spec(\tilde R_{\dbC})$ is $e_k\circ z\in\scrN^\prime(dbC)=\scrN(\dbC)$. Let $W_{\tilde R_{\dbC}}:=H^1_{dR}(B_{\tilde R_{\dbC}}/\tilde R_{\dbC})$. Let $\psi_{\tilde R_{\dbC}}$ be the perfect alternating form on $W_{\tilde R_{\dbC}}$ defined by $p_{B_{\tilde R_{\dbC}}}$. Let $t_{\alpha}^{\tilde R_{\dbC}}\in\scrT(W_{\tilde R_{\dbC}})$ be the de Rham component of $w_{\alpha}^{\tilde R_{\dbC}}$. Let $\nabla_B$ be the Gauss--Manin connection on $W_{\tilde R_{\dbC}}$ defined by $B_{\tilde R_{\dbC}}$. 
We denote by $\psi^*$ the alternating form on $W^*$ defined naturally by $\psi$. It is well known that there are isomorphisms $I_B:(W_{\tilde R_{\dbC}},\psi_{\tilde R_{\dbC}},(t_{\alpha}^{\tilde R_{\dbC}})_{\alpha\in\scrJ})\arrowsim (W^*\otimes_{\dbQ} \tilde R_{\dbC},\psi^*,(v_{\alpha})_{\alpha\in\scrJ})$ under which $\nabla_B$ becomes the connection on $W^*\otimes_{\dbQ} \tilde R_{\dbC}$ annihilating $x\otimes 1$, $\forall x\in W^*$ (for instance, see [De3, \S2 and \S6]). We fix such an isomorphism $I_B$ and we view it as an identification. 

By induction on $s\in\dbN$ we show there is a unique morphism of $\dbC$-schemes
$$J_s:\Spec(R_{\dbC}/I_{\dbC}^s)\to\Spec(\tilde R_{\dbC})$$ 
having the following two properties:

\medskip
{\bf (i)} {\it the kernel of the composite $\dbC$-homomorphism $\tilde R_{\dbC}\to R_{\dbC}/I_{\dbC}^s\twoheadrightarrow R_{\dbC}/I_{\dbC}$ is the ideal $(y_1,...,y_d)$ of $\tilde R_{\dbC}$;} 

\smallskip
{\bf (ii)} {\it there is an isomorphism $Q_s$ between $J_s^*((B_{\tilde R_{\dbC}},p_{B_{\tilde R_{\dbC}}},(t_{\alpha}^{\tilde R_{\dbC}})_{\alpha\in\scrJ}))$ and the reduction mod $I_{\dbC}^s$ of $(A_{R_{\dbC}},p_{A_{R_{\dbC}}},(\tilde t_{\alpha})_{\alpha\in\scrJ})$ which modulo $I_{\dbC}/I_{\dbC}^s$ is defined by $1_{A_{\dbC}}$.} 

\medskip
As $\scrN_{E(G,X)}$ is a finite, \'etale scheme over a closed subscheme of $\scrM_{E(G,X)}$, the deformation $(B_{\tilde R_{\dbC}},p_{B_{\tilde R_{\dbC}}})$ of the principally polarized abelian variety $(A,p_A)_{\dbC}$ is versal and the Kodaira--Spencer map of $\nabla_B$ is injective having as image a free $\tilde R_{\dbC}$-module of rank $d$. This implies the uniqueness of $J_s$. The existence of $J_1$ is obvious. The passage from the existence of $J_s$ to the existence of $J_{s+1}$ goes as follows. Let $J_{s+1}^\prime:\Spec(R_{\dbC}/I_{\dbC}^{s+1})\to\Spec(\tilde R_{\dbC})$ be an arbitrary morphism of $\dbC$-schemes lifting $J_s$. Let $\nabla_{B,s+1}$ be the connection on $W_{\tilde R_{\dbC}}\otimes_{\tilde R_{\dbC}} R_{\dbC}/I_{\dbC}^{s+1}=W^*\otimes_{\dbQ} R_{\dbC}/I_{\dbC}^{s+1}$ which is the extension of $\nabla_B$ via $J_{s+1}^\prime$ (the identification being defined by $I_B$). 

We recall that for any $u\in\dbN$ and for every abelian scheme $\pi_C:C\to Z$ over a smooth $\dbC$-scheme $Z$, we have a natural isomorphism of complex sheaves
$$R^u\pi_{C^{\an}*}(\dbC)\arrowsim R^u \pi_{C^{\an}*}(\Omega^{.}_{C^{\an}/Z^{\an}})^{\nabla_C^{\an}}\leqno (1)$$ 
on the complex manifold $Z^{\an}$. Here $\pi_{C^{\an}*}:C^{\an}\to Z^{\an}$ (resp. $\nabla^{\an}_C$) is the morphism of complex manifolds associated naturally to $\pi_C$ (resp. is the connection on $R^u \pi_{C^{\an}*}(\Omega^{.}_{C^{\an}/Z^{\an}})$ induced by the Gauss--Manin connection on $R^u \pi_{C*}(\Omega^{.}_{C/Z})$). 

For any $\beta\in\dbG_{m\dbC}(\dbC)$, there are isomorphisms of $(W^*\otimes_{\dbQ} R_{\dbC},(v_{\alpha})_{\alpha\in\scrJ})$ taking $\psi^*$ into $\beta\psi^*$. So based on the proof of 4.2.1 and the construction of $M_{R}$ we get that there are isomorphisms
$$I_A:(M_{R}\otimes_{R} R_{\dbC},p_{M_{R}},(\tilde t_{\alpha})_{\alpha\in\scrJ})\arrowsim (W^*\otimes_{\dbQ} R_{\dbC},\psi^*,(v_{\alpha})_{\alpha\in\scrJ}).$$
We consider a deformation $\tilde A_Z$ of $A_{R_{\dbC}/I_{\dbC}^{s+1}}$ over a smooth $\dbC$-scheme $Z$. Let $\tilde A_Z^t$ be the dual abelian scheme of $\tilde A_Z$. By applying (1) for $u=1$ and for $C=\tilde A_Z$ and based on 4.3 (ii), we get that we can choose $I_A$ such that under it the Gauss--Manin connection $\nabla_{A,s+1}$ on $M_{R}\otimes_{R} R_{\dbC}/I_{\dbC}^{s+1}$ becomes the connection on $W^*\otimes_{\dbC} R_{\dbC}/I_{\dbC}^{s+1}$ annihilating $x\otimes 1$, $\forall x\in W^*$. We will view the reduction $I_{A,s+1}$ of $I_A$ mod $I_{\dbC}^{s+1}$ as an identification. 

We fix an isomorphism $D_{s+1}$ between $J_{s+1}^{\prime*}((W_{\tilde R_{\dbC}},\psi_{\tilde R_{\dbC}},(t_{\alpha}^{\tilde R_{\dbC}})_{\alpha\in\scrJ}))=(W^*\otimes_{\dbQ} R_{\dbC}/I_{\dbC}^{s+1},\break
\psi^*,(v_{\alpha})_{\alpha\in\scrJ})$ and $(M_{R}\otimes_{R} R_{\dbC}/I_{\dbC}^{s+1},p_{M_{R}},(\tilde t_{\alpha})_{\alpha\in\scrJ})=(W^*\otimes_{\dbQ} R_{\dbC}/I_{\dbC}^{s+1},\psi^*,(v_{\alpha})_{\alpha\in\scrJ})$ with the properties that:

\medskip
{\bf  (iii)} {\it modulo $I_{\dbC}^s$ it is the isomorphism defined by $Q_s$;}

\smallskip
{\bf (iv)} {\it it respects the Gauss--Manin connections, i.e. it takes $\nabla_{B,s+1}$ into $\nabla_{A,s+1}$.}

\medskip
Let $F_{A,s+1}^1$ and $F_{B,s+1}^1$ be the Hodge filtrations of $W^*\otimes_{\dbQ} R_{\dbC}/I_{\dbC}^{s+1}$ defined naturally by $A_{R_{\dbC}}$ and respectively by $B_{\tilde R_{\dbC}}$. The direct summands $F^1_{A,s+1}$ and $D_{s+1}(F^1_{B,s+1})$ of $W^*\otimes_{\dbQ} R_{\dbC}/I_{\dbC}^{s+1}$ coincide modulo $I_{\dbC}^s/I_{\dbC}^{s+1}$, cf. (iii). Moreover, there are direct sum decompositions $W^*\otimes_{\dbQ} R_{\dbC}/I_{\dbC}^{s+1}=F^1_{A,s+1}\oplus F^0_{A,s+1}=F^1_{B,s+1}\oplus F^0_{B,s+1}$ defined naturally by cocharacters $\mu_{A,s+1}$ and respectively $\mu_{B,s+1}$ of the subgroup $G_{R_{\dbC}/I_{\dbC}^{s+1}}$ of $GL(W^*\otimes_{\dbQ} R_{\dbC}/I_{\dbC}^{s+1})$. Argument: the existence of $\mu_{A,s+1}$ is a direct consequence of the existence of the cocharacter $\mu:\dbG_{mW(k)}\to\tilde G$ (see paragraph before 4.2.1) and of the definition of $F^1_{R}$ (see 4.3) while the existence of $\mu_{B,s+1}$ is well known. As $F^1_{A,s+1}$ and $D_{s+1}(F^1_{B,s+1})$ coincide modulo $I_{\dbC}^s/I_{\dbC}^{s+1}$ (cf. (iii)), we can choose $\mu_{A,s+1}$ and $\mu_{B,s+1}$ such that $D_{s+1}^{-1}\mu_{A,s+1}D_{s+1}$ and $\mu_{B,s+1}$ coincide modulo $I_{\dbC}^s/I_{\dbC}^{s+1}$. So based on [SGA3, Vol. II, p. 47--48], there is $g_{s+1}\in\Ker(G(R_{\dbC}/I_{\dbC}^{s+1})\to G(R_{\dbC}/I_{\dbC}^s))$ such that $D_{s+1}^{-1}\mu_{A,s+1}D_{s+1}=g_{s+1}\mu_{B,s+1}g_{s+1}^{-1}$. So $D_{s+1}(g_{s+1}(F^1_{B,s+1}))=F^1_{A,s+1}$.${}^1$ $\vfootnote{1}{The original approach of Faltings was using the strictness of filtrations of morphisms between Hodge $\dbR$-structures to show the existence of the element $g_{s+1}$.}$

As the Kodaira--Spencer map of $\nabla_B$ is injective and its image is a free $\tilde R_{\dbC}$-module of rank $d$ and as $\nabla_{B,s+1}(v_{\alpha})=0$, $\forall\alpha\in\scrJ$, we get that we can replace $J_{s+1}^\prime$ by another morphism $J_{s+1}:\Spec(R_{\dbC}/I_{\dbC}^{s+1})\to\Spec(\tilde R_{\dbC})$ lifting $J_s$ and such that under it and $I_B$ the Hodge filtration $F^1_{B,s+1}$ gets replaced by $g_{s+1}(F^1_{B,s+1})$. So $D_{s+1}$ becomes the de Rham realization of an isomorphism $Q_{s+1}$ between $J_{s+1}^*((B_{\tilde R_{\dbC}},p_{B_{\tilde R_{\dbC}}},(t_{\alpha}^{\tilde R_{\dbC}})_{\alpha\in\scrJ}))$ and the reduction mod $I_{\dbC}^{s+1}$ of $(A_{R_{\dbC}},p_{A_{R_{\dbC}}},(\tilde t_{\alpha})_{\alpha\in\scrJ})$ and which lifts $Q_s$. So $J_{s+1}$ has the desired properties. This ends the induction.  

Let $J_{\infty}:\Spec(R_{\dbC})\to\Spec(\tilde R_{\dbC})$ be the morphism defined by $J_s$'s ($s\in\dbN$). The isomorphism $Q_s$ is uniquely determined by (i) and (ii) and so $Q_{s+1}$ lifts $Q_s$. So we get the existence of an isomorphism 
$$Q_{\infty}:J_{\infty}^*((B_{\tilde R_{\dbC}},p_{B_{\tilde R_{\dbC}}},(t_{\alpha}^{\tilde R_{\dbC}})_{\alpha\in\scrJ}))\arrowsim (A_{R_{\dbC}},p_{A_{R_{\dbC}}},(\tilde t_{\alpha})_{\alpha\in\scrJ})$$ 
which modulo $I_{\dbC}$ is defined by $1_{A_{\dbC}}$. So $\tilde t_{\alpha}\in \scrT(M)\otimes_{W(k)} R_{\dbC}$ is the de Rham component of a Hodge cycle on $A_{R_{\dbC}}$, $\forall\alpha\in\scrJ$. This ends the argument for Faltings' result. 

\medskip\noindent
{\bf 4.3.2. End of the proof of 1.4.} The existence of $Q_{\infty}$ implies that the $W(k)$-epimorphism $h_y^{\bigg}:O_y^{\bigg}\to \tilde R$ (see paragraph before 4.3.1) factors through $O_y$. By reasons of dimensions, the resulting $W(k)$-epimorphism $h_{y}:O_y\twoheadrightarrow \tilde R$ is an isomorphism. So $\scrN^{\prime}_{W(k)}$ is formally smooth at $z$ and so $z$ factors through $\scrN$. So we can also view $y$ as a $k$-valued point of $\scrN_{W(k)}$. So the $O_{(v)}$-scheme $\scrN$ together with the natural action of $G(\dbA_f^{(p)})$ on it is the weak integral canonical model of $\Sh_H(G,X)$ over $O_{(v)}$. So 1.4 1) holds. As $h_y$ is an isomorphism, the natural $W(k)$-homomorphism $O_y^{\bigg}\to O_y$ is also an isomorphism. So the natural $W(k)$-morphism $\scrN_{W(k)}\to\scrM_{W(k)}$ is a formally closed embedding at the $k$-valued point $y$ of $\scrN_{W(k)}$. As the morphism $q_Y$ of the beginning of 4.3 was arbitrary, the role of $z$ is that of an arbitrary $W(k)$-valued of $\scrN^\prime$ (and so based on 1.4 1)) of $\scrN$. So the $W(k)$-morphism $\scrN_{W(k)}\to\scrM_{W(k)}$ is a formally closed embedding at any $k$-valued point of $\scrN_{W(k)}$. So 1.4 2) also holds. This ends the proof of 1.4.

\medskip\noindent
{\bf 4.3.3. Simple properties.} We denote also by $q_{R}$ its factorization through $\scrN$ or $\scrN^\prime$. As $h_y$ is an isomorphism and as we have an epimorphism $h_R:R\twoheadrightarrow\tilde R$, the morphism $q_R:\Spec(R)\to\scrN$ is formally smooth. We have a canonical identification $H^1_{dR}(q_R^*(\scrA^\prime)/R)=M_R$ and under it the perfect form on $M_R$ defined by the principal polarization $q_R^*(\scrP_{\scrA^\prime})$ of $q_R^*(\scrA^\prime)$ is $p_{M_R}$, cf. the definition of $q_R:\Spec(R)\to\scrM$. 

Also the pull back of $w_{\alpha}^{\scrA^\prime}$ via $q_{R[{1\over p}]}$ is a Hodge cycle on $A_{R[{1\over p}]}$ having $\tilde t_{\alpha}\in\scrT(M)\otimes_{W(k)} R[{1\over p}]$ as its de Rham component. This follows either from the existence of $Q_{\infty}$ or from the fact that there is no non-trivial tensor of $\scrT(M)\otimes_{W(k)} (x_1,...,x_l)[{1\over p}]$ fixed by $\Phi$. 

\bigskip\noindent
{\bf 4.4. Lemma.} {\it We recall that we use the notations of 1.3. We assume that one of the following two conditions holds:

\medskip
{\bf (i)} there is a smooth, affine group scheme $G^0_{\dbZ_{(p)}}$ over $\dbZ_{(p)}$ extending $G$, whose special fibre $G^0_{\dbF_p}$ has the same rank as $G$ and such that there is a homomorphism $G^0_{\dbZ_{(p)}}\to G_{\dbZ_{(p)}}$ extending the identity automorphism of $G$;

\smallskip
{\bf (ii)} the group scheme $G_{\dbZ_{(p)}}$ is quasi-reductive for $(G,X,v)$ in the sense of 1.2 3). 

\medskip
Then $e(v)=1$ and the $k(v)$-scheme $\scrN_{k(v)}$ is non-empty.}

\medskip
\proof
We first assume (i) holds. Each torus of $G^0_{\dbF_p}$ lifts to a torus of $G^0_{\dbZ_p}$ (cf. [SGA3, Vol. II, p. 47--48]) and so $G_{\dbZ_p}^0$ has tori of rank equal to the rank of $G$. Let $T^0_{\dbZ_{(p)}}$ be a torus of $G^0_{\dbZ_{(p)}}$ of the same rank as $G$ and such that there is $h\in X$ factoring through $T^0_{\dbR}$. Its existence is implied by [Ha, Lemma 5.5.3]. The pair $(T^0_{\dbQ},\{h\})$ is a Shimura pair. Each prime of $E(T^0_{\dbQ},\{h\})$ dividing $v$ is unramified over $p$ (cf. [Mi3, 4.6 and 4.7]) and so we have $e(v)=1$. The intersection $H^0:=H\cap T^0_{\dbZ_{(p)}}(\dbQ_p)$ is the unique hyperspecial subgroup $T^0_{\dbZ_{(p)}}(\dbZ_p)$ of $T^0_{\dbZ_{(p)}}(\dbQ_p)$. So the integral canonical model $\scrT^0$ of $\Sh_{H^0}(T^0_{\dbQ},\{h\})$ over $E(T^0_{\dbQ},\{h\})_{(p)}\otimes_{\dbZ_{(p)}} O_{(v)}$ exists and has non-trivial fibres of dimension $0$ (cf. [Va1, 3.2.8]). So $\scrT^0$ as an $O_{(v)}$-scheme is formally \'etale. So as $\scrN$ has the smooth extension property (cf. 1.4 1)) and as $\scrT^0$ is formally \'etale over $O_{(v)}$, the functorial morphism $\Sh_{H^0}(T^0_{\dbQ},\{h\})\to \Sh_H(G,X)$ of $E(G,X)$-schemes extends uniquely to a morphism $\scrT^0\to\scrN$ of $O_{(v)}$-schemes. So the $k(v)$-scheme $\scrN_{k(v)}$ is non-empty as $\scrT^0_{k(v)}$ is so. 

We now assume (ii) holds. Let $G^1_{\dbZ_p}$ and $\mu_v$ be as in 1.2 3). Let $T^1_{\dbF_p}$ be a maximal torus of $G^1_{\dbF_p}$. Due to the existence of $\mu_v$, $T^1_{\dbF_p}$ has positive rank. The torus $T^1_{\dbF_p}$ lifts to a torus $T^1_{\dbZ_p}$ of $G^1_{\dbZ_p}$, cf. [SGA3, Vol. II, p. 47--48]. Let $T^0_{0\dbQ_p}$ be a maximal torus of $G_{\dbQ_p}$ having $T^1_{\dbQ_p}$ as a subtorus. Let $T^0$ be a maximal torus of $G$ such that there is $h\in X$ factoring through $T^0_{\dbR}$ and such that $T^0_{\dbQ_p}$ is $G_{\dbZ_p}(\dbZ_p)$-conjugate to $T^0_{0\dbQ_p}$. Again, its existence is implied by [Ha, Lemma 5.5.3]. Not to introduce extra notations, we will assume $T^0_{0\dbQ_p}=T^0_{\dbQ_p}$. 

The intersection $H^0:=H\cap T^0(\dbQ_p)$ is not necessarily the maximal compact, open subgroup of $T^0(\dbQ_p)$ and the subgroup $T^0(\dbQ)H^0$ of $T^0(\dbQ_p)$ is not necessarily $T^0(\dbQ_p)$ itself. However, the intersection $T^1_{\dbZ_p}(\dbQ_p)\cap H$ is the unique hyperspecial subgroup $T^1_{\dbZ_{p}}(\dbZ_p)$ of $T^1_{\dbQ_{p}}(\dbQ_p)$. We fix an $O_{(v)}$-monomorphism $W(k(v))\hookrightarrow\dbC$ as in 1.2 3). As $\mu_h$ and $\mu_v$ are $G(\dbC)$-conjugate and as $G^1_{\dbC}$ is a normal subgroup of $G_{\dbC}$, $\mu_h$ factors through the intersection $T^1_{\dbC}=T^0_{\dbC}\cap G^1_{\dbC}$. So as $T^1_{\dbZ_p}$ splits over $W(\dbF)$, we get that the field of definition $E(T^0_{\dbQ},\{h\})$ of $\mu_h$ is a number subfield of $\dbC$ containing $E(G,X)$ and such that all primes of it dividing $v$ are unramified over $p$. So $e(v)=1$. From class field theory (see [Lan, Th. 4 of p. 220]) and the reciprocity map of [Mi2, p. 163--164] we easily get that each connected component of $\Sh_{H^0}(T^0_{\dbQ},\{h\})_{\dbC}$ is the spectrum of an abelian extension of $E(T^0_{\dbQ},\{h\})$ unramified over any prime of $E(T^0_{\dbQ},\{h\})$ dividing $v$. So as in [Va1, 3.2.8] we argue that there is an integral canonical model $\scrT^0$ of $\Sh_{H^0}(T^0_{\dbQ},\{h\})$ over $E(T^0_{\dbQ},\{h\})_{(p)}\otimes_{\dbZ_{(p)}} O_{(v)}$. As above, this implies that the $k(v)$-scheme $\scrN_{k(v)}$ is non-empty. This ends the proof.

\bigskip\noindent
{\bf 4.5. An application.} We assume $G_{\dbQ_p}$ is unramified. Let $\tilde H$ be a hyperspecial subgroup of $G_{\dbQ_p}(\dbQ_p)$. We show that we can modify the $\dbZ$-lattice $L$ of $W$ and if needed $f$ such that we have $H=\tilde H$. Let $\tilde G_{\dbZ_{(p)}}$ be the reductive group scheme over $\dbZ_{(p)}$ having $G$ as its generic fibre and having $\tilde H$ as its group of $\dbZ_p$-valued points (see beginning of \S3). Let $\tilde L_{(p)}$ be a $\dbZ_{(p)}$-lattice of $W$ such that the monomorphism $G\hookrightarrow GL(W)$ extends to a homomorphism $\tilde G_{\dbZ_{(p)}}\to GL(\tilde L_{(p)})$, cf. [Ja, 10.9 of Part I]. Let $\tilde L$ be the $\dbZ$-lattice of $W$ such that we have $\tilde L[{1\over p}]=L[{1\over p}]$ and $\tilde L\otimes_{\dbZ} \dbZ_{(p)}=\tilde L_{(p)}$. If $\psi$ induces a perfect form on $\tilde L$, then $\tilde L$ is the searched for $\dbZ$-lattice. Argument: as $\tilde H$ is a maximal compact subgroup of $G_{\dbQ_p}(\dbQ_p)$, the monomorphism $\tilde H\hookrightarrow G_{\dbQ_p}(\dbQ_p)\cap GL(\tilde L\otimes_{\dbZ} \dbZ_p)(\dbZ_p)$ is an isomorphism. If $\psi$ does not induces a perfect form on $\tilde L$, then we modify $f$ as follows. 

Let $L_1^\prime:=\tilde L\oplus \tilde L^*$. Let $W_1:=L_1^\prime\otimes_{\dbZ} \dbQ$. Let $\psi_1^\prime$ be a perfect alternating form on $L_1^\prime$ such that the group scheme $SL(\tilde L)$, when viewed naturally as a subgroup of $SL(L_1^\prime)$, is in fact a subgroup of $Sp(L_1^\prime,\psi_1^\prime)$. So $\tilde L$ and $\tilde L^*$ are both maximal isotropic $\dbZ$-lattices of $W_1$ with respect to $\psi_1^\prime$. Let $\tilde G^0_{\dbZ_{(p)}}$ be the Zariski closure in $\tilde G_{\dbZ_{(p)}}$ of the identity component of the subgroup of $G$ fixing $\psi$ (one can check that the subgroup of $G$ fixing $\psi$ is in fact connected). It is a reductive subgroup of $SL(\tilde L)$ and so also of $GSp(L_1^\prime,\psi_1^\prime)$. The subgroup of $GSp(L_1^\prime,\psi_1^\prime)$ generated by $Z(GL(L_1^\prime))$ and $\tilde G^0_{\dbZ_{(p)}}$ is a reductive group scheme (cf. 2.2.5 2)) and it is naturally identified with $\tilde G_{\dbZ_{(p)}}$. We fix an $h\in X$. The Hodge structure on $W_1$ defined by $h$ is of type $\{(-1,0),(0,-1)\}$. So as in the proof of 3.3 we get that there is a perfect alternating form $\psi_1$ on $W_1$ normalized by $G$, inducing a perfect form on $L_1^\prime\otimes_{\dbZ} \dbZ_{(p)}$ and such $2\pi i\psi_1$ is a polarization of the Hodge $\dbQ$--structure on $W_1$ defined by $h$. We get an injective map $f_1:(G,X)\hookrightarrow (GSp(W_1,\psi_1),S_1)$. Let $L_1$ be a $\dbZ$-lattice of $W_1$ such that $\psi_1$ induces a perfect alternating form on $L_1$ and $L_1\otimes_{\dbZ} \dbZ_{(p)}=L_1^\prime\otimes_{\dbZ} \dbZ_{(p)}$. As above we argue that $\tilde H=G_{\dbQ_p}(\dbQ_p)\cap GL(L_1\otimes_{\dbZ} \dbZ_p)(\dbZ_p)$. As $\tilde G_{\dbZ_{(p)}}$ is a reductive group scheme, the condition 4.4 (i) holds in the context of $(f_1,L_1,v)$. So from 1.4 1) and 4.4 we get:

\medskip\noindent
{\bf 4.5.1. Corollary.} {\it Let $(G,X)$ be a Shimura pair of Hodge type. Let $v$ a prime of the reflex field $E(G,X)$ dividing a prime $p$ with the property that the group $G_{\dbQ_p}$ is unramified. Then for any hyperspecial subgroup $\tilde H$ of $G_{\dbQ_p}(\dbQ_p)$, there is a weak integral canonical model $\scrN$ of $\Sh_{\tilde H}(G,X)$ over $O_{(v)}$.}
\medskip
This Corollary can be viewed as a complete solution to Langlands conjecture of [La, p. 411] for Shimura varieties of Hodge type. 

\medskip\noindent
{\bf 4.5.2. A more general form of 4.5.1.} We recall that we use the notations of 1.3. Let $\tilde H$ be a maximal compact, open subgroup of $G_{\dbQ_p}(\dbQ_p)$. Let $\tilde G_{\dbZ_{p}}$ be a smooth, affine group scheme over $\dbZ_{p}$ having $G_{\dbQ_p}$ as its generic fibre and such that $\tilde H$ is its group of $\dbZ_p$-valued points (cf. [Ti2, p. 52]). Let $\tilde G_{\dbZ_{(p)}}$ be the smooth, affine group scheme over $\dbZ_{(p)}$ having $G$ as its generic fibre and whose extension to $\dbZ_p$ is $\tilde G_{\dbZ_{p}}$, cf. [Va1, 3.1.3.1]. As above we argue that there is an injective map $f_1:(G,X)\hookrightarrow (GSp(W_1,\psi_1),S_1)$ and a $\dbZ$-lattice $L_1$ of $W_1$ such that $\psi_1$ induces a perfect alternating form on $L_1$ and the monomorphism $f_1:G\hookrightarrow GSp(W_1,\psi_1)$ extends to a homomorphism $\tilde G_{\dbZ_{(p)}}\to GSp(L_1\otimes_{\dbZ} \dbZ_{(p)},\psi_1)$ which is not necessarily a closed embedding. The maximal property of $\tilde H$ implies $\tilde H=G_{\dbQ_p}(\dbQ_p)\cap GSp(L_1\otimes_{\dbZ} \dbZ_p,\psi_1)(\dbZ_p)$. So if the special fibre $\tilde G_{\dbF_{p}}$ of $\tilde G_{\dbZ_{p}}$ has a torus of the same rank as $G$, then condition 4.4 (i) holds and so from 1.4 1) and 4.4 we again get that $\Sh_{\tilde H}(G,X)$ has a weak integral canonical model over $O_{(v)}$.

\bigskip\noindent
{\bf 4.6. A complement to 4.3.} We use the notations of 4.3 used in the proof of 1.4. Let 
$$\scrC_{\text{univ}}:=(M_{R},F^1_{R},\Phi,\nabla,(\tilde t_{\alpha})_{\alpha\in\scrJ},p_{M_{R}}).$$ 
We recall that $M_{R}:=M\otimes_{W(k)} R$, that $F^1_{R}:=F^1\otimes_{W(k)} R$ and that $\Phi:=g_{\text{univ}}(\phi\otimes\Phi_{R})$, where $g_{\text{univ}}\in\tilde G^\prime(R)$ is the universal element. We also recall that there is a principally quasi-polarized $p$-divisible group $(D,p_D)$ over $W(k)$ whose principally quasi-polarized filtered $F$-crystal over $k$ is $(M,F^1,\phi,p_M)$ (see 4.2), that $R_Y=W(k)[[T_1,...,T_m]]$, that $Y=\Spec(R_Y)$ and that $z_Y\in Y(W(k))$ is defined by the ideal $(T_1,...,T_m)$ of $R_Y$. Let $\Phi_{R_Y}$ be the Frobenius lift of $R_Y$ compatible with $\sigma$ and taking $T_i$ into $T_i^p$, $\forall i\in\{1,...,m\}$.

If $\tilde F^1$ is a direct summand of a free $R_Y$-module $\tilde M$ of finite rank, then we endow $\scrT(\tilde M)$ with the tensor product filtration $(F^i(\scrT(\tilde M)))_{i\in\dbZ}$ defined by the filtrations $(F^i(\tilde M))_{i\in\{0,1\}}$ and $(F^i(\tilde M^*)_{i\in\{-1,0\}}$ of $\tilde M$ and respectively $\tilde M^*$, where $F^1(\tilde M):=\tilde F^1$, $F^0(\tilde M):=\tilde M$, $F^{-1}(\tilde M^*):=\tilde M^*$ and $F^0(\tilde M^*):=\{x\in M^*|x(\tilde F^1)=\{0\}\}$. If $\tilde\Phi$ is a $\Phi_{R_Y}$-linear monomorphism of $\tilde M$, then we let $\tilde \Phi$ act in the natural tensor way on $\scrT(\tilde M)[{1\over p}]$; so if $e\in \tilde M^*$, then $\tilde\Phi(e)\in \tilde M^*[{1\over p}]$ is such that we have $\tilde\Phi(e)(\tilde\Phi(x))=\Phi_{R_Y}(e(x))\in R_Y$, $\forall x\in\tilde M$.

Let $(\tilde M,\tilde F^1,\tilde\Phi,\tilde\nabla,(\tilde t_{\alpha}^Y)_{\alpha\in\scrJ},p_{\tilde M})$ be a principally quasi-polarized filtered $F$-crystal over $\Spec(R_Y/pR_Y)$ endowed with a family of tensors $(\tilde t_{\alpha}^Y)_{\alpha\in\scrJ}$ of $\scrT(\tilde M)[{1\over p}]$ such that the following three axioms hold:

\medskip
{\bf (i)} {\it $\tilde\Phi$ induces an $R_Y$-linear isomorphism $(\tilde M+{1\over p}\tilde F^1)\otimes_{R_Y} {}_{\Phi_{R_Y}} R_Y\arrowsim\tilde M$;}

\smallskip
{\bf (ii)} {\it each $\tilde t_{\alpha}^Y$ is fixed by $\tilde\Phi$, is annihilated by $\tilde\nabla$ and belongs to $F^0(\scrT(\tilde M))[{1\over p}]$;}

\smallskip
{\bf (iii)} {\it modulo $I_Y:=(T_1,...,T_m)$ it is $(M,F^1,\phi,(\tilde t_{\alpha})_{\alpha\in\scrJ},p_M)$.}

\medskip\noindent
{\bf 4.6.1. Theorem.} {\it There is a morphism $i_Y:Y\to\Spec(R)$ of $W(k)$-schemes such that $g_{\text{univ}}\circ i_Y\circ z_Y$ is the identity section of $\tilde G^\prime$ and $i_Y^*(\scrC_{R})$ is isomorphic to $\scrC_Y$ under an isomorphism which modulo $I_Y$ becomes the identity automorphism of $1_M$.}

\medskip
\proof
If $\tilde G$ is smooth, then the Theorem is implied by [Fa, Th. 10 and rm. iii) after it]. We now treat the general case following the proof of [Va6, 5.1]. Let $(D_Y,p_{D_Y})$ be the unique principally quasi-polarized $p$-divisible group over $Y$ lifting $(D,p_D)$ and whose principally quasi-polarized filtered $F$-crystal is $(\tilde M,\tilde F^1,\tilde\Phi,\tilde\nabla,p_{\tilde M})$; its existence and uniqueness is argued in the same way we argued them for $(D_{R},p_{D_{R}})$ of 4.3. 

By induction on $s\in\dbN$ we show that there is a morphism $i_{Y,s}:\Spec(R_Y/I_Y^s)\to\Spec(R)$ of $W(k)$-schemes such that $i_{Y,s}^*((D_{R},p_{D_{R}}))$ is isomorphic to $(D_Y,p_{D_Y})$ modulo $I_Y^s$ under an isomorphism $DIS_s$ having the following two properties:

\medskip
{\bf (i)} {\it it lifts the identity automorphism of $(D,p_D)$;}

\smallskip
{\bf (ii)} {\it it defines an isomorphism $IS_s$ between $i_{Y,s}^*(\scrC_{R})$ and $\scrC_Y$ modulo $I_Y^s$ which modulo $I_Y$ is the identity automorphism of $1_M$.} 

\medskip
As $\Phi_Y(T_i)=T_i^p$, such an isomorphism $IS_s$ is unique. The uniqueness of $DIS_1$ is obvious. If $s=1$ we take $i_{Y,s}$ to be defined by the $W(k)$-epimorphism $R\twoheadrightarrow R/(x_1,...,x_l)=W(k)=R_Y/I_Y$ and we take $DIS_1$ and $IS_1$ to be defined by the identity automorphism of $(D,p_D)$ and respectively by $1_M$. 

The passage from $s$ to $s+1$ goes as follows. We endow the ideal $J_{Y,s}:=I_{Y}^s/I_{Y}^{s+1}$ of $R_Y/I_{Y}^{s+1}$ with the trivial PD structure; so $J_{Y,s}^{[2]}=\{0\}$. Based on Grothendieck--Messing deformation theory, the uniqueness of $DIS_{s+1}$ is implied by the uniqueness of $DIS_s$ and of $IS_{s+1}$. So to end the induction we are left to show that $IS_{s+1}$ and $DIS_{s+1}$ exist. 

Let $\tilde i_{Y,s+1}:\Spec(R_Y/I_Y^{s+1})\to\Spec(R)$ be an arbitrary morphism of $W(k)$-schemes through which $i_{Y,s}$ factors naturally. We write 
$$\tilde i_{Y,s+1}^*(M_{R},F_{R},\Phi,\nabla,p_{M_{R}})=(M\otimes_{W(k)} R_Y/I_Y^{s+1},F^1\otimes_{W(k)} R_Y/I_Y^{s+1},\Phi_{s+1},\nabla_{s+1},p_{M_{R}}).$$ 
Due to existence of $DIS_s$, from Grothendieck--Messing deformation theory we get that there is a direct summand of $M\otimes_{W(k)} R_Y/I_Y^{s+1}$ lifting $F^1\otimes_{W(k)} R_Y/I_Y^s$ and such that the quadruple $(\tilde M,\tilde F,\tilde\Phi,\tilde\nabla,p_{\tilde M})$ modulo $I_Y^{s+1}$ is isomorphic to the quadruple 
$(M\otimes_{W(k)} R_Y/I_Y^{s+1},F^1_{s+1},\Phi_{s+1},\nabla_{s+1},p_{M_{R}})$ under an isomorphism $\tilde{IS}_{s+1}$ lifting the one defined by $IS_s$. As $\Phi_{s+1}(M\otimes_{W(k)} J_{Y,s})=\{0\}$ and as $\tilde{IS}_{s+1}$ lifts $IS_s$, under $\tilde{IS}_{s+1}$ the tensor $\tilde t_{\alpha}^Y$ modulo $I_Y^{s+1}$ becomes $\tilde t_{\alpha}$.

The rest of the inductive argument is entirely as in the last four paragraphs of [Va6, proof of 5.1]. Briefly, let $U^{\text{big}}$ and $U$ be the maximal unipotent, closed subgroups of $GL(M)$ and respectively of $\tilde G$ on which $\mu$ acts via the identity character of $\dbG_{mW(k)}$ (see loc. cit.). Let $u_{s+1}\in \Lie(U^{\text{big}})\otimes_{W(k)}J_{Y,s}$ be the unique element such that we have $(1_{M\otimes_{W(k)} R_Y/I_Y^{s+1}}+u_{s+1})(F^1\otimes_{W(k)} R_Y/I_Y^{s+1})=F^1_{s+1}$. As in loc. cit. we argue that $u_{s+1}\in\Lie(U)\otimes_{W(k)} J_{Y,s}$. The image of the Kodaira--Spencer map of $\nabla$ is the tensorization with $R$ of the direct summand $\Lie(U)$ of $\Hom(F^1,F^0)=\Lie(U^{\text{big}})$, cf. 4.3 (iii). So as in loc. cit. we can replace $\tilde i_{Y,s+1}$ by another morphism $i_{Y,s+1}:\Spec(R_Y/I_Y^{s+1})\to\Spec(R)$ through which $i_{Y,s}$ still factors and for which $F^1_{s+1}$ gets replaced by (i.e. becomes) $F^1\otimes_{W(k)} R_Y/I_Y^{s+1}$. So from Grothendieck--Messing deformation theory we get that $i_{Y,s+1}^*((D_{R},p_{D_{R}}))$ is isomorphic to $(D_Y,p_{D_Y})$ modulo $I_Y^{s+1}$ under an isomorphism $DIS_{s+1}$ which lifts $DIS_s$ and which defines an isomorphism $IS_{s+1}$ between $i_{Y,s+1}^*(\scrC_{R})$ and $\scrC_Y$ modulo $I_Y^{s+1}$. As $DIS_{s+1}$ lifts $DIS_s$, the uniqueness of $I_s$ implies that $IS_{s+1}$ lifts $IS_s$ and so also $IS_1$. This ends the induction. 

We take $i_Y:Y\to\Spec(R)$ such that it lifts $i_{Y,s}$, $\forall s\in\dbN$. From the very definition of $i_{Y,1}$ we get that $g_{\text{univ}}\circ i_Y\circ z_Y$ is the identity section of $\tilde G^\prime$. Moreover, $i_Y^*(\scrC_{R})$ is isomorphic to $\scrC_Y$ under an isomorphism lifting $IS_s$, $\forall s\in\dbN$. This ends the proof.

\bigskip
\noindent
{\boldsectionfont \S5. Proof of 1.5}
\bigskip

We recall from 4.2 that $k$ is an algebraically closed field of characteristic $p$ and that the notations $L$, $L_{(p)}$, $f:(G,X)\hookrightarrow (GSp(W,\psi),S)$, $E(G,X)$, $v$, $k(v)$, $O_{(v)}$, $K_p=GSp(L,\psi)(\dbZ_p)$, $G_{\dbZ_{(p)}}$, $H=G_{\dbZ_{(p)}}(\dbZ_p)$, $\scrM$, $\scrN^\prime$, $\scrN$, $(\scrA^\prime,\scrP_{\scrA^\prime})$ are as in 1.3. For simplicity we will also assume that $k$ is of countable transcendental degree. In the next two Chapters, $\pi_{G_0}$, $(v_{\alpha})_{\alpha\in\scrJ}$, $(w_{\alpha}^{\scrA^\prime})_{\alpha\in\scrJ}$ will be as in 4.2 and to a point $z\in\scrN^\prime(W(k))$ we will associate $(A,p_A,(\tilde w_{\alpha})_{\alpha\in\scrJ})$, $(M,F^1,\phi,p_M,(\tilde t_{\alpha})_{\alpha\in\scrJ})$, $M=F^1\oplus F^0$ and $\mu:\dbG_{mW(k)}\to\tilde G$ as in 4.2. In this Chapter we first prove 1.5 (see 5.1 to 5.4) and then at the end of it we include two refinements of this proof needed in \S6 (see 5.5 and 5.6).

Let $R_0:=W(k)[[x]]$, where $x$ is an independent variable. Let $\Phi_{R_0}$ be the Frobenius lift of $R_0$ compatible with $\sigma$ and taking $x$ into $x^p$.  
 
\bigskip\noindent
{\bf 5.1. Basic notations and facts.} 
We start proving 1.5 by introducing some basic notations and facts. We have $e(v)=1$. We recall that $\scrN$ is an open subscheme of $\scrN^\prime$, cf. 3.1.5. So $\scrN_{k(v)}$ is also an open subscheme of $\scrN^\prime_{k(v)}$. Moreover, the open embedding $\scrN\hookrightarrow\scrN^\prime$ is a pro-\'etale cover of an open embedding between $O_{(v)}$-schemes of finite type (cf. proof of 3.1.2) and the $k(v)$-scheme $\scrN_{k(v)}$ is non-empty (cf. 4.4). So to show that $\scrN_{k(v)}$ is a non-empty, open closed subscheme of $\scrN^\prime_{k(v)}$ we just need to show that for any commutative diagram of the form  
$$
\CD
\Spec(k) @>{}>> \Spec(k[[x]]) @<{}<<
\Spec(k((x))) \\
@VV{y}V @VV{q}V @VV{q_{k((x))}}V \\
\scrN^\prime @<{}<<  \scrN^\prime_{k(v)} @<{}<< \scrN_{k(v)},
\endCD
$$ 
the morphism $y:\Spec(k)\to\scrN^\prime$ factors through the subscheme $\scrN_{k(v)}$ of $\scrN^\prime$.  

We consider the principally quasi-polarized filtered $F$-crystal
$$(M_0,\Phi_0,\nabla_0,p_{M_0})$$ 
over $k[[x]]$ of $q^*((\scrA^\prime,\scrP_{\scrA^\prime})\times_{\scrN^\prime} \scrN^\prime_{k(v)})$. So $M_0$ is a free $R_0$-module of rank equal to $\dim_{\dbQ}(W)$, $\Phi_0$ is a $\Phi_{R_0}$-linear endomorphism of $M_0$ and $\nabla_0$ is an integrable and nilpotent mod $p$ connection on $M_0$ such that we have $\nabla_0\circ\Phi_0=(\Phi_0\otimes d\Phi_{R_0})\circ\nabla_0$. 

Let $O$ be the local ring of ${R_0}$ which is a DVR of mixed characteristic. Let $\hat O$ be the completion of $O$. Let $k_1:=\overline{k((x))}$. We fix a Teichm\"uller lift 
$$\Spec(W(k_1))\to\Spec(R_0)$$ 
with respect to $\Phi_{R_0}$; under it $W(k_1)$ gets naturally the structure of an $O$-algebra and so also of an $\hat O$-algebra. Let $\Phi_{\hat O}$ be the Frobenius lift of $\hat O$ defined by $\Phi_{R_0}$ via localization and completion. 

As $\scrN$ is formally smooth over $O_{(v)}$, there is a lift $\tilde z_1:\Spec(\hat O)\to\scrN$ of the morphism $\Spec(k((x)))\to\scrN$ defined naturally by $q$. Let 
$$(\tilde A_1,p_{\tilde A_1},(\tilde w_{1\alpha})_{\alpha\in\scrJ}):=\tilde z_1^*(\scrA^\prime,\scrP_{\scrA^\prime},(\tilde w_{\alpha})_{\alpha\in\scrJ}).$$ 
Let $\tilde t_{1\alpha}$ be the de Rham realization of $\tilde w_{1\alpha}$. We identify canonically $M_0\otimes_{R_0} \hat O=H^1_{dR}(\tilde A_1/\hat O)$ (cf. [Be, 2.3 of Ch. V]) and so we view $\tilde t_{1\alpha}$ as a tensor of $\scrT(M_0\otimes_{R_0} \hat O)[{1\over p}]$. 

For $\alpha\in\scrJ$ let $n(\alpha)\in\dbN\cup\{0\}$ be such that $v_{\alpha}\in W^{*\otimes n(\alpha)}\otimes_{\dbQ} W^{\otimes n(\alpha)}\subset\scrT(W^*)$. We get the existence of $n_{\alpha}\in\dbN\cup\{0\}$ such that we have
$$p^{n_{\alpha}}\tilde t_{1\alpha}\in (M_0^{\otimes n(\alpha)}\otimes_{R_0} M_0^{*\otimes n_{\alpha}})\otimes_{R_0} \hat O\subset \scrT(M_0\otimes_{R_0} \hat O).$$
\noindent
{\bf 5.1.1. Proposition.} {\it We have $p^{n_{\alpha}}\tilde t_{1\alpha}\in M_0^{\otimes n(\alpha)}\otimes_{R_0} M_0^{*\otimes n(\alpha)}\subset\scrT(M_0)$, $\forall\alpha\in\scrJ$.}

\medskip
\proof
The tensor $p^{n_{\alpha}}\tilde t_{1\alpha}$ is fixed by the natural $\sigma_{k_1}$-linear endomorphism of $\scrT(M_0\otimes_{R_0} B(k_1))$ defined by $\Phi_0$ (see 4.2). So $p^{n_{\alpha}}\tilde t_{1\alpha}$ is also fixed by the natural $\Phi_{\hat O}$-linear endomorphism of $\scrT(M_0\otimes_{R_0} \hat O)[{1\over p}]$ defined by $\Phi_0$. 

The field $k((x))$ has $\{x\}$ as a $p$-basis, i.e. $\{1,x,...,x^{p-1}\}$ is a basis of $k((x))$ over $k((x))^p=k((x^p))$. So the $p$-adic completion of the $\hat O$-module $\Omega_{\hat O/W(k)}$ of relative differentials is naturally isomorphic to $\hat Odx$, cf. [BM, 1.3.1]. 

The de Rham component of $w_{\alpha}^{\scrA^\prime}$ is annihilated by the Gauss--Manin connection of $\scrA^\prime$ (this is a property of Hodge cycles, for instance it follows from [De3, 2.5] applied in the context of a quotient of $\Sh_H(G,X)$ by a small compact, open subgroup of $G(\dbA_f^{(p)})$). So the tensor $p^{n_{\alpha}}\tilde t_{1\alpha}$ is annihilated by the Gauss--Manin connection on $\scrT(H^1_{dR}(\tilde A_1/\hat O))=\scrT(M_0\otimes_{R_0} \hat O)$ of $\tilde A_1$ and so also by the $p$-adic completion of this last connection. In other words, $p^{n_{\alpha}}\tilde t_{1\alpha}$ is annihilated by the connection $\nabla_0:M_0\otimes_{R_0} \hat O\to M_0\otimes_{R_0} \hat Odx$ which is the natural extension of the connection $\nabla_0$ on $M_0$ and so denoted also by $\nabla_0$. 

As the field $k((x))$ has a $p$-basis, any $F$-crystal over $k((x))$ is uniquely determined by its evaluation at the thickening naturally associated to the closed embedding $\Spec(k((x)))\hookrightarrow\Spec(\hat O)$ (cf. [BM, 1.3.3]). So the natural identification $(M_0^{\otimes n(\alpha)}\otimes_{R_0} M_0^{*\otimes n_{\alpha}})\otimes_{R_0} \hat O=\End(M_0^{\otimes n(\alpha)}\otimes_{R_0} \hat O)$ allows us to view naturally $p^{n_{\alpha}}\tilde t_{1\alpha}$ as an endomorphism of the $F$-crystal over $k((x))$ defined by the tensor product of $n(\alpha)$-copies of $(M_0\otimes_{R_0} \hat O,\Phi_0\otimes\Phi_{\hat 0},\nabla_0)$. So from 4.1 we get that $p^{n_{\alpha}}\tilde t_{1\alpha}$ can be viewed naturally as an endomorphism of the $F$-crystal over $k[[x]]$ defined by the tensor product of $n(\alpha)$-copies of $(M_0,\Phi_0,\nabla_0)$. So we have $p^{n_{\alpha}}\tilde t_{1\alpha}\in M_0^{\otimes n(\alpha)}\otimes_{R_0} M_0^{*\otimes n(\alpha)}\subset\scrT(M_0)$. This ends the proof of the Proposition.

\medskip\noindent
{\bf 5.1.2. Some group schemes.} To state the next Proposition we need several notations pertaining to group schemes. Let $G^1_{\dbZ_p}$ be a reductive subgroup of $G_{\dbZ_p}$ as in 1.2 3). Warning: it is not necessarily the pull back to $\Spec(\dbZ_p)$ of a closed subgroup of $G_{\dbZ_{(p)}}$. However, by enlarging $(v_{\alpha})_{\alpha\in\scrJ}$, we can assume that the tensor $\pi_{G^1_{\dbQ_p}}$ of 4.2 is a $\dbQ_p$-linear combination of $v_{\alpha}$'s. Let $\tilde\pi^1$ be the corresponding $\dbQ_p$-linear combination of $\tilde t_{1\alpha}$'s. Let $n^1\in\dbN\cup\{0\}$ be such that $p^{n^1}\tilde\pi^1\in \End(\End(M_0))$, cf. 5.1.1. We consider the product decomposition 
$$G^{1\ad}_{\dbQ_p}=\prod_{i\in I_p} G^{1i\ad}_{\dbQ_p}$$ 
in $\dbQ_p$-simple, adjoint groups. Let $G^{1i\der}_{\dbQ_p}$ be the normal subgroup of $G^{1\der}_{\dbQ_p}$ having $G^{1i\ad}_{\dbQ_p}$ as its adjoint. Similarly  we get that there is $n^{1i}\in\dbN\cup\{0\}$ such that $p^{n^{1i}}\tilde\pi^{1i}\in \End(\End(M_0))$; here $\tilde\pi^{1i}$ corresponds to $\pi_{G^{1i\der}_{\dbQ_p}}$ in the same way $\tilde\pi^1$ corresponded to $\pi_{G^1_{\dbQ_p}}$.

By enlarging $(v_{\alpha})_{\alpha\in\scrJ}$, we can assume that each element of $\End(L_{(p)})$ fixed by $G_{\dbZ_{(p)}}$ is $v_{\alpha_0}$ for some $\alpha_0\in\scrJ$. Let $Z^1_{\dbZ_p}$ be the center of the centralizer of $Z^0(G^1_{\dbZ_p})$ in $GL(L_{(p)}\otimes_{\dbZ_{(p)}} \dbZ_p)$. Let $AL(Z^1_{\dbZ_p})$ be the semisimple $\dbZ_p$-subalgebra of $\End(L_{(p)}\otimes_{\dbZ_{(p)}} \dbZ_p)$ whose elements are the elements of $\Lie(Z^1_{\dbZ_p})$. Each $e\in AL(Z^1_{\dbZ_p})$ is a $\dbZ_p$-linear combination of endomorphisms of $L_{(p)}$ fixed by $G_{\dbZ_{(p)}}$ and so is identified naturally with a $\dbZ_p$-endomorphism $e$ of $\scrA^\prime$. For simplicity  we denote also by $e\in\End(M_0)$ the crystalline realization of the $\dbZ_p$-endomorphism $q^*(e)$ of $q^*(\scrA^\prime\times_{\scrN^\prime} \scrN^\prime_{k(v)})$. 

Let $\eta_0$ be the field of fractions of $R_0$. Let $\tilde G_{0\eta_0}$ be the subgroup of $GL(M_0)_{\eta_0}$ fixing $p^{n_{\alpha}}\tilde t_{1\alpha}$, $\forall\alpha\in\scrJ$; this makes sense due to 5.1.1. Let $\tilde G_{0}$ be the Zariski closure of $\tilde G_{0\eta_0}$ in $GL(M_0)$. Let $\tilde G_{0\eta_0}^1$ be the normal, connected subgroup of $\tilde G_{0\eta_0}$ whose Lie algebra is $\tilde\pi^1(\Lie(\tilde G_{0\eta_0}))$. The uniqueness of $\tilde G_{0\eta_0}^1$ is implied by [Bo, 7.1 of Ch. 1]. So to check the existence of $\tilde G_{0\eta_0}^1$ and to show that $\tilde G_{0\eta_0}^1$ is reductive, it suffices to show that $\tilde G_{0B(k_1)}^1$ exists and is reductive. But $\tilde G_{0B(k_1)}^1$ corresponds via Fontaine comparison theory for $A_{1B(k_1)}$ to $G^1_{\dbQ_p}$ and so it exists and is reductive. Similarly we argue that $\tilde G_{0\eta_0}$ is a reductive group. Let 
$$\tilde G_0^1$$ 
be the Zariski closure of $\tilde G_{0\eta_0}^1$ in $GL(M_0)$.

\bigskip\noindent
{\bf 5.2. Key Proposition.} {\it The closed subscheme $\tilde G_0^1$ of $GL(M_0)$ is a reductive subgroup.}

\medskip
\proof
We first show that the Zariski closure $Z^0(\tilde G_0^1)$ of $Z^0(\tilde G_{0\eta_0}^1)$ in $GL(M_0)$ is a subtorus of $GL(M_0)$. As $R_0$ is strictly henselian, the semisimple $R_0$-subalgebra $AL(Z^1_{\dbZ_p})\otimes_{\dbZ_p} R_0)$ of $\End(M_0)$ is a product of copies of $R_0$ and so is the Lie algebra of a split subtorus $\tilde Z^1_{R_0}$ of $GL(M_0)$. But $Z^0(\tilde G_{0\eta_0}^1)$ is naturally a subtorus of $\tilde Z^1_{\eta_0}$ and so $Z^0(\tilde G_0^1)$ is a subtorus of $\tilde Z^1_{R_0}$ and so also of $GL(M_0)$. Based on 2.2.5, to prove the Proposition we just need to show that the Zariski closure $\tilde G_0^{1\der}$ of $\tilde G_{0\eta_0}^{1\der}$ in $GL(M_0)$ is a semisimple group scheme. Let $V$ be an arbitrary local ring of $R_0$ which is a DVR. Based on 2.2.3 and 2.2.4, to prove the last statement, it suffices to show that $\tilde G^{0\der}_{1V}$ is a semisimple subgroup of $GL(M_0)_V$. We have to consider two Cases: $V=O$ and $V$ is of equal characteristic 0. 

\medskip
{\bf Case 1.} 
Let $V$ be $O$. Let $z_1\in\scrN(W(k_1))\subset\scrN^\prime(W(k_1))$ be the composite of $\tilde z_1\in\scrN(\hat O)$ with the the natural morphism $\Spec(W(k_1))\to\Spec(\hat O)$. Let $(A_1,p_{A_1}):=z_1^*(\scrA^\prime,\scrP_{\scrA^\prime})=(\tilde A_1,p_{\tilde A_1})_{W(k_1)}$. Let $(M_1,F^1_1,\phi_1,p_{M_1})$ be the principally quasi-polarized filtered $F$-crystal over $k_1$ of $(A_1,p_{A_1})$. Let $\tilde G_1$, $\mu_1:\dbG_{mW(k_1)}\to \tilde G_1$, $R_{1}$ and 
$$\scrC_{1\text{univ}}=(M_{1R_{1}},F^1_{1R_{1}},\Phi_1,\nabla_1,(\tilde t_{1\alpha})_{\alpha\in\scrJ},p_{M_{1R_{1}}})$$ 
be the analogues of $\tilde G$, $\mu:\dbG_{mW(k)}\to\tilde G$, $R$ and respectively $\scrC_{\text{univ}}$ of 4.2, 4.3 and 4.6 but obtained working with $z_1\in\scrN(W(k_1))$ instead of some $z\in\scrN^\prime(W(k))$. We can identify naturally $M_1=M_0\otimes_{R_0} W(k_1)$ and so we also view each tensor $\tilde t_{1\alpha}$ as a tensor of $\scrT(M_1)[{1\over p}]$ and we view the group scheme $\tilde G^1_{0W(k_1)}$ as a normal, flat, closed subgroup of the flat, closed subgroup $\tilde G_1=\tilde G_{0W(k_1)}$ of $GL(M_1)$. We consider also the morphism
$$q_{R_{1}}:\Spec(R_{1})\to\scrM$$ 
which is the analogue of $q_{R}$ of 4.3 but obtained working with $z_1\in\scrN(W(k_1))$; it factor through $\scrN$ (see 4.3.2 for the case of $q_R$). 

We first show that if $p=2$, then condition 4.2.2 (*) holds for the quadruple $(M_1,F^1_1,\phi_1,\tilde G_1)$ and for the flat, closed subgroup $\tilde G^1_{0W(k_1)}$ of $\tilde G_1=\tilde G_{0W(k_1)}$.

We fix an arbitrary $O_{(v)}$-monomorphism $W(k_1)\hookrightarrow\dbC$ and an isomorphism $\rho_{1\dbC}:(M_1\otimes_{W(k_1)} \dbC,(\tilde t_{1\alpha})_{\alpha\in\scrJ})\arrowsim (W^*\otimes_{\dbQ} \dbC,(v_{\alpha})_{\alpha\in\scrJ})$  such that $F^1_1\otimes_{W(k_1)} \dbC$ is mapped into the Hodge filtration of $W^*\otimes_{\dbQ} \dbC$ defined by a cocharacter $\mu_h:\dbG_{m\dbC}\to G_{\dbC}$ introduced in 1.1 (see proof of 4.2.1; see also [Va1, 4.1]). From 1.2 3) we get that each such $\mu_h$ factors through $G^1_{0\dbC}$. So there is a cocharacter $\mu_0^{1\prime}$ of $\tilde G^1_{0\dbC}$ which is $\tilde G^1_{0\dbC}(\dbC)$-conjugate to any $\rho_{1\dbC}^{-1}\mu_h\rho_{1\dbC}$ and which acts via the inverse of the identical cocharacter of $\dbG_{m\dbC}$ on $F^1_1\otimes_{W(k_1)} \dbC$ and trivially on $M_1\otimes_{W(k_1)} \dbC/F^1\otimes_{W(k_1)} \dbC$. It is easy to see that $\mu_{1\dbC}$ and $\mu_0^{1\prime}$ are $\tilde G_0(\dbC)$-conjugate (to be compared with [Va1, 5.3.1]). So $\mu_1$ factors through $\tilde G^1_{0W(k_1)}$. 

Let $\tilde G^{1\prime}_{0W(k_1)}$ and $\tilde G^\prime_1$ be the universal smoothenings of $\tilde G^{1}_{0W(k_1)}$ and respectively $\tilde G_1$. We have a natural homomorphism $\tilde G^{1\prime}_{0W(k_1)}\to \tilde G^\prime_1$, cf. the universal property of $\tilde G^\prime_1$. Under this homomorphism, the identity component $\tilde C^1_1$ of $\tilde G^{1\prime}_{0k_1}$ maps into the identity component $\tilde C_1$ of $\tilde G^\prime_{1k_1}$. As $\tilde\pi^1$ is fixed by $\phi_1$, $\Lie(\tilde G^1_{B(k_1)})$ is normalized by $\phi$. Moreover $\tilde G^1_{B(k_1)}$ is a reductive group and $\mu_1$ factors through $\tilde G^1_{0W(k_1)}$. So from 4.2.3 we get that 4.2.2 (*) holds for $(M_1,F^1_1,\phi_1,\tilde G_1)$ and for the smooth subgroup $\tilde G^1_{0W(k_1)}$ of $\tilde G_1$ iff we have: 

\medskip
{\bf (i)} {\it all slopes of $(M_1,\tilde g_1\phi_1)$ are positive, for any element $\tilde g_1\in \tilde G^{1\prime}_{0W(k_1)}(W(k_1))$ which mod $p$ belongs to $\tilde C_1^1(k_1)$.}

\medskip
We show that the assumption that (i) does not hold leads to a contradiction. This assumption implies that there is an open, Zariski dense subscheme $\tilde U_1$ of $\tilde C_1$ such that for any $\tilde g_0\in\tilde G_0^{\prime}(W(k_1))$ lifting a $k_1$-valued point of $\tilde U_1$, the $F$-crystal $(M_1,\tilde g_0\phi_1)$ has slope $0$ with positive multiplicity. So the pull backs of $\scrC_{1\text{univ}}$ via geometric points of $R_{1}/pR_{1}$ mapping into $\tilde U_1$ have slope $0$ with positive multiplicity. The morphism $q_{R_1}:\Spec(R_1)\to\scrM$ factors through $\scrN$ and it is associated naturally to $\scrC_{1\text{uni}}$ (and the special fibre of $z_1$). We conclude that there are geometric points of $\scrN_{k(v)}$ such that the $p$-ranks of the corresponding pull backs of $\scrA^\prime$ are positive. But this contradicts the extra hypothesis of 1.5 for $p=2$. So for $p=2$, the condition 4.2.2 (*) holds in the context of $(M_1,F^1_1,\phi_1,\tilde G_1)$. 

So we can apply 4.2.2 2) to $z_1\in\scrN(W(k_1))$ even for $p=2$. So from the proof of 4.2.2 2) we get that regardless who $p$ is, there is an isomorphism $\rho_1:(M_1,(t_{1\alpha})_{\alpha\in\scrJ})\arrowsim (L_{(p)}^*\otimes_{\dbZ_{(p)}} W(k_1),(v_{\alpha})_{\alpha\in\scrJ})$. Such an isomorphism is unique up to its (left) composite with an element of $G_{\dbZ_{(p)}}(W(k))$. So as $G^1_{\dbZ_p}$ is a normal subgroup of $G_{\dbZ_p}$ and as $\tilde G^1_{0B(k_1)}$ corresponds to $G^1_{\dbQ_p}$ via Fontaine comparison theory for $A_{1B(k_1)}$, we get that the isomorphism $GL(M_1)\arrowsim GL(L_{(p)}^*\otimes_{\dbZ_{(p)}} W(k_1))$ induced naturally by $\rho_1$ takes $\tilde G^{1}_{0W(k_1)}$ into $G^1_{W(k_1)}$. So $\tilde G^{1}_{0W(k_1)}$ is a reductive group scheme. So $\tilde G^{1\der}_{0V}$ is a semisimple group scheme.

\medskip
{\bf Case 2.} Let now $V$ be of equal characteristic $0$. Let $L_V^i$ be the image of $\End(M_0\otimes_{R_0} V)$ through $p^{n^{1i}}\tilde\pi^{1i}$. Let $L_V:=\oplus_{i\in I_p} L_V^i$. The $\tilde\pi^{1i}$'s are commuting projectors of $\End(M_0[{1\over p}])$ whose images are pairwise intersecting just in $\{0\}$. So $L_V$ as a $V$-module is a direct summand of $\End(M_0\otimes_{R_0} V)$ which when tensored with $\eta_0$ is $\Lie(\tilde G^{1\der}_{0\eta_0})$. So $L_V$ is a Lie subalgebra of $\End(M_0\otimes_{R_0} V)$. As $\tilde\pi^{1i}$ is a projector of $\End(M_0\otimes_{R_0} V)$, the restriction $TR^i_V$ to $L^i_V$ of the trace form on $\End(M_0\otimes_{R_0} V)$ is perfect. The Killing form on $L^i_V$ is a non-zero rational multiple of $TR^i_V$ as this is so after tensorization with $\eta_0$. So the Killing form on $L^i_V$ is perfect. So also the Killing form on $L_V$ is perfect. So the tensorizations of $L_V$ with $V$-algebras which are fields are semisimple (cf. [Hu, 5.1]) and so any derivation of them is inner (cf. [Hu, 5.3]). 

Let $\text{Aut}(L_V)$ be the group scheme of Lie automorphisms of $L_V$. It is an affine group scheme over $V$ whose Lie algebra is $L_V$ (due tho the fact that all mentioned derivations are inner). So the identity components of the fibres of $\text{Aut}(L_V)$ are semisimple groups. We easily get that the identity component $\text{Aut}^0(L_V)$ of $\text{Aut}(L_V)$ is an adjoint group scheme whose generic fibre is $\tilde G_{0\eta_0}^{1\ad}$ (argument: an inner automorphism can not specialize to an outer automorphism). So the normalization of $\text{Aut}^0(L_V)$ in $\tilde G_{0\eta_0}^{1\der}$ is a semisimple group scheme $\tilde G_{0V}^{2\der}$ over $V$ extending $\tilde G_{0\eta_0}^{1\der}$. To the Lie monomorphism $L_V\hookrightarrow \End(M_0\otimes_{R_0} V)$ corresponds a homomorphism $h:\tilde G_{0V}^{2\der}\to GL(M_0)_V$ extending the monomorphism $\tilde G_{0\eta_0}^{1\der}=\tilde G_{0\eta_0}^{2\der}\to GL(M_0)_{\eta_0}$ (this is so as the DVR $V$ is of equal characteristic 0). As $\Ker(h_{\eta_0})$ is trivial, from 2.2.2 we get that $h$ is a closed embedding. So $\tilde G_{0V}^{1\der}=\tilde G_{0V}^{2\der}$ is a semisimple group scheme. This ends the argument for the Proposition.

\bigskip\noindent
{\bf 5.3. Applications of 5.2.} Let $F^1_0/pF^1_0$ be the kernel of $\Phi_0$ mod $p$; it is a free module over $k[[x]]=R_0/pR_0$ of rank equal to half the rank of $M_0$. As $\mu_1$ factors through $\tilde G^1_{0W(k_1)}$, the subgroup of $\tilde G^1_{0k_1}$ normalizing $F^1_0/pF^1_0\otimes_{k[[x]]} k_1$ is parabolic and so it is the natural pull back of a parabolic subgroup $\tilde P^1_{0k((x))}$ of $\tilde G^1_{0k((x))}$. The $k[[x]]$-scheme of parabolic subgroups of $\tilde G^1_{0k[[x]]}$ is a projective $k[[x]]$-scheme, cf. [SGA3, Vol. III, Cor. 3.5 of p. 445].  So the Zariski closure $\tilde P^1_{0k[[x]]}$ of $\tilde P^1_{0k((x))}$ in $\tilde G^1_{0k[[x]]}$ is a parabolic subgroup of $\tilde G^1_{0k[[x]]}$. So as $\tilde G^1_0$ is a split reductive group scheme and as $\mu_{1k_1}$ factors through $\tilde G^1_{0k_1}$, there is a cocharacter $\mu_{0k[[x]]}:\dbG_{mk[[x]]}\to\tilde G^1_{0k[[x]]}$ factoring through $\tilde P^1_{0k[[x]]}$ and producing a direct sum decomposition $M_0/pM_0=F^1_0/pF^1_0\oplus F^0_0/pF^0_0$ (with $\beta\in\dbG_{mk[[x]]}(k[[x]])$ acting through $\mu_{0k[[x]]}$ on $F^i_0/pF^i_0$ via the multiplication with $\beta^{-i}$, $i\in\{0,1\}$). We choose a cocharacter
$$\mu_0:\dbG_{mR_0}\to \tilde G^1_0$$ 
lifting $\mu_{0k[[x]]}$, cf. [SGA3, Vol. II, p. 47--48]. Let $M_0=F^1_0\oplus F^0_0$ be the direct sum decomposition such that $\beta\in\dbG_{mR_0}(R_0)$ acts through $\mu_0$ on $F^i_0$ via the multiplication with $\beta^{-i}$, $i\in\{0,1\}$; the notations match, i.e. the reduction mod $p$ of $F^i_0$ is as defined above. 

Let 
$$(M,F^1,\phi,\tilde G,(\tilde t_{\alpha})_{\alpha\in\scrJ},p_{M}):=(M_0,F^1_0,\Phi_0,\nabla_0,\tilde G_0,(\tilde t_{1\alpha})_{\alpha\in\scrJ},p_{M_0})\otimes_{R_0} W(k),$$
where the $W(k)$-epimorphism $R_0\twoheadrightarrow W(k)$ has $(x)$ as its kernel. If $p=2$, then all slopes of $(M_1,\phi_1)$ are positive (cf. the extra hypothesis of 1.5 for $p=2$) and so via specialization we get that all slopes of $(M,\phi)$ are positive. So either $p>2$ or $p=2$ and moreover all slopes of $(M,\phi)$ are positive. So there is a unique $p$-divisible group $D$ over $W(k)$ whose filtered $F$-crystal over $k$ is $(M,F^1,\phi)$, cf. [Va6, 2.2.3]. To $p_M$ corresponds a principal quasi-polarization $p_D$ of $D$. So the sextuple 
$$\scrC_0:=(M_0,F^1_0,\Phi_0,\nabla_0,(\tilde t_{1\alpha})_{\alpha\in\scrJ},p_{M_0})$$ 
is a principally quasi-polarized filtered $F$-crystal over $k[[x]]$ endowed with a family of crystalline tensors for which the axioms 4.6 (i) to (iii) hold. 

\bigskip\noindent
{\bf 5.4. End of the proof of 1.5.}
Let $$(A_k,p_{A_k}):=y^*(\scrA^\prime,\scrP_{\scrA^\prime}).$$ 
As in 4.3, starting form $(D,p_D)$ and $(M_0,F^1_0,\Phi_0,\nabla_0,p_{M_0})$ we construct a principally quasi-polarized $p$-divisible group $(D_{R_0},p_{D_{R_0}})$ over $R_0$ and a  homomorphism 
$$q_{R_0}:\Spec(R_0)\to\scrM$$ 
lifting the composite of $y$ withe the morphism $\scrN\to\scrM$ and such that the principally quasi-polarized $p$-divisible group of $q_{R_0}^*(\scrA^\prime,\scrP_{\scrA^\prime})$ is $(D_{R_0},p_{D_{R_0}})$. We consider the composite
$$z_2:\Spec(W(k_1))\to\scrM$$ 
of the dominant morphism $\Spec(W(k_1))\to\Spec(R_0)$ of 5.1 with $q_{R_0}$. We will next check that $q_{R_0}$ factors through $\scrN^\prime$. For this, it suffices to check that $z_2$ factors through $\scrN^\prime$. 

Let $(A_2,p_{A_2})$ be the pull back via $z_2$ of the universal principally polarized abelian scheme over $\scrM$. The principally quasi-polarized filtered $F$-crystal of $(A_2,p_{A_2})$ is canonically identified with $(M_1,F^1_2,\phi_1,p_{M_1})$, where $F^1_2$ is a direct summand of $M_1$ of rank equal to half the rank of $M_1$. Let $(F^i_2(\scrT(M_1))_{i\in\dbZ}$ be the tensor product filtration of $\scrT(M_1)$ defined by $F^1_2$ (see 4.6).  For $\alpha\in\scrJ$ the tensor $\tilde t_{1\alpha}\in\scrT(M_0)[{1\over p}]$ is annihilated by $\nabla_0$, is fixed by $\Phi_0$ and belongs to $F^0_0(\scrT(M_0))[{1\over p}]$, where $(F^i_0(\scrT(M_0))_{i\in\dbZ}$ is the tensor product filtration of $\scrT(M_0)$ defined by $F^1_0$. This implies that $\tilde t_{1\alpha}\in F^0_2(\scrT(M_1))[{1\over p}]$, $\forall\alpha\in\scrJ$. So as before 4.2.1 we argue that the canonical split cocharacter of $(M_1,F^1_2,\phi_1)$ defined in [Wi, p. 512] factors through $\tilde G_{0W(k_1)}$; let $\mu_2:\dbG_{mW(k_1)}\to \tilde G_{0W(k_1)}$ be the resulting factorization. As for $\mu_1$ we argue that we can also view $\mu_2$ as a cocharacter $\mu_2:\dbG_m\to\tilde G_{0W(k_1)}^1$. 

If $G^1_{\dbZ_p}=G_{\dbZ_p}$, then $\tilde G_0=\tilde G_0^1$ is a reductive group scheme over $R_0$ (cf. 5.2) and so $\tilde G$ itself is a reductive group scheme over $W(k)$. So the fact that $z_2$ factors through $\scrN^\prime$ is proved in [Va5, 6.4.1 i) and ii)]. If $G^1_{\dbZ_p}\neq G_{\dbZ_p}$, then the proof of loc. cit. applies entirely provided we work with the triple $(M_1,\phi_1,\tilde G^1_{0W(k_1)})$ which in [Va5] is called a Shimura $F$-crystal over $k_1$. This is so due to the following two reasons. First $\mu_1$ and $\mu_2$ are conjugate under an element of $\tilde G^1_{0W(k_1)}(W(k_1))$ which mod $p$ normalizes the kernel $F^1_1/pF^1_1=F^1_2/pF^1_2$ of $\phi_1$ mod $p$, cf. [Va5, 3.1.2]. Second the fact that the image of the Kodaira--Spencer map of the connection $\nabla_1$ (of Case 1 of 5.2) can be identified with the maximal direct summand of $\Lie(\tilde G_{0W(k_1)})$ and so also of $\Lie(\tilde G_{0W(k_1)}^1)$ on which $\mu_1$ acts via the identity character of $\dbG_{mW(k_1)}$ (see 4.3 (iii) for the case of $\nabla$), applies entirely as in [Va5, proof of 6.4.1] to show that there is a morphism $z_3:\Spec(W(k_1))\to\Spec(R_{1})$ such that the pull back of $\scrC_{1\text{univ}}$ through it is $(M_1,F^1_2,\phi_1,(\tilde t_{1\alpha})_{\alpha\in\scrJ},p_{M_1})$. We denote also by $z_3\in\scrM(W(k_1))$ the composite of $z_3:\Spec(W(k_1))\to \Spec(R_{1})$ with $q_{R_{1}}:\Spec(R_{1})\to\scrM$ (of Case 1 of 5.2); it factors through $\scrN$ (as $q_{R_{1}}$ does) and so we also view $z_3$ as a $W(k_1)$-valued point of $\scrN$ or $\scrN^\prime$. As for $p=2$ the special fibre of $A_1$ has $p$-rank $0$, for $p$ arbitrary we get that the principally quasi-polarized $p$-divisible groups of the pull backs of $(\scrA^\prime,\scrP_{\scrA^\prime})$ via $z_2$ and $z_3$ are the same (cf. [Va6, 2.2.3]). So from Serre--Tate deformation theory we get that $z_2$ and $z_3$ are the same $W(k_1)$-valued points of $\scrM$ and so also of $\scrN^\prime$. 

So $z_2$ factors through $\scrN^\prime$. So $q_{R_0}$ factors through $\scrN^\prime$ and so (cf. 1.4 1)) also through $\scrN$. So the morphisms $q:\Spec(k[[x]])\to\scrN_{k(v)}^\prime$ and $y:\Spec(k)\to\scrN^\prime$ factor through $\scrN_{k(v)}$. This ends the argument that $\scrN_{k(v)}$ is a non-empty, open closed subscheme of $\scrN_{k(v)}^\prime$ and so ends the proof of 1.5.

\bigskip\noindent
{\bf 5.5. Lemma.} {\it We use the above notations of 5.1.2, of Case 1 of 5.2 and of 5.3. We do not assume that the hypotheses of 1.5 hold but we just assume that $e(v)=1$ and that $\tilde G_0$ is a reductive subgroup of $GL(M_0)$. Then there is an isomorphism 
$$\rho_{1z}:(M_1,(\tilde t_{1\alpha})_{\alpha\in\scrJ},p_{M_1})\arrowsim (M\otimes_{W(k)} W(k_1),(\tilde t_{\alpha})_{\alpha\in\scrJ},p_M).$$}
{\it Proof:} Let $\scrC_{\text{univ}}$ be constructed as in 4.6 starting from $(M,F^1,\phi,(\tilde t_{\alpha})_{\alpha\in\scrJ},p_M)$. Its construction makes sense even if we do not know if the composite $z$ of $q_{R_0}:\Spec(R_0)\to\scrM$ with the natural $W(k)$-epimorphism $\Spec(W(k))\hookrightarrow\Spec(R_0)$ defined by the locus of $x=0$, factors or not through $\scrN^\prime$. From 4.6.1 we get that $\scrC_0$ of 5.3 is the pull back of $\scrC_{\text{univ}}$ via a morphism $i_{\Spec(R_0)}:\Spec(R_0)\to\Spec(R)=\Spec(W(k)[[x_1,...,x_l]])$ of $W(k)$-schemes which at the level of $W(k)$-algebras maps the ideal $(x_1,...,x_l)$ of $R$ into the ideal $(x)$ of $R_0$. The existence of $i_{\Spec(R_0)}$ implies the existence of $\rho_{1z}$. This ends the proof.

\bigskip\noindent
{\bf 5.6. A refinement.} In this section we assume $p=2$. The extra hypothesis of 1.5 that the pull backs of $\scrA^\prime$ via geometric points of $\scrN_{k(v)}$ have all slopes positive is too restrictive; it has been used in Case 1 of 5.2 as well as in 5.4. So we now explain how we can considerably weaken it under some extra hypotheses. Let $\tilde v_{1\alpha}$ be the $p$-component of the \'etale component of the Hodge cycle $\tilde w_{1\alpha}$ on $\tilde A_{1B(k_1)}=A_{1B(k_1)}$ (cf. notations before 5.1.1). We list three additional conditions:

\medskip
{\bf (i)} {\it the group scheme $G_{\dbZ_{(p)}}$ is reductive and we have $G^1_{\dbZ_p}=G_{\dbZ_p}$;}

\smallskip
{\bf (ii)} {\it there is an isomorphism $(M_1,(\tilde t_{1\alpha})_{\alpha\in\scrJ})\arrowsim (H^1_{\acute et}(A_{1B(k_1)},\dbZ_p)\otimes_{\dbZ_p} W(k_1),(\tilde v_{1\alpha})_{\alpha\in\scrJ})$;}

\smallskip
{\bf (iii)} {\it the quadruple $(M,F^1,\phi,p_M)$ is the principally quasi-polarized filtered $F$-crystal of a principally polarized abelian scheme $(A,p_A)$ over $W(k)$ lifting $(A_k,p_{A_k})$ and such that there is an isomorphism $(H^1_{\acute et}(A_{B(k)},\dbZ_p)\otimes_{\dbZ_p} W(k),(\tilde v_{\alpha})_{\alpha\in\scrJ})\arrowsim (M,(\tilde t_{\alpha})_{\alpha\in\scrJ})$; here $\tilde v_{\alpha}\in\scrT(H^1_{\acute et}(A_{B(k)},\dbZ_p))[{1\over p}]$ corresponds to $\tilde t_{\alpha}$ via Fontaine's comparison theory of $A_{B(k)}$.}

\medskip\noindent
{\bf 5.6.1. Theorem.} {\it We assume that $p=2$ and that the above three conditions 5.6 (i) to (iii) hold. Warning: we do not assume that the pull backs of $\scrA^\prime$ via geometric points of $\scrN_{k(v)}$ have all slopes positive. Then with the notations of 5.1, the morphism $y:\Spec(k)\to\scrN^\prime$ factors through $\scrN_{k(v)}$.}

\medskip
\proof
Due to 5.6 (i) and (ii), the proof of 5.2 applies entirely to give us that the Zariski closure $\tilde G_0=\tilde G_0^1$ of $\tilde G_{0\eta_0}=\tilde G_{0\eta_0}^1$ in $GL(M_0)$ is a reductive subgroup of $GL(M_0)$. So an isomorphism $\rho_{1z}$ as in 5.5 exists and the constructions of $\mu_0$, $q_{R_0}$ and $z_2$ can be performed as in 5.3 and 5.4, provided we take $(D,p_D)$ to be the principally quasi-polarized $p$-divisible group of the principally polarized abelian scheme $(A,p_A)$ over $W(k)$ mentioned in 5.6 (iii). A $W(k_1)$-valued point $z_3$ of $\scrN$ as in 5.4 always exists, cf. [Va5, proof 6.4.1 ii)]. But if the special fibre of $A_1$ has positive $p$-rank, then (as $p=2$) we can not conclude that $z_3$ and $z_2$ are the same $W(k_1)$-valued points of $\scrM$. However, as the filtered $F$-crystals over $k_1$ associated to $A_2$ and $A_3$ coincide to $(M_1,F^1_2,\phi_1)$, from Fontaine comparison theory we get that we can identify canonically $(H^1_{\acute et}(A_{2B(k_1)},\dbQ_p),p_{A_2})=(H^1_{\acute et}(A_{3B(k_1)},\dbQ_p),p_{A_3})$; here we denote also by $p_{A_2}$ and $p_{A_3}$ the perfect bilinear forms on $(H^1_{\acute et}(A_{2B(k_1)},\dbQ_p),p_{A_2})=(H^1_{\acute et}(A_{3B(k_1)},\dbQ_p),p_{A_3})$ defined by $p_{A_2}$ and $p_{A_3}$. Under this identification we have natural inclusions $pH^1_{\acute et}(A_{2B(k_1)},\dbZ_p)\subset H^1_{\acute et}(A_{3B(k_1)},\dbZ_p)\subset {1\over p}H^1_{\acute et}(A_{2B(k_1)},\dbZ_p$ of $\Gal(B(k_1))$-modules (cf. [Va6, second paragraph of the proof of 2.2.3] applied to the $p$-divisible groups of $A_2$ and $A_3$; loc. cit. holds with no restrictions on these two $p$-divisible groups). So the multiplication by $p^2$ endomorphism of $A_2$ factors as a composite homomorphism $A_2\to A_3\to A_2$ and in fact the principally polarized abelian schemes $(A_2,p_{A_2})$ and $(A_3,p_{A_3}):=z_3^*(\scrA^\prime,\scrP_{\scrA^\prime})$ are naturally $\dbZ[{1\over p}]$-isomorphic. So we get that $\tilde t_{1\alpha}$ is the de Rham component of a Hodge cycle $\tilde w_{2\alpha}$ on $A_2$. This implies that $\tilde t_{1\alpha}$ is also the de Rham component of a Hodge cycle $\tilde w_{0\alpha}$ on $q_{R_0}^*(\scrA^\prime)\times_{R_0} R_0[{1\over p}]$. 

Based on 5.6 (ii) and (iii) and due to the existence of $\rho_{1z}$, there is an isomorphism
between $H^1_{\acute et}(A_{3B(k_1)},\dbZ_p)\otimes_{\dbZ_p} W(k_1)=H^1_{\acute et}(A_{B(k_1)},\dbZ_p)\otimes_{\dbZ_p} W(k_1)$ and $H^1_{\acute et}(A_{2B(k_1)},\dbZ_p)\otimes_{\dbZ_p} W(k_1)=H^1_{\acute et}(A_{1B(k_1)},\dbZ_p)\otimes_{\dbZ_p} W(k_1)$
which takes the $p$-component $\tilde v_{3\alpha}$ of the \'etale component of $\tilde w_{3\alpha}:=z_3^*(w_{\alpha}^{\scrA^\prime})$ into the $p$-component $\tilde v_{2\alpha}$ of the \'etale component of $\tilde w_{2\alpha}$. So as $G_{\dbF_p}$ is a reductive group (cf. 5.6 (i)) and so connected, from Lang's theorem for affine, connected groups over finite fields we get that any torsor of $G_{\dbZ_p}$ is trivial. 
So there is an isomorphism
$$(H^1_{\acute et}(A_{3B(k_1)},\dbZ_p),(\tilde v_{3\alpha})_{\alpha\in\scrJ})\arrowsim (H^1_{\acute et}(A_{2B(k_1)},\dbZ_p),(\tilde v_{2\alpha})_{\alpha\in\scrJ}).\leqno (2)$$
\indent
Let $W(k_1)\hookrightarrow\dbC$ be as in Case 1 of 5.2. Let $L_2:=H_1(A_{2\dbC},\dbZ)$ and $L_3:=H_1(A_{3\dbC},\dbZ)$. We identify $W=L_3\otimes_{\dbZ} \dbQ=L_2\otimes_{\dbZ} \dbQ$ in such a way that a $\dbG_{m\dbQ}(\dbQ)$-multiple of $\psi$ and $v_{\alpha}$ are the Betti realizations of the principal polarizations of $A_{3\dbC}$ and $A_{2\dbC}$ and respectively of the Hodge cycles $\tilde w_{3\alpha}$ and $\tilde w_{2\alpha}$, $\forall\alpha\in\scrJ$ (cf. proof of 4.2.1). This identification is unique up to an element of $G(\dbQ)$. We have $G(\dbQ_p)=G(\dbQ)H$ (cf. [Mi2, 4.9]) and $L_3\otimes_{\dbZ} \dbZ_p$ is $G_{\dbQ_p}(\dbQ_p)$-conjugate to $L\otimes_{\dbZ} \dbZ_p=L_{(p)}\otimes_{\dbZ_{(p)}} \dbZ_p$ (cf. proof of 4.2.1). From this last sentence and (2) we get the existence of $g\in G(\dbQ)$ such that $g(L_3\otimes_{\dbZ} \dbZ_{(p)})=L_2\otimes_{\dbZ} \dbZ_{(p)}$. 

So as we have $L_3[{1\over p}]=L_2[{1\over p}]$ (cf. the existence of the mentioned $\dbZ[{1\over p}]$-isomorphism) we get that the $B(k_1)$-valued points of $\scrM$ defined naturally by $z_3$ and $z_2$ are in the same orbit under the action of the element $g^{-1}$ of $G(\dbA_f^{(p)})$ on $\scrM$ (see [Va1, 4.1]). So $z_2$ factors through $\scrN^\prime$ as $z_3$ does. So $z_2$ factors also through $\scrN$, cf. 1.4 1). So $q_{R_0}$ factors through $\scrN^\prime$ and so as in the end of 5.4 we get that $y$ factors through $\scrN_{k(v)}$. This ends the proof. 

\bigskip
\noindent
{\boldsectionfont \S6. Proof of 1.5 and complements}
\bigskip

In this Chapter we first prove 1.6 (see 6.1 and 6.2) and then we apply 1.6 2) to provide new examples of N\'eron models (see 6.3). We end up with few extra complements on integral models (see 6.4). As in \S5, $k$ is an algebraically closed field of countable transcendental degree and of characteristic $p$ and the following notations $L$, $L_{(p)}$, $f:(G,X)\hookrightarrow (GSp(W,\psi),S)$, $E(G,X)$, $v$, $k(v)$, $O_{(v)}$, $K_p=GSp(L,\psi)(\dbZ_p)$, $G_{\dbZ_{(p)}}$, $H=G_{\dbZ_{(p)}}(\dbZ_p)$, $\scrM$, $\scrN^\prime$, $\scrN$, $(\scrA^\prime,\scrP_{\scrA^\prime})$ are as in 1.3.

We start the proof of 1.6 by first showing that 1.6 2) implies 1.6 1). So we assume that $G_{\dbZ_{(p)}}$ is a reductive group scheme and that one of the conditions (i) to (iii) of 1.6 1) holds. As in the proof of 3.2.1 we argue that $e(v)=1$. So from 1.6 2) we get that $\Sh_{H}(G,X)$ has an integral canonical model over $O_{(v)}$ and (as $e(v)=1$) from 3.1.4 we get that this model is $\scrN^\prime$ itself. So $\scrN^\prime$ is formally smooth and so we have $\scrN=\scrN^\prime$. So 1.6 1) holds. 

So to prove 1.6 it suffices to prove 1.6 2). So let $(G_1,X_1,H_1)$ and $O(G_1,X_1,p)$ be as in 1.6 2). So $O(G_1,X_1,p)$ is an \'etale $\dbZ_{(p)}$-algebra. Let $G_{1\dbZ_{(p)}}$ be as in the beginning of \S3. A prime $v_1$ of $E(G_1,X_1)$ dividing $p$ is $p$-compact for $(G_1,X_1)$ iff the prime $v^{\ad}_1$ of $E(G_1^{\ad},X_1^{\ad})$ divided by $v_1$ is $p$-compact for $(G_1^{\ad},X_1^{\ad})$. We will prove 1.6 2) in the following two main Steps. 

\bigskip\noindent
{\bf 6.1. Step 1.} In this Step 1 we assume that $(G_1,X_1)$ is a simple, adjoint Shimura pair. We will apply 3.3 to show that we can choose an injective map $f:(G,X)\hookrightarrow (GSp(W,\psi),S)$ as in 3.3 and a prime $v$ of $E(G,X)$ dividing a prime $v_1$ of $O(G_1,X_1,p)$ such that $\scrN=\scrN^{\prime}$ is the integral canonical model of $\Sh_{H}(G,X)$ over $O_{(v)}$. If either $(G_1,X_1)$ is of $A_n$ type or $(G_1,X_1)$ has no compact factors and moreover it is of $C_n$ or $D_n^{\dbH}$ type, then this follows from 3.3 (iv) and 3.2.1. If $p\Ge 5$, this follows from [Va1, 6.4.1] and 3.1.4 regardless of the possible choices we could make in the proof of 3.3. 

So we now consider the case when $p$ is arbitrary and $(G_1,X_1)$ has compact factors and it is of $B_n$, $C_n$, $D_n^{\dbH}$ or $D_n^{\dbR}$ type. If $p=2$ and $(G_1,X_1)$ is of $B_n$, $D_n^{\dbH}$ or $D_n^{\dbR}$ type, then we also assume that $v_1$ is $p$-compact for $(G_1,X_1)$. As $(G_1,X_1)$ is a simple, adjoint Shimura pair having compact factors, the $\dbQ$--rank of $G^{\ad}=G_1$ is $0$ and so based on [BHC] we get that $\Sh(G,X)$ is a pro-\'etale cover of a smooth, projective $E(G,X)$-scheme. Let $N\in\dbN\setminus\{1,2\}$ be relatively prime to $p$. Let $K(N)$ be the compact, open subgroup of $GSp(L,\psi)(\widehat\dbZ)$ acting trivially on $L/NL$. Let $H_0$ be a compact, open subgroup of $G(\dbA_f^{(p)})\cap K(N)$ such that $\scrN^\prime$ is a pro-\'etale cover of $\scrN^\prime/H_0$, cf. 3.1.1 1). We consider a smooth, projective, toroidal compactification $\scrM^{\stc}/K(N)$ of $\scrM/K(N)$ such that the abelian scheme defined by the universal principally abelian scheme over $\scrM/K(N)$ extends to a semiabelian scheme $\scrS(N)$ over $\scrM^{\stc}/K(N)$ (cf. [FC, 6.7 of p. 120]). The fibres of $\scrS(N)$ over points of the complement of  $\scrM/K(N)$ in $\scrM^{\stc}/K(N)$ are semiabelian scheme which are not abelian schemes. 

We show that we can choose $f:(G,X)\hookrightarrow (GSp(W,\psi),S)$ in 3.3 such that the following two things hold:

\medskip
{\bf (a)}  {\it the $k(v)$-scheme $\scrN_{k(v)}/H_0$ is an open closed subscheme of $\scrN^\prime_{k(v)}/H_0$;}

\smallskip
{\bf (b)} {\it the $k(v)$-scheme $\scrN_{k(v)}/H_0$ is proper.} 

\medskip
Before proving these two properties we show how they imply that $\scrN=\scrN^\prime$ and that $\scrN/H_0$ is a smooth, projective $O_{(v)}$-scheme. The connected components of $\scrN^\prime_{B(\dbF)}$ are permuted transitively by $G(\dbA_f^{(p)})$, cf. [Va1, 3.3.2]. As $\scrN$ is $G(\dbA_f^{(p)})$-invariant (see 3.1.2) and as its special fibre is non-empty (see 4.4), from the previous sentence we get that each connected component $\scrC$ of the $W(\dbF)$-scheme $\scrN^\prime_{W(\dbF)}/H_0$ has smooth points of characteristic $p$ and so $\scrC_{B(\dbF)}$ is an absolutely irreducible $B(\dbF)$-variety. Let $\scrC^{\stc}$ be the normalization in the field of fractions of $\scrC$ of the pull back to $\Spec(W(\dbF))$ of $\scrM^{\stc}/K(N)$. It is a normal, integral, projective $W(\dbF)$-scheme having $\scrC$ as an open, Zariski dense subscheme. From (a) and (b) we get that either $\scrC$ is smooth and equal to $\scrC^{\stc}$ or the special fibre of $\scrC^{\stc}$ is not connected. As the theorem on formal functions implies that this last thing can not hold (see [Har, Cor. 11.3 of p. 279]), we get that $\scrC$ is a smooth, projective $W(\dbF)$-scheme. So $\scrN=\scrN^\prime$ and $\scrN/H_0$ is a smooth, projective $O_{(v)}$-scheme.

We prove (a) and (b) in two Cases depending on the parity of $p$. We use the notations of the proof of 3.3; so $F_1$, $\tilde G_1$, $F_2$, $E_2$, $I_p$, $F_1\otimes_{\dbQ} \dbQ_p=\prod_{i\in I_p} F_{1i}$, $F_2\otimes_{\dbQ} \dbQ_p=\prod_{i\in I_p} F_{2i}$, $i_0\in I_p$, $F_3$, $T_{\dbZ_{(p)}}$, $T^\prime_{\dbQ}=Z^0(G)$, $T^c_{\dbZ_{(p)}}$ and $\scrB_{(p)}$ are as in the proof of 3.3. Warning: for both two Cases below we choose $F_2$ such that $[F_2:F_1]>>0$ and we also choose $T^\prime_{\dbQ}$ to be $T_{\dbQ}$. 

\medskip\noindent
{\bf 6.1.1. Case 1.} We assume that either $p\Ge 3$. So 6.1 (a) holds, cf. 1.5. We now show that for $[F_2:F_1]>>0$ there is no morphism 
$$q^{\stc}:\Spec(k[[x]])\to\scrC^{\stc}$$ 
whose generic fibre factors through $\scrN_{\dbF}/H_0$ and such that 
$$B:=q^{\stc*}(\scrS(N)\times_{\scrM/K(N)} \scrC^{\stc})$$ 
is a semiabelian scheme over $k[[x]]$ whose special fibre $B_k$ over $k$ is not an abelian variety. 

It suffices to show that the assumption that such a $q^{\stc}$ exists leads to a contradiction for $[F_2:F_1]>>0$. As in \S5 let $k_1:=\overline{k((x))}$ and let $z_1\in\scrN_{W(\dbF)}(W(k_1))$ be such that the special fibre of its composite with the quotient morphism $\scrN_{W(\dbF)}\to\scrN_{W(\dbF)}/H_0$ factors through the point $\Spec(k((x)))\to\scrN_{\dbF}/H_0$ defined naturally by $q^{\stc}$. Let $(M_1,F^1_1,\phi_1,(\tilde t_{1\alpha})_{\alpha\in\scrJ})$, $A_1$ and $\mu_1$ be obtained as in Case 1 of 5.2.       .

As $p\Ge 3$ from 4.2.1 and 4.2.2 1) we get the existence of an isomorphism 
$$\rho_1:(M_1,(\tilde t_{1\alpha})_{\alpha\in\scrJ})\arrowsim (L_{(p)}^*\otimes_{\dbZ_{(p)}} W(k_1),(v_{\alpha})_{\alpha\in\scrJ}).$$ 
We recall that $L_{(p)}$ and $L_{(p)}^*$ are naturally $E_{2(p)}$-modules and that $\scrB_{(p)}$ is an $E_{2(p)}$-subalgebra of $\End(L_{(p)})$ fixed by $G_{\dbZ_{(p)}}$, cf. 3.3.1. So as $E_{2(p)}$ is an $F_{2(p)}$-algebra we get that there is a unique direct sum decomposition
$$L_{(p)}^*\otimes_{\dbZ_{(p)}}\dbZ_p=\oplus_{i\in I_p} L^i$$ 
in $G_{\dbZ_p}$-modules such that each simple factor $F_{2i}$ of $F_2\otimes_{\dbQ} \dbQ_p$ acts trivially on $L^j[{1\over p}]$ if $j\neq i$ and there is no element of $L^i[{1\over p}]$ fixed by $F_{2i}$. This implies that each $L^i$ is self dual with respect to $\psi^*$. If $i\neq i_0$, the rank of $L^i$ is bounded in terms just of $(G_1,X_1)$ and so independently of $F_2$, cf. properties (vi) and (viii) of the proof of 3.3 and the fact that field $F_3$ depends only on $F_1$ and not on $F_2$. On the other hand, the rank of $L^{i_0}$ goes to infinity when $[F_2:F_1]$ goes to infinity (cf. (vi) of the proof of 3.3). 

We view $G_{W(k_1)}$ naturally as a closed subgroup of $GL(L_{(p)}^*\otimes_{\dbZ_{(p)}} W(k_1))$. We recall (see 4.3.1) that we denote by $\psi^*$ the perfect alternating form on $L_{(p)}^*$ and so also on $L_{(p)}^*\otimes_{\dbZ_{(p)}} W(k_1)$ defined by $\psi$; it is normalized by $G_{W(k_1)}$. Let $\mu^0:=\rho_1\mu_1\rho_1^{-1}:\dbG_{mW(k_1)}\to G_{W(k_1)}$. So $\rho_1\phi_1\rho_1^{-1}$ is a $\sigma_{k_1}$-linear endomorphism of $L_{(p)}^*\otimes_{\dbZ_{(p)}} W(k_1)$ of the form $g\sigma_{k_1}\mu^0({1\over p})$, where $g\in G_{\dbZ_{(p)}}(W(k_1))$ and $\sigma_{k_1}$ is identified here with the $\sigma_{k_1}$-linear automorphism of $L_{(p)}^*\otimes_{\dbZ_{(p)}} W(k_1)$ fixing $L_{(p)}^*\otimes_{\dbZ_{(p)}}\dbZ_p$. From the second paragraph of Case 1 of 5.2 we get that under any fixed  $O_{(v)}$-embedding $W(k_1)\hookrightarrow \dbC$, the cocharacters $\mu^0_{\dbC}$ and $\mu_h$'s of $G_{\dbC}$ are $G(\dbC)$-conjugate (here $h\in X$).

We recall from the proof the proof of 3.3 that the centralizer of $G_{\dbZ_{(p)}}$ in $GL(L_{(p)})$ is a reductive group scheme $C_{\dbZ_p}$. So the centralizer of $G_{\dbZ_p}$ in $GL(L_{(p)}^*\otimes_{\dbZ_{(p)}} \dbZ_p)$ is also a reductive group scheme and so has a maximal torus $T(C_{\dbZ_p})$ which is the group scheme of invertible elements of a commutative, semisimple $\dbZ_p$-algebra $F(C_{\dbZ_p})$. As $\scrB_{(p)}$ is an $E_{2(p)}$-algebra, we get that $F(C_{\dbZ_p})$ is in fact an $E_{2(p)}\otimes_{\dbZ_{(p)}} \dbZ_p$-semisimple algebra. Let $J_0$ be the set of distinct characters of the action of $T(C_{\dbZ_p})_{W(k_1)}$ on $L^{i_0}\otimes_{\dbZ_p} W(k_1)$. The Frobenius automorphism $\sigma_{k_1}$ acts naturally on $J_0$. Let $M^j$ be the direct summand of $L^{i_0}\otimes_{\dbZ_p} W(k_1)$ on which $T(C_{\dbZ_p})_{W(k_1)}$ acts via the character $j\in J_0$. We get a direct sum decomposition
$$L^{i_0}\otimes_{\dbZ_p} W(k_1)=\oplus_{j\in J_0} M^j$$
into irreducible $G_{W(k_1)}$-modules permuted by $\sigma_{k_1}$. 

We fix an orbit $o$ of the action of $\sigma_{k_1}$ on $J_0$. As $T_{\dbQ}=T^\prime_{\dbQ}$, the torus $T^c_{W(k_1)}$ is a subtorus of $G_{W(k_1)}$. The action of $T^c_{W(k)}$ on $L^{i_0}\otimes_{\dbZ_p} W(k_1)$ is via two non-trivial and distinct irreducible cocharacters. So there is a disjoint decomposition 
$$J_0=J_1\cup J_2$$ 
such that for $j_1$ and $j_2\in J_0$, the $T^c_{W(k)}$-modules $M^{j_1}$ and $M^{j_2}$ are isomorphic iff $(j_1,j_2)\in J_1^2\cup J_2^2$. As $T^c_{W(k)}$ fixes $\psi$, for $(j_1,j_2)\in J_1^2\cup J_2^2$ we have $\psi^*(M^{j_1},M^{j_2})=\{0\}$.

There is a simple factor of the adjoint group of the image of $G_{W(k_1)}$ in $GL(L^{i_0}\otimes_{\dbZ_p} W(k_1))$ in which $\mu^0$ has a trivial image. This is implied by the fact that $i_0\in I_p$ is a compact element in the sense of the proof of 3.3 and by the $G(\dbC)$-conjugacy property of the last sentence of the third paragraph above. So there is $j_0\in o$ such that the image of $\mu^0$ in $GL(M^{j_0})$ is contained in $Z(GL(M^{j_0}))$. As $E_2\otimes_{F_2} F_{2i_0}$ is a quadratic field extension of $F_{2i_0}$ and as $F(C_{\dbZ_p})[{1\over p}]$ is an $E_2\otimes_{\dbQ} \dbQ_p$-algebra, the orbit $o$ contains at least two characters $j_1\in J_1$ and $j_2\in J_2$ such that $M^{j_1}$ and $M_{j_2}$ are not perpendicular with respect to $\psi^*$. So there is $i\in\{1,2\}$ and $j_i\in o\cap J_i$, such that the image of $\mu^0$ in $GL(M^{j_i})$ is trivial. 
So the multiplicity of the slope $1$ for $(L^{i_0}\otimes_{\dbZ_p} W(k_1),\tilde gg\sigma_{k_1}\mu^0({1\over p}))$ is $0$ for any $\tilde g\in G_{W(k_1)}(W(k_1))$. So as $L^{i_0}$ is self dual with respect to $\psi^*$, the multiplicity of the slope $0$ for $(L^{i_0}\otimes_{\dbZ_p} W(k_1),\tilde gg\sigma_{k_1}\mu^0({1\over p}))$ is also $0$ for any $\tilde g\in G_{W(k_1)}(W(k_1))$. 

So the $p$-rank of $A_{1k_1}$ (and so also of the generic fibre $B_{k((x))}$ of $B$) is at most
$$n_0:={1\over 2}\sum_{i\in I_p\setminus\{i_0\}} \dim_{\dbZ_p}(L^i)$$ 
and so bounded above just in terms of $(G_1,X_1)$ and so independently of $[F_2:F_1]$. 

Let $T(B_k)$ be the toric part (i.e. the maximal torus) of the identity component of $B_k$. As $B_k$ is not an abelian scheme, we have $1\Le\dim_k(T(B_k))$. From [FC, Cor. 5.11 of p. 70] we get that $\forall m\in\dbN$ the finite group scheme $T(B_k)[p^m]_{k((x))}$ is naturally identified with a closed subgroup of $B_{k((x))}[p^n]$ and that this identification is compatible with the $F_{2(p)}$-actions. So the $p$-rank of $A_{1k_1}$ (and so also of $B_{k((x))}$) is at lest $\dim_k(T(B_k))$. So we have $\dim_k(T(B_k))\Le n_0$. But $F_{2(p)}$ acts on $B$ and so also on the $\dbZ$-module $X_*(T(B_k))$ of cocharacters of $T(B_k)$. So as $1\Le\dim_{\dbZ}(X_*(T(B_k)))\Le n_0$, by taking $[F_2:F_1]>n_0$ we get that this action is trivial. So $F_{2(p)}$ acts trivially on $T(B_k)$. As $T_0[p^m]_{k((x))}$ is naturally identified with a closed subgroup of $B_{k((x))}[p^m]$ and as this identification is compatible with the $F_{2(p)}$-actions, we get that $F_{2(p)}$ acts trivially on a direct summand of $M_1$ of positive rank. So via $\rho_1$, $F_{2(p)}$ acts trivially on a direct summand of $L_{(p)}^*\otimes_{\dbZ_{(p)}} W(k_1)$ of positive rank. But this contradicts 3.3.1. So a morphism $q^{\stc}$  as in the beginning of 6.1.1 dos not exist. This implies that 6.1 (b) holds for $p\Ge 3$. 

\medskip\noindent
{\bf 6.1.2. Case 2.} Let $p=2$. Let $(T^0_{\dbQ},\{h\})$, $H^0$ and $\scrT^0$ be as in the proof of 4.4 (i). We fix a point $z\in \im(\scrT^0(W(\dbF))\to\scrN(W(\dbF)))$. The statements of 4.2.2 1) and 2) hold for $k=\dbF$ and for $z\in\scrT^0(W(\dbF))$, cf. the toric part of the hypotheses of 4.2.2 applied in the context of the morphism $\scrT^0\to\scrM$ and so of the injective map $(T^0_{\dbQ},\{h\})\hookrightarrow (GSp(W,\psi),S)$ of Shimura pairs. Let $\scrN^0_{k(v)}$ be the connected component of $\scrN_{k(v)}$ to which $z$ specializes, cf. 1.4 1). Let $H_{00}$ be the maximal closed subgroup of $H_0$ normalizing $\scrN^0_{k(v)}$. 

If $q_{R}:\Spec(R)\to\scrN$ is as in 4.3.3 and if $\scrC_{\text{univ}}$ is as in 4.6, then the pull back of $(M_{R},(\tilde t_{\alpha})_{\alpha\in\scrJ})$ via any $W(k_1)$-valued point of $\Spec(R)$, is isomorphic to $(M\otimes_{W(k)} W(k_1),(\tilde t_{\alpha})_{\alpha\in\scrJ})$. Moreover the pull back of $w_{\alpha}^{\scrA^\prime}$ via $q_{R[{1\over p}]}$ is a Hodge cycle on $A_{R[{1\over p}]}$ having $\tilde t_{\alpha}\in\scrT(M)\otimes_{W(k)} R[{1\over p}]$ as its de Rham component, cf. 4.3.3. So as 4.2.2 1) and 2) hold for $z$ and due to 4.2.1, we get that the logical analogues of 4.2.2 1) and 2) also hold for any $W(k_1)$-valued point of $\scrN$ factoring through the formally smooth morphism $q_{R}:\Spec(R)\to\scrN$. 

We take $\scrC$ to be the connected component of $\scrN^\prime_{W(\dbF)}/H_0$ whose special fibre $\scrC_{\dbF}$ contains $\scrN^0_{\dbF}/H_{00}$. Let $q^{\stc}:\Spec(k[[x]])\to\scrC^{\stc}$ be a dominant morphism; so its generic fibre factors through a generic point of $\scrN^0_{\dbF}/H_{00}$. We consider two Subcases and we use the notations of Case 1.

\medskip\noindent
{\bf 6.1.2.1. Subcase 1.} We consider the Subcase when $v_1$ is $p$-compact for $(G_1,X_1)$. The prime $v_1$ is $p$-compact iff each $i\in I_p$ is a compact element. So as in Case 1 we argue that $(L^i,g\sigma_{k_1}\mu^0({1\over p}))$ has no slope 1, $\forall i\in I_p$. So also $(M_1,\phi_1)$ has no slope 1. So the $p$-rank of $A_{1k_1}$ is $0$. So as $q^{\stc}$ is dominant, the pull backs of $\scrA^\prime$ via geometric points of $\scrN^0_{k(v)}$ have $p$-ranks equal to $0$. So from the proof of 1.5 applied just to the connected component $\scrN^0_{k(v)}$ of $\scrN_{k(v)}$, we get that $\scrN^0_{k(v)}$ is an open closed subscheme of $\scrN^{\prime}_{k(v)}$. Moreover, either from loc. cit. or from [Oo, 1.3] we get that $A_{1k_1}$ can no specialize to a semiabelian variety which is not an abelian variety. So $\scrN^0_{\dbF}/H_{00}$ is a smooth, open closed subscheme of $\scrC^{\stc}_{\dbF}$. So as in Case 1 we argue that $\scrC=\scrC^{\stc}$ is smooth and projective. So as the connected components of $\scrN$ or $\scrN^\prime$ are permuted transitively by $G(\dbA_f^{(p)})$ (cf. [Va1, 3.2.2]), we get that $\scrN=\scrN^\prime$ and that $\scrN_{k(v)}/H_0$ is a proper scheme. So also 6.1 (a) and (b) hold in this case. 

\medskip\noindent
{\bf 6.1.2.2. Subcase 2.} We consider the Subcase when $(G_1,X_1)$ is of $C_n$ type and has compact factors. We will first use 5.6.1 to show that $\scrN^0_{k(v)}$ is an open closed subscheme of $\scrN^{\prime}_{k(v)}$. 

Let $q:\Spec(k[[x]])\to\scrN_{k(v)}^\prime$ be a morphism whose generic fibre factors through the generic point of $\scrN^0_{k(v)}$. Let $R_0=W(k)[[x]]$, $\Phi_{R_0}$, $(M_0,\Phi_0,\nabla_0,p_{M_0})$ and $(\tilde t_{1\alpha})_{\alpha\in\scrJ}$ be as in 5.1. Let $G^1_{\dbZ_p}:=G_{\dbZ_p}$. Let $\tilde G_0=\tilde G_0^1$ be as in 5.1.2; it is a reductive group scheme over $R_0$ (cf. 5.2). Let $(M_1,F^1_1,\phi_1,(\tilde t_{1\alpha})_{\alpha\in\scrJ})$, $A_1$, $\mu_1$, $\mu_0:\dbG_{mR_0}\to\tilde G_0$ and $(M,F^1,\phi,\tilde G,(\tilde t_{\alpha})_{\alpha\in\scrJ},p_M)$ be as in Case 1 of 5.2 and 5.3. Let $y$ be as in 5.1 and let $(A_k,p_{A_k})$ be as in 5.4. 

In order to apply 5.6.1, we will check that conditions 5.6 (i) to (iii) hold. The group scheme $G_{\dbZ_{(p)}}$ is reductive (cf. 3.3 (ii)) and we have $G_{\dbZ_p}^1=G_{\dbZ_p}$; so 5.6 (i) holds. The point $z_1:\Spec(W(k_1))\to\scrN$ is such that its generic fibre factors through the generic point of $\scrN_{k(v)}^0$ and so implicitly through the special fibre $q_{R/pR}$ of $q_R$. So 4.2.2 2) applies in the context of $z_1$ (see beginning of 6.1.2) and so 5.6 (ii) holds. 

We now check that 5.6 (iii) holds. Let $\rho_1$ be as in Case 1. As 5.6 (i) and (ii) hold, from the first paragraph of the proof of 5.6.1 we get that $\tilde G_0$ is a reductive subgroup of $GL(M_0)$. So let $\rho_{1z}$ be as in 5.5. So $\rho_1\rho_{1z}^{-1}$ is an isomorphism $(M\otimes_{W(k)} W(k_1),(\tilde t_{\alpha})_{\alpha\in\scrJ})\arrowsim (L_{(p)}^*\otimes_{\dbZ_{(p)}} W(k_1),(v_{\alpha})_{\alpha\in\scrJ})$. As $k=\overline{k}$, we get the existence of an isomorphism $\rho:(M,(\tilde t_{\alpha})_{\alpha\in\scrJ})\arrowsim (L_{(p)}^*\otimes_{\dbZ_{(p)}} W(k),(v_{\alpha})_{\alpha\in\scrJ})$. Let $M_i:=\rho^{-1}(L^i\otimes_{\dbZ_p} W(k))$. Let $\tilde G$ be the reductive subgroup of $GL(M)$ which is the closed subscheme of $\tilde G_0$ defined by the zero locus $x=0$. We get a direct sum decomposition of $\tilde G$-modules
$$M=\oplus_{i\in I_p} M_i$$ 
\indent
We choose the principally polarized abelian scheme $(A,p_A)$ over $W(k)$ lifting $(A_k,p_{A_k})$ such that the principally quasi-polarized $p$-divisible group $(D,p_D)$ of $(A,p_A)$ is a direct sum $\oplus_{i\in I_p} (D_i,p_{D_i})$, where the filtered $F$-crystal of $D_i$ is $(M_i,F^1\cap M_i,\phi)$. Implicitly $H^1_{\acute et}(A_{B(k)},\dbZ_p)$ is a direct sum $\oplus_{i\in I_p} H^{1i}$, where $H^{1i}$ is the dual of the Tate-module of $D_i$. As we have $G^{\ad}_{\dbQ_p}=G_{1\dbQ_p}=\prod_{i\in I_p} \Res_{F_{1i}/\dbQ}\tilde G_{1F_{1i}}$, we can define the closed, flat subgroup $G^{i\der}_{\dbZ_p}$ of $GL(H^1_{\acute et}(A_{B(k)},\dbZ_p))$ whose generic fibre corresponds to the semisimple subgroup of $G^{\der}_{\dbQ_p}$ having $\Res_{F_{1i}/\dbQ}\tilde G_{1F_{1i}}$ as its adjoint, via the Fontaine comparison theory of $A_{B(k)}$. Also, we can define the closed, flat subgroup $Z^0_{\dbZ_p}$ of $GL(H^1_{\acute et}(A_{B(k)},\dbZ_p))$ whose generic fibre corresponds to $Z^0(G_{\dbQ_p})$ via the same theory. 

We view naturally $\scrB_{(p)}\otimes_{\dbZ_{(p)}} \dbZ_p$ as a $\dbZ_p$-algebra of $\dbZ_p$-endomorphisms of $\scrA^\prime$. So $\scrB_{(p)}^{\text{opp}}\otimes_{\dbZ_{(p)}} \dbZ_p$ acts naturally on $M$. We consider the product decomposition $\scrB_{(p)}^{\text{opp}}\otimes_{\dbZ_{(p)}} \dbZ_p=\prod_{i\in I_p}\scrB_i$ such that $\scrB_i[{1\over p}]$ is an $F_i$-algebra. So $\scrB_i$ acts faithfully on $M_i$ and fixes $M/M_i$. We show that $G^{i\der}_{\dbZ_p}$ is a semisimple group scheme over $\dbZ_p$, $\forall i\in I_p$.

If $i\in I_p$ is a compact element, then as in Case 1 we argue that $(M_i,\tilde g\phi)$ has no slope 1, $\forall\tilde g\in\tilde G(W(k))$. So $(D_i,p_{D_i})$ is uniquely determined, cf. [Va6, 2.2.3]. Let $(t_{\alpha})_{\alpha\in\scrJ_i}$ be a family of tensors of $\scrT(M_i)$ fixed by $\phi$ and such that the image of $\tilde G$ in $GL(M_i)$ is the Zariski closure of the subgroup of $GL(M_i[{1\over p}])$ fixing $t_{\alpha}$, $\forall\alpha\in\scrJ_i$. Its existence is implied by [Va6, 2.5.3]. From [Va6, 1.3 a)] applied to $(D_i,(t_{\alpha})_{\alpha\in\scrJ_i})$ we get that $G^{i\der}_{W(k)}$ is isomorphic to the image of $\tilde G^{\der}$ in $GL(M_i)$ and so $G^{i\der}_{\dbZ_p}$ is a semisimple group scheme. 

We now consider an $i\in I_p$ which is not a compact element. Let $\tilde G^i_0$ be the subgroup of $GL(M_i)$ generated by the image of $\tilde G^{\der}$ in $GL(M_i)$ and by $Z(GL(M_i))$. It is a reductive group scheme, cf. 2.2.5 2). 
Let $p_{M_i}$ be the perfect alternating form on $M_i$ defined by $p_M$. Let $\mu$ be the cocharacter of $\tilde G$ which is the reduction modulo $x$ of the cocharacter $\mu_0:\dbG_{mR_0}\to\tilde G_0$ of 5.3. Let $\mu_i:\dbG_{mW(k)}\to GL(M_i)$ be the cocharacter defined by $\mu$. As $i\in I_p$ is not a compact element, $\mu_i$ factors through $\tilde G^i_0$. The cocharacter $\mu_0:\dbG_{mR_0}\to\tilde G_0$ is uniquely determined up to inner conjugation under an element $g_0\in\tilde G_0(R_0)$ which normalizes the kernel $F^1_0/pF^1_0$ of $\Phi_0$ mod $p$. Let $S^i$ be the semisimple $\dbZ_p$-algebra of elements of $\End(M_i)$ fixed by $\tilde G^i_0$ and by $\phi$. We identify its opposite $S^{i\text{opp}}$ with a $\dbZ_p$-algebra of endomorphisms of the special fibre $D_{ik}$ of $D_i$. The subgroup of $GSp(M_i,p_{M_i})$ fixing $S^i\otimes_{\dbZ_p} W(k)$ is $\tilde G^i_0$. As all simple factors of the adjoint group of $\tilde G^{i\der}_0$ are of $C_n$ Lie type, it is well known that the formal moduli space of the pair $((D_i,p_{D_i})_k,S^{i\text{opp}})$ is formally smooth and representable by a scheme of formal power series over $W(k)$. 

So the pair $((D_i,p_{D_i})_k,S^{i\text{opp}})$ has a lift $((D_i^\prime,p_{D_i^\prime}),S^{i\text{opp}})$ to $\Spec(W(k))$. Let $F^{1\prime}$ be the direct summand of $M_i$ which is the Hodge filtration of $D_i^\prime$. Let $\mu_i^\prime:\dbG_{mW(k)}\to GL(M_i)$ be the inverse of the canonical split cocharacter of $(M_i,F^{1\prime},\phi)$ defined in [Wi, p. 512]. As the group $\tilde G^i_0$ is the subgroup of $GL(M_i,p_{M_i})$ centralizing $S^i\otimes_{\dbZ_p} W(k)$ and as $\mu_i^\prime$ fixes each element of $S^i$ and normalizes $p_{M_i}$ (cf. the functorial aspects of [Wi, p. 513]), $\mu_i^\prime$ factors through $\tilde G^i_0$. So $\mu_i^\prime$ is conjugate with $\mu_i$ under an element $g_0^i\in\tilde G^i_0(W(k))$ normalizing $F^1/pF^1\cap M_i/pM_i$, cf. [Va5, 3.1.2]. It is easy to see that we can choose the $g_0^i$'s for $i\in I_p$ not a compact element, such that there is $g_0\in\tilde G_0(R_0)$ normalizing $F^1_0/pF^1_0$ and with the property that its image in $\tilde G_i^0(W(k))$ is $g_0^i$, for any $i\in I_p$ which is not a compact element. Not to introduce extra notations we will assume that $\mu_i=\mu_i^\prime$, for any $i\in I_p$ which is not a compact element.

So we can also assume that $(D_i,p_{D_i})=(D_i^\prime,p_{D_i^\prime})$. We denote also by $S^i$ the $\dbZ_p$-algebra of endomorphisms of $H^{1i}$ corresponding to the endomorphisms $S^{i\text{opp}}$ of $D_i=D_i^\prime$. Let $p_{H^{1i}}$ be the perfect alternating form on $H^{1i}$ which is the \'etale realization of $p_{D_i}=p_{D_i^\prime}$ (equivalently of $p_{M_i}$). As $(G_1,X_1)$ is of $C_n$ type, it is well known that the subgroup of $GL(H^{1i})$ fixing $S^i$ and normalizing $p_{H^{1i}}$ is reductive. Its derived group is $G^{i\der}_{\dbZ_p}$. 

We conclude that $G^{i\der}_{\dbZ_p}$ is a semisimple group scheme over $\dbZ_p$, $\forall i\in I_p$. The group scheme $Z^0_{\dbZ_p}$ is the Zariski closure of a subtorus of the generic fibre of the torus of $GL(H^1_{\acute et}(A_{B(k)},\dbZ_p))$ which is the center of the centralizer of $\prod_{i\in I_p} \scrB_i$ in $GL(H^1_{\acute et}(A_{B(k)},\dbZ_p))$. So $Z^0_{\dbZ_p}$ is a torus. From 2.2.5 2) we get that  the Zariski closure $G_{\dbZ_p}^{\acute et}$ in $GL(H^1_{\acute et}(A_{B(k)},\dbZ_p))$ of the subgroup of $GL(H^1_{\acute et}(A_{B(k)}/\dbQ_p))$ fixing $\tilde v_{\alpha}$, $\forall\alpha\in\scrJ$, is reductive.

The center of the centralizer of $\prod_{i\in I_p} \scrB^i$ in $GSp(H^1_{\acute et}(A_{B(k)},\dbZ_p),\oplus_{i\in I_p} p_{H^{1i}})$ is $Z^0_{\dbZ_p}$, cf. 3.3.2. So as in [Va5, 7.8.3.1] (working under the assumptions [Va5, 7.8.3 1) and 2a)]) we get that there is an isomorphism $(M,(\tilde t_{\alpha})_{\alpha\in\scrJ})\arrowsim (H^1_{\acute et}(A_{B(k)},\dbZ_p)\otimes_{\dbZ_p} W(k),(\tilde v_{\alpha})_{\alpha\in\scrJ})$. So 5.6 (iii) holds.

So all hypotheses of 5.6.1 hold and so from 5.6.1 we get that the morphism $q:\Spec(k[[x]])\to\scrN^\prime_{k(v)}$ factors through $\scrN_{k(v)}$. This implies that $\scrN^0_{k(v)}$ is an open closed subscheme of $\scrN^{\prime}_{k(v)}$. As in Case 1 we check that for $[F_2:F_1]>>0$, there is no morphism $q^{\stc}:\Spec(k[[x]])\to\scrC^{\stc}$ whose generic fibre factors through $\scrN^0_{\dbF}/H_0$ and such that the natural pull back of $\scrS(N)$ to $\Spec(k[[x]])$ is a semiabelian scheme which is not an abelian scheme. The rest of the arguments needed to prove 6.1 (a) and (b) are as in Subcase 1. So 6.1 (a) and (b) hold for $p=2$ as well.

\medskip\noindent
{\bf 6.1.3. More on the $A_n$ type.} We assume $(G_1,X_1)$ is a simple, adjoint Shimura pair of $A_n$ type and with compact factors. Let $F_2$ and $(G_2,X_2)$ be defined as in the proof of 3.3. We apply 3.3 to $(G_2,X_2)$ instead of $(G_1,X_1)$ and we denote also by $f:(G,X)\hookrightarrow (GSp(W,\psi),S)$ the resulting embedding. So to end this Step 1, we point out that entirely as above (for the $B_n$, $C_n$, $D_n^{\dbH}$ and $D_n^{\dbR}$ types) we argue that for $[F_2:F_1]>>0$ the $O_{(v)}$-scheme $\scrN/H_0$ is smooth and projective.

\bigskip\noindent
{\bf 6.2. Step 2.} In this Step 2 we prove 1.6 2) in general. So we work with a general Shimura quadruple $(G_1,X_1,H_1,v_1)$ satisfying the hypotheses of 1.6 2). If $G_1$ is a torus, then the fact that $\Sh_{H_1}(G_1,X_1)$ has an integral canonical model over $O_{(v_1)}$ is well known (for instance, see [Va1, 3.2.8]). So for the rest of the proof (i.e. until 6.3) we will assume that the group $G_1^{\ad}$ is non-trivial. Let $G_{1\dbZ_{(p)}}$ be the reductive group scheme over $\dbZ_{(p)}$ extending $G_1$ and such that $H_1$ is its group scheme of $\dbZ_p$-valued points, cf. beginning of \S3. 

Let $H_1^{\ad}:=G^{\ad}_{1\dbZ_{(p)}}(\dbZ_p)$. Let  $v_1^{\ad}$ be the prime of $E(G_1^{\ad},X_1^{\ad})$ divided by $v_1$. We refer to the Shimura quadruple $(G_1^{\ad},X_1^{\ad},H_1^{\ad},v_1^{\ad})$ as the adjoint Shimura quadruple of $(G_1,X_1,H_1,v_1)$. We consider its product decomposition
$$(G_1^{\ad},X_1^{\ad},H_1^{\ad},v_1^{\ad})=\prod_{i\in I_0} (G_1^{i},X_1^{i},H_1^i,v_1^i)$$ 
into simple factors. So each $G_1^{i}$ is a simple, adjoint $\dbQ$--group extending to an adjoint group scheme $G^{i}_{1\dbZ_{(p)}}$ over $\dbZ_{(p)}$ such that we have $G_{1\dbZ_{(p)}}^{\ad}=\prod_{i\in I_0} G^i_{\dbZ_{(p)}}$.

For each $i\in I_0$ we apply 3.3 and 6.1. So we consider an injective map  
$f^i:(G^i,X^i)\hookrightarrow (GSp(W^i,\psi^i),S^i)$ such that there is a $\dbZ$-lattice $L^i$ of $W^i$ with the property that we have a perfect alternating form $\psi^i:L^i\otimes_{\dbZ} L^i\to\dbZ$, the analogues of 3.3 (i) to (iv) hold for $f^i$  but with $(G_1,X_1)$ being replaced by $(G_1^i,X_1^i)$ and moreover (cf. 6.1) the integral canonical model $\scrN^i$ of $\Sh_{H^i}(G^i,X^i)$ over $O(G_1^i,X_1^i,p)$ exists and is obtained in the same $\scrN^\prime$ was obtained in 1.3. Here $H^i$ is the hyperspecial subgroup of $G^i_{\dbQ_p}$ which is the group of $\dbZ_p$-valued points of the Zariski closure of $G^i$ in $GL(L^i\otimes_{\dbZ} \dbZ_{(p)})$.
 
Let $(G_3,X_3)\to (G_1,X_1)$ be a map of Shimura pairs such that the following five properties hold (cf. [MS, proof of 3.4] and [Va1, 3.2.7 10)]):

\medskip
{\bf (i)} {\it it is a cover in the sense of [Va1, 2.4], i.e. $G_3$ surjects onto $G_1$ and $\Ker(G_3\to G_1)$ is a torus such that $H^1(K,\Ker(G_3\to G_1))=\{0\}$, for any field $K$ of characteristic $0$;}

\smallskip
{\bf (ii)} {\it we have $E(G_3,X_3)=E(G_1,X_1)$;}

\smallskip
{\bf (iii)} {\it the group $G_3^{\der}$ is simply connected if $p>2$ and is $\prod_{i\in I_0} G^{i\der}$ if $p=2$;}

\smallskip
{\bf (iv)} {\it the group $G_{3\dbQ_p}$ is unramified;}

\smallskip
{\bf (v)} {\it we have a product decomposition $(G_3,X_3)=\prod_{i\in I_0} (G_3^i,X_3^i)$ such that the adjoint of $(G_3^i,X_3^i)$ is $(G_1^i,X_1^i)$.}

\medskip
Let $G_{3\dbZ_{(p)}}$ be the reductive group scheme over $\dbZ_{(p)}$ extending $G_3$ and having $G^{\ad}_{1\dbZ_{(p)}}$ as its adjoint. Let $H_3:=G_{3\dbZ_{(p)}}(\dbZ_p)$.

Let $i\in I_0$. If $p=2$ or if $(G^i_1,X^i_1)$ is not of $D_n^{\dbH}$ type, then let $(G^{i0},X^{i0}):=(G^i,X^i)$ and $\scrN^{i0}:=\scrN^i$. If $p>2$ and $(G^i_1,X^i_1)$ is of $D_n^{\dbH}$ type, then we perform the following construction. As above we consider a cover $(G^{i0},X^{i0})\to (G^i,X^i)$ with $G^{i0}_{\dbQ_p}$ unramified, with $G^{i0\der}$ simply connected and with $E(G^{i0},X^{i0})=E(G^i,X^i)$. Let $H^{i0}$ be the hyperspecial subgroup of $G^{i0}_{\dbQ_p}(\dbQ_p)$ which is the analogue of $H_3$ but obtained in the context of the cover $(G^{i0},X^{i0})\to (G^i,X^i)$. As $p>2$, the logical analogue of 4.2.2 2) always holds for any $z\in\scrN^i(W(k))$. So the triple $(f^i,L^i,v^i)$ is a standard Hodge situation in the sense of [Va5, 1.4], for any prime $v^i$ of $E(G^i,X^i)$ dividing $p$. So the integral canonical model $\scrN^{i0}$ of $\Sh_{H^{i0}}(G^{i0},X^{i0})$ over $O(G^{i0},X^{i0},p)=O(G^i,X^i,p)$ exists and is a pro-\'etale cover of $\scrN^i$, cf. [Va5, 7.9.3].

So regardless of who $p$ and $(G^i_1,X^i_1)$'s are we get:

\medskip
{\bf (vi)} {\it we have $O(G^{i0},X^{i0},p)=O(G^i,X^i,p)$ and $\scrN^{i0}$ is a pro-\'etale cover of a smooth, (quasi-) projective  $O(G^{i0},X^{i0},p)$-scheme iff $\scrN^i$ is a pro-\'etale cover of a smooth, (quasi-) projective $O(G^{i},X^{i},p)$-scheme.}

\medskip
Let $(G_4,X_4,H_4):=\prod_{i\in I_0} (G^{i0},X^{i0},H^{i0})$. The product of the $O(G_4,X_4,p)$-schemes $\scrN^{i0}_{O(G_4,X_4,p)}$'s is the integral canonical model $\scrN_4$ of $\Sh_{H_4}(G_4,X_4)$ over $O(G_4,X_4,p)$. We have $G_4^{\der}=G_3^{\der}$. So as in [Va1, 6.2.3] we argue that:

\medskip
{\bf (vii)}  {\it $\Sh_{H_3}(G_3,X_3)$ has an integral canonical model $\scrN_3$ over $O(G_3,X_3,p)$ and that for any pair of primes $v_4$ and $v_3$ of $O(G_4,X_4,p)$ and respectively of $O(G_3,X_3,p)$ dividing the same prime of $O(G_1^{\ad},X_1^{\ad},p)=O(G_4^{\ad},X_4^{\ad},p)=O(G_3^{\ad},X_3^{\ad},p)$, the connected components of $\scrN_{4O_{(v_4)}^{\sh}}$ and $\scrN_{3O_{(v_3)}^{\sh}}$ are isomorphic.}

\medskip
From (vii) we get  

\medskip
{\bf (viii)} {\it $\scrN_3$ is a pro-\'etale cover of a smooth, (quasi-) projective $O(G_3,X_3,p)$-scheme iff $\scrN_4$ is a pro-\'etale cover of a smooth, (quasi-) projective $O(G_4,X_4,p)$-scheme.}

\medskip
From 3.1.1 4) and 3.1.5 we also we get that:

\medskip
{\bf (ix)} {\it $\scrN^i$ is a pro-\'etale cover of a smooth, quasi-projective $O(G^i,X^i,p)$-scheme.}

\medskip
We now check the following two things:

\medskip
{\bf (a)} {\it the integral canonical model $\scrN_1^{\ad}$ of $\Sh_{H_1^{\ad}}(G_1^{\ad},X_1^{\ad})$ over $O(G_1^{\ad},X_1^{\ad},p)$ exists and is a pro-\'etale cover of a smooth, quasi-projective $O(G_1^{\ad},X_1^{\ad},p)$-scheme;}

{\bf (b)} {\it if $(G^i_1,X^i_1)$ has compact factors for all $i\in I_0$, then $\scrN_1^{\ad}$ is a pro-\'etale cover of a smooth, projective $O(G_1^{\ad},X_1^{\ad},p)$-scheme.}

\medskip
The same argument showing that $\scrN_4$ exists, shows that it suffices to prove (a) and (b) under the extra assumption that $I_0$ has only one element $i$. So $(G_1,X_1)$ is a simple, adjoint Shimura pair and so we can appeal to the notations of 6.1. However, for the sake of uniformity of notations we will continue to use the notation $\scrN_1^{\ad}$ instead of $\scrN_1$. As $I_0$ has only one element $i$, we have $\scrN_4=\scrN^{i0}$. To check (a) and (b) we consider three Cases.

\medskip\noindent
{\bf 6.2.1. Case 1: $(G_1,X_1)$ is of $A_n$ type.} As 3.3 (iv) holds for $F^i$ and $L^i$, the triple $(f^i,L^i,v^i)$ is a standard PEL situation for any prime $v^i$ of $E(G^i,X^i)$ dividing $p$. So from [Va4, 5.3.3] we get that $\scrN_1^{\ad}$ exists and that $\scrN^i$ is a pro-\'etale cover of a non-empty, open closed subscheme of $\scrN^{\ad}_1$. From this and 6.2 (ix) we get that a non-empty, open closed subscheme of $\scrN_1^{\ad}$ is a pro-\'etale cover of a smooth, quasi-projective $O(G_1,X_1,p)$-scheme. So as the connected components of $\scrN_1^{\ad}$ are permuted transitively by $G_1(\dbA_f^{(p)})$ (cf. [Va1, 3.3.2]), we get that $\scrN_1^{\ad}$ is a pro-\'etale cover of a smooth, quasi-projective $O(G_1,X_1,p)$-scheme. So 6.2 (a) holds. If $(G_1,X_1)=(G^i_1,X^i_1)$ has compact factors, then we can choose $f^i$ such that $\scrN^i$ is a pro-\'etale cover of a projective $O(G^i,X^i,p)$-scheme (cf. 6.1.3). So similarly we argue that $\scrN_1^{\ad}$ is a pro-\'etale cover of a smooth, projective $O(G_1,X_1,p)$-scheme; so 6.2 (b) also holds. 

\medskip\noindent
{\bf 6.2.2. Case 2: $p\Ge 3$ and $(G_1,X_1)$ is of $B_n$, $C_n$, $D_n^{\dbH}$ or $D_n^{\dbR}$ type.} In this Case the order of the center of $G^{\sc}=G_1^{\sc}$ is a power of $2$ and so prime to $p$. So from [Va1, 6.2.2 a)] applied in the context of the cover $(G_3,X_3,H_3)\to (G_1,X_1,H_1)$ we get that $\scrN_1^{\ad}$ exists and that $\scrN_3$ is a pro-\'etale cover of $\scrN_1^{\ad}$. From 6.2 (vii), (viii) and (ix) we get that $\scrN_3$ is a pro-\'etale cover of a smooth, quasi-projective $O(G_3,X_3,p)$-scheme. From the last two sentences we get that $\scrN_1^{\ad}$ is a pro-\'etale cover of a smooth, quasi-projective $O(G_1,X_1,p)$-scheme and so that 6.2 (a) holds. If $(G_1,X_1)=(G^i_1,X^i_1)$ has compact factors, then from 6.1 we get that we can choose $f^i$ such that $\scrN^i$ is a pro-\'etale cover of a projective $O(G^i,X^i,p)$-scheme. So based on 6.2 (vi) and (viii) and the equality $\scrN_4=\scrN^{i0}$ we get that 6.2 (b) holds.

\medskip\noindent
{\bf 6.2.3. Case 3: $p=2$ and $(G_1,X_1)$ is of $B_n$, $C_n$, $D_n^{\dbH}$ or $D_n^{\dbR}$ type.} Let $\tilde G^\prime_{\dbZ_{(p)}}$ be as in 3.3.2 and 3.3.3. Let $\tilde G^\prime:=\tilde G^\prime_{\dbQ}$. Let $\tilde X^\prime$ be such that the monomorphism $G\hookrightarrow\tilde G^\prime$ extends to an injective map $(G,X)\hookrightarrow (\tilde G^\prime,\tilde X^\prime)$ of Shimura pairs. Let $(G_2^\prime,X_2^\prime):=(\tilde G^{\ad},\tilde X^{\ad})$. Let $\tilde H^\prime:=\tilde G^\prime_{\dbZ_{(p)}}(\dbZ_p)$ and $H_2^\prime:=\tilde H^{\prime\ad}=G^\prime_{2\dbZ_{(p)}}(\dbZ_p)$. Each simple factor of $G_{2\dbC}^{\prime\ad}$ is of $A_s$ Lie type, where $s\in\{2n,2^n,2^{n-1}\}$ (cf. 3.3.2 and 3.3.3). So based on 6.2.1 we get that the integral canonical model $\tilde\scrN^\prime$ (resp. $\scrN_2^\prime$) of $\Sh_{\tilde H^\prime}(\tilde G^\prime,\tilde X^\prime)$ (resp. of $\Sh_{H_2^\prime}(G_2^\prime,X_2^\prime)$) over the localization of $E(\tilde G^\prime,\tilde X^\prime)_{(p)}$ (resp. of $E(G_2^\prime,X_2^\prime)_{(p)}$) with respect to the finite primes divided by the finite primes of $O(G_1,X_1,p)$, exists.

We denote by $\tilde f^\prime$ and $f_2^\prime$ the injective maps $\tilde f:(G,X)\hookrightarrow (\tilde G^\prime,\tilde X^\prime)$ and $f_2:(G^{\ad},X^{\ad})=(G_1,X_1)\hookrightarrow (G_2^\prime,X_2^\prime)$ defined naturally by $f$. Let $\scrN^{\ad}_1$ be the normalization of $\scrN^\prime_{2O(G_1,X_1,p)}$ in $\Sh_{H^{\ad}}(G^{\ad},X^{\ad})$. We know that $\scrN^{\ad}_1$ is a normal integral model of $\Sh_{H^{\ad}}(G^{\ad},X^{\ad})=\Sh_{H_1}(G_1,X_1)$ over $O(G_1,X_1,p)$ having the extension property, cf. 3.1.1 1). Moreover we have a commutative diagram  
$$
\spreadmatrixlines{1\jot}
\CD
\scrN @>{\tilde s^\prime}>> \tilde\scrN^\prime_{O(G_1,X_1,p)}\\
@V{q}VV @VV{q_2}V\\
\scrN^{\ad}_1@>{s_2^\prime}>>\scrN^\prime_{2O(G_1,X_1,p)},
\endCD
$$
where $\tilde s^\prime$ and $s_2^\prime$ are naturally defined by $\tilde f^\prime$ and respectively by $f_2^\prime$ and where $q$ and $q_2$ are the natural morphisms whose existence is implied by the fact that $\scrN_2^\prime$ has the extension property. We know that $q_2$ is a pro-\'etale cover of its image (cf. 6.2.1 and [Va4, 5.3.1]), that $\tilde s^\prime$ and $s_2^\prime$ are finite morphisms (cf. 3.1.1 3)), that the generic fibres of $\tilde s^\prime$ and $s_2^\prime$ are closed embeddings (cf. [Va1, 3.2.14]) and that the generic fibre of $q$ is a pro-\'etale morphism (this last thing is well known). We easily get that $q$ is a pro-\'etale cover of its image $IM$ and that $IM$ is an open closed subscheme of $\scrN^{\ad}_1$. So $IM$ is also a regular, formally smooth $O(G_1,X_1,p)$-scheme. The connected components of $\scrN_1^{\ad}$ are permuted transitively by $G_1(\dbA_f^{(p)})$ (cf. [Va1, 3.3.2]) and so $\scrN^{\ad}_1$ is the union of the $G_1(\dbA_f^{(p)})$-translates of $IM$. So $\scrN^{\ad}_1$ is a regular, formally smooth  $O(G_1,X_1,p)$-scheme and so an integral canonical model of $\Sh_{H^{\ad}}(G^{\ad},X^{\ad})=\Sh_{H_1}(G_1,X_1)$ over $O(G_1,X_1,p)$. As $\scrN^\prime_{2O(G_1,X_1,p)}$ is a pro-\'etale cover of a smooth, quasi-projective $O(G_1,X_1,p)$-scheme (cf. 6.2.1) and as $s_2^\prime$ is a finite morphism, we get that $\scrN_1^{\ad}$ is a pro-\'etale cover of a smooth, quasi-projective $O(G_1,X_1,p)$-scheme. So 6.2 (a) holds. 

If $(G_1,X_1)=(G_1^i,X_1^i)$ has compact factors, then so does $(G_2^\prime,X_2^\prime)$ and so (cf. 6.2.1) $\tilde\scrN^\prime_{O(G_1,X_1,p)}$ is a pro-\'etale cover of a smooth, projective $O(G_1,X_1,p)$-scheme. It is easy to see that this implies that $\scrN^{\ad}_1$ is a pro-\'etale cover of a smooth, projective $O(G_1,X_1,p)$-scheme. So 6.2 (b) also holds. This ends the argument for 6.2 (a) and (b).

\medskip\noindent
{\bf 6.2.4. End of the proof of 1.6 2).} We now come back to the general case of an arbitrary finite set $I_0$ and we will check the following extra thing:

\medskip
{\bf (*)} {\it the normalization of $\scrN_1$ of $\scrN_1^{\ad}$ in $\Sh_{H_1}(G_1,X_1)$ is the integral canonical model of $\Sh_{H_1}(G_1,X_1)$ over $O(G_1,X_1,p)$.} 

\medskip
As in the proof of [Va4, 5.4.1], it suffices to check (*) under the extra assumption that $G_3=G_1$. Based on 6.2 (v), to prove (*) we can also assume that $I_0$ has only one element $i$ and so that $(G_1,X_1)$ is a simple, adjoint Shimura pair. But in this case, (*) follows easily from 6.2 (vii) and the three Cases considered in 6.2.1 to 6.2.3. 

From 6.2 (b) and from (*) we get that if $(G^i_1,X^i_1)$ has compact factors for all $i\in I_0$, then $\scrN_1$ is a pro-\'etale cover of a smooth, projective $O(G_1,X_1,p)$-scheme. This ends the proof of 1.6 2) and so it also ends the proof of 1.6.  

\bigskip\noindent
{\bf 6.3. New examples of N\'eron models.} Let $D$ be a DVR. Let $K_D$ be the field of fractions of $D$. Let $Z_{K_D}$ be a smooth, separated $K_D$-scheme of finite type. We recall (cf. [BLR, p. 12]) that a N\'eron model of $Z_{K_D}$ over $D$ is a smooth, separated $D$-scheme $Z$ of finite type having $Z_{K_D}$ as its generic fibre and which satisfies the following universal property, called the N\'eron mapping property: 

\medskip
{\bf (NMP)} {\it for any smooth $D$-scheme $Y$ and for every $K_D$-morphism $u_{K_D}:Y_{K_D}\to Z_{K_D}$, there is a unique morphism $u:Y\to Z$ of $D$-schemes extending $u_{K_D}$.} 

\medskip
A classical result of N\'eron says that any abelian variety over $K_D$ has a N\'eron model over $D$, cf. [N\'e]. This result has an analogue for the case of $K_D$-torsors of smooth group schemes over $K_D$ of finite type, cf. [BLR, Cor. 4 of p. 158]. On [BLR, p. 15] it is stated that the importance of the notion of N\'eron models ``seems to be restricted" to ``torsors under group schemes". It was a deep insight of Milne which implicitly pointed out that N\'eron models are very important in the study of Shimura varieties, cf. the definitions [Mi2, 2.1, 2.2, 2.5 and 2.9]. The goal of this section is to bring to a concrete fruition Milne's insights. In other words, we will use integral canonical models of certain Shimura varieties to provide large classes of projective varieties over certain $K_D$'s which have projective N\'eron models and which often are not even embeddable into abelian varieties (or torsors of smooth group schemes) over $K_D$.

Let $(G_1,X_1)$ be a Shimura pair of preabelian type such that each simple factor of $(G_1^{\ad},X_1^{\ad})$ has compact factors. Let $p\in\dbN$ be a prime such that the group $G_{1\dbQ_p}$ is unramified. Let $H_1$ be a hyperspecial subgroup of $G_{1\dbQ_p}(\dbQ_p)$. If $p=2$ we assume that the set $S(G_1,X_1,2)$ defined in 1.6 2) is non-empty. Let $O(G_1,X_1,p)$ be as in 1.6 2). Let $\scrN_1$ be the integral canonical model of $\Sh_{H_1}(G_1,X_1)$ over $O(G_1,X_1,p)$, cf. 1.6 2). Let $H_{1,p}$ be a compact, open subgroup of $G_1(\dbA_f^{(p)})$ such that $\scrN_1$ is a pro-\'etale cover of $\scrN_1/H_{1,p}$. The $O(G_1,X_1,p)$-scheme $\scrN_1/H_{1,p}$ is smooth and projective, cf. 1.6 2). 

\medskip\noindent
{\bf 6.3.1. Theorem.} {\it We take $D$ to be an $O(G_1,X_1,p)$-algebra which is a DVR of mixed characteristic $(0,p)$ and of index of ramification at most $\max\{1,p-2\}$. Let 
$$Z:=\scrN_1/H_{1,p}\times_{O(G_1,X_1,p)} D.$$ 
Then the $D$-scheme $Z$ is the N\'eron model of its generic fibre $Z_{K_D}$.}

\medskip
\proof
Let $Y$ be a smooth $D$-scheme. Let $u_{K_D}:Y_{K_D}\to Z_{K_D}$ be a morphism of $K_D$-schemes. The $D$-scheme $Z$ is projective and so there is an open subscheme $U$ of $Y$ such that the complement of $U$ in $Y$ has codimension at most $2$ in $Y$ and $u_{K_D}$ extends uniquely to a morphism $v_U:U\to Z$. The scheme $U_1:=U\times_Z \scrN_{1D}$ is a pro-\'etale cover of $U$ and so from the classical purity theorem (see [SGA1, p. 275]) we get that it extends to a pro-\'etale cover $Y_1$ of $Y$. The $D$-scheme $Y_1$ is healthy regular, cf. 1.2.1. So as $\scrN_1$ has the extension property, the natural morphism $Y_{1K_D}\to\scrN_1$ extends uniquely to a morphism $Y_1\to\scrN_1$ of $D$-schemes. This implies that $v_U$ extends uniquely to a morphism $u:Y\to Z$ of $D$-schemes. This ends the proof.

\medskip\noindent
{\bf 6.3.2. Remark.} If $H_{1,p}$ is small enough, then $\Sh(G_1,X_1)_{\dbC}/H_{1,p}\times H_1$ is of general type (see [Mi1, 1.2 of \S2]). So $Z$ is not among the N\'eron models studied in [BLR]. Often the albanese variety of any connected component $\scrC_1$ of $\Sh(G_1,X_1)_{\dbC}/H_{1,p}\times H_1$ is trivial and so $Z_{\overline{K_D}}$ is not even embeddable into an abelian variety. For instance, this happens if $G^{\ad}_{1\dbR}$ is $SU^{\ad}(a,b)_{\dbR}\times_{\dbR} SU(a+b,0)^{\ad}_{\dbR}$, with $a$, $b>2$ (we have $H^{1,0}(\scrC_1,\dbC)=0$, cf. [CL, p. 11--12]). This remark was first hinted at in [Va2].

\bigskip\noindent
{\bf 6.4. Complements.} {\bf 1)} One can use [Va3] to show that the part of 6.1.2.2 pertaining to $p=2$ and $(G_1,X_1)$ a simple, adjoint Shimura pair of $C_n$ type and having compact factors, can be entirely adapted for the $D_n^{\dbH}$ types. 

{\bf 2)} Let the quadruple $(f,L,v,G_{\dbZ_{(p)}})$ be as in 1.3. We assume that $G_{\dbZ_{(p)}}$ is a reductive group scheme and that one of the conditions (i) and (ii) of 1.6 2) holds; so $p\Ge 3$. Then from 4.2.2 2) and 1.6 1) we get that the triple $(f,L,v)$ is a standard Hodge situation in the sense of [Va5, 1.4]. One can combine 4.2.2 2) with (iii) of 1.6 1) to get a weaker form of this for $p=2$.

{\bf 3)} The arguments of 6.2 can be used to show starting form 4.5 that for $p\Ge 3$ (resp. for $p=2$) and for any triple $(G_1,X_1,H_1)$ with $(G_1,X_1)$ a Shimura pair of preabelian (resp. of abelian) type and with $H_1$ a hyperspecial subgroup of $G_{1\dbQ_p}(\dbQ_p)$, $\Sh_{H_1}(G_1,X_1)$ has a ``good" smooth integral model $\scrN_1$ over $E(G_1,X_1)_{(p)}$. For instance, if $G_1=G_1^{\ad}$, then here by ``good" we mean that a pro-\'etale cover of an open closed subscheme of $\scrN_1$ is a weak integral canonical model of $\Sh_H(G,X)$ over $E(G,X)_{(p)}$, where $(G,X)$ is a Shimura pair of Hodge type having $(G_1,X_1)$ as its adjoint and such that $G_{\dbQ_p}$ is unramified and where $H$ is the unique hyperspecial subgroup of $G_{\dbQ_p}(\dbQ_p)$ satisfying $H^{\ad}=H_1$. 

Warning: the methods of [Va1] and of this paper can be used (in a way entirely similar to the proof of 6.3.1) to show that $\scrN_1$ is a weak integral canonical model only if $\scrN$ is a pro-\'etale cover of a smooth, projective scheme over $E(G,X)_{(p)}$ (to be compared with the proof of [Va1, 6.2.2 a)] where the arguments did appeal to stronger extension properties than the smooth extension property, like the extension property). However, in part II we will show that always $\scrN_1$ is in fact an integral canonical model. 

{\bf 4)} We consider a Shimura pair $(G_1,X_1)$ of abelian type. We assume that $p=2$, that all simple factors of $(G_1^{\ad},X_1^{\ad})$ are of some $D_n^{\dbH}$ type and do not have compact factors and that $G_{1\dbQ_p}$ is unramified. Let $H_1$ be an arbitrary hyperspecial subgroup of $G_{1\dbQ_p}(\dbQ_p)$. Then entirely as in 6.2.2.2 we argue based on 3.1.4 (as a substitute of 6.1), that $\Sh_{H_1}(G_1,X_1)$ has an integral canonical model over $E(G_1,X_1)_{(p)}$.

{\bf 5)} One can easily adapt 6.1 (resp. 3.3 and 6.1) to show that in the context of 1.5, if $p=2$ (resp. if $p\Ge 3$ and each simple factor of $(G^{\ad},X^{\ad})$ has compact factors), then $\scrN=\scrN^\prime$ is a pro-\'etale scheme over a smooth, projective $O_{(v)}$-scheme. Also for $p\Ge 3$ one can combine 1.5 with [Va1, Th. 5.1] to get variants of 1.5 in which in fact we have $\scrN=\scrN^\prime$ without assuming that each simple factor of $(G^{\ad},X^{\ad})$ has compact factors.

{\bf 6)} We now explain why we view 1.5 and its variants of 5) as a generalization of [Mo] for the unramified context with $p>2$. We refer to 1.5 with $p>2$. Let $G^1_{\dbZ_p}$ be as in 1.2 3). Let $G^{2}_{\dbQ_p}$ be the normal, semisimple subgroup of $G^{\der}_{\dbQ_p}$ such that $G^{\ad}_{\dbQ_p}$ is naturally isomorphic to $G^{1\ad}_{\dbQ_p}\times_{\dbQ_p} G^{2\ad}_{\dbQ_p}$. Let $g_1$, ..., $g_n\in G^{2}_{\dbQ_p}(\dbQ_p)$. Let 
$$\tilde H:=H\cap g_1Hg_1^{-1}\cap...\cap g_nHg_n^{-1}.$$
 Using Hodge quasi products as in [Va4, 2.4] and Segre embeddings as in [Va1, Example 3 of 2.5], one gets that there is an injective map $\tilde f:(G,X)\hookrightarrow (GSp(W,\psi),S)$ and a $\dbZ$-lattice $\tilde L$ of $\tilde W$ such that the following things hold:

\medskip
-- {\it we have a perfect alternating form $\psi:\tilde L\otimes_{\dbZ} \tilde L\to\dbZ$, $\tilde H=\{h\in G(\dbQ_p)|h(\tilde L\otimes_{\dbZ} \dbZ_p)=\tilde L\otimes_{\dbZ} \dbZ_p\}$, and $G^1_{\dbZ_p}$ is a closed subgroup of the Zariski closure of $G_{\dbQ_p}$ in $GL(\tilde L\otimes_{\dbZ} \dbZ_p)$.} 

\medskip
From 1.5 we get that $\Sh_{\tilde H}(G,X)$ has a weak integral model $\tilde\scrN$ over $O_{(v)}$. As a scheme $\tilde\scrN$ is a quasi-finite scheme over a product $\scrN^n$ of $n$-copies of the $O_{(v)}$-scheme $\scrN$ (cf. 3.1.1 1)) in such a way that the natural projection morphisms $\tilde\scrN\to\scrN$ are formally \'etale (this follows easily from 1.4 2) and 4.3 (iii)). If each simple factor of $(G^{\ad},X^{\ad})$ has compact factors, then $\tilde\scrN$ is a pro-\'etale scheme over a smooth, projective $O_{(v)}$-scheme (cf. 5)) and the natural projection morphisms $\tilde\scrN\to\scrN$ are \'etale covers. 

Let $H_{\infty}:=\cap_{g\in G^2_{\dbQ_p}(\dbQ_p)} gHg^{-1}$. The projective limit of $\tilde\scrN$'s is the weak integral canonical model of $\Sh_{H_{\infty}}(G,X)$ over $O_{(v)}$. If $\bigl(G_{\dbQ_p}/Z^0(G_{\dbQ_p})G^1_{\dbQ_p}\bigr)(\dbQ_p)$ has no compact subgroups normalized by $G^2_{\dbQ_p}(\dbQ_p)$, then $H_{\infty}=H\cap (Z^0(G_{\dbQ_p})G^1_{\dbQ_p})(\dbQ_p)$. 

One can also vary $L$ subject to the fact that $\psi$ induces a perfect form on it and $G^1_{\dbZ_p}$ is a reductive subgroup of $GL(L\otimes_{\dbZ} \dbZ_p)$. So in the previous paragraph, often one can replace $H_{\infty}$ by many other compact subgroups of it containing $G^1_{\dbZ_p}(\dbZ_p)$. A quite similar approach was used in [Mo] for the particular case when $(G_1^{\ad},X_1^{\ad})$ is a Shimura curve and $v_1^{\ad}$ is a prime such that $G_1^{\ad}$ splits over $E(G_1^{\ad},X_1^{\ad})_{v_1^{\ad}}$. Warning: in our case $v_1$ is always unramified over an odd prime $p$ while loc. cit. does not require this restriction.

\bigskip
\references{37}
{\nspace{

\Ref[BB]
W. Baily and A. Borel,
\sl Compactification of arithmetic quotients of bounded
symmetric domains,
\rm Ann. of Math. {\bf 84} (1966), p. 442--528.
\Ref[Be] P. Berthelot, 
\sl Comologie cristalline des sch\'emas de caract\'eristique $p>0$, 
\rm LNM {\bf 407} (1974), Springer--Verlag.
\Ref[BHC]
A. Borel and Harish-Chandra,
\sl Arithmetic subgroups of algebraic subgroups,
\rm Ann. of Math. 75 (1962), 485--535.
\Ref[Bl]
D. Blasius,
\sl A $p$-adic property of Hodge cycles on abelian varieties, 
\rm Proc. Symp. Pure Math. 55, Part 2, p. 293--308.

\Ref[BLR]
S. Bosch, W. L\"utkebohmert, and M. Raynaud,
\sl N\'eron models,
\rm Springer-Verlag, 1990.

\Ref[BrT]
F. Bruhat and J. Tits, 
\sl Groupes r\'eductifs sur un corps local,
\rm Publ. Math. I.H.E.S. {\bf 60} (1984), p. 5--184.

\Ref[Bo]
A. Borel,
\sl Linear algebraic groups,
\rm Grad. Texts Math. 126, Springer-Verlag, 1991.

\Ref[CL] L. Clozel, 
\sl Produits dans la cohomologie holomorphe des vari\'et\'es de Shimura II: Calculs et applications, 
\rm J. reine angew. Math. {\bf 444} (1993), p. 1--15.

\Ref[dJ1] J. de Jong, 
\sl Crystalline Dieudonn\'e module theory via formal and rigid geometry, 
\rm Publ. Math. IHES {\bf 82} (1995), p. 5--96.

\Ref[dJ2] 
J. de Jong, 
\sl Homomorphisms of Barsotti-Tate groups and crystals in positive characteristic,
\rm Inv. Math. {\bf 134} (1998), no. 2, p. 301--333.

\Ref[De1]
P. Deligne,
\sl Travaux de Shimura,
\rm S\'eminaire  Bourbaki 389, LNM 244 (1971), p. 123--163.

\Ref[De2]
P. Deligne,
\sl Vari\'et\'es de Shimura: Interpr\'etation modulaire, et
techniques de construction de mod\`eles canoniques,
\rm Proc. Symp. Pure Math. 33, Part 2 (1979), p. 247--290.

\Ref[De3]
P. Deligne,
\sl Hodge cycles on abelian varieties,
\rm Hodge cycles, motives, and Shimura varieties, LNM 900, Springer-Verlag, 1982, p. 9--100.

\Ref[Dr] V. G. Drinfeld, 
\sl Elliptic modules, 
\rm Math. USSR Sbornik {\bf 23} (1974), p. 561--592.

\Ref[Fa]
G. Faltings,
\sl Integral crystalline cohomology over very ramified
valuation rings,
\rm J. of Am. Math. Soc., Vol. 12, No. 1 (1999), p. 117--144.

\Ref[FC]
G. Faltings and C.-L. Chai,
\sl Degeneration of abelian varieties,
\rm Springer, Heidelberg, 1990.

\Ref[Fo] 
J.-M. Fontaine, 
\sl Le corps des p\'eriodes p-adiques, 
\rm J. Ast\'erisque {\bf 223}, Soc. Math. de France, Paris 1994, p. 59--101.

\Ref[Ha]
G. Harder,
\sl \"Uber die Galoiskohomologie halbeinfacher Matrizengruppen II,
\rm  Math. Z. 92 (1966), p. 396--415.

\Ref[Har]
R. Hartshorne,
\sl Algebraic geometry,
\rm Grad. Text Math. 52, Springer--Verlag, 1977.

\Ref[He]
S. Helgason,
\sl Differential geometry, Lie groups, and symmetric spaces,
\rm Academic Press, New-York, 1978.

\Ref[HT]
M. Harris and R. Taylor,
\sl On the geometry and cohomology of some simple Shimura varieties,
\rm preprint, Aug. 30, 1999.

\Ref[Hu]
J. E. Humphreys, 
\sl Introduction to Lie algebras and representation theory,
\rm Grad. Texts Math. {\bf 9}, Springer--Verlag 1975.

\Ref[Ja]
J. C. Jantzen,
\sl Representations of algebraic groups,
\rm Academic Press, 1987.

\Ref[Ko]
R. E. Kottwitz,
\sl Points on some Shimura Varieties over finite fields,
\rm J. of Am. Math. Soc., Vol. 5, No. 2, 1992, p. 373--444.

\Ref[La]
R. Langlands,
\sl Some contemporary problems with origin in the Jugendtraum,
\rm Mathematical developments arising from Hilbert's problems, Am. Math. Soc., Providence, RI, 1976, p. 401--418.

\Ref[Lan]
S. Lang,
\sl Algebraic Number Theory,
\rm Grad. Texts Math. {\bf 110}, Springer--Verlag 1994.

\Ref[LR]
R. Langlands and M. Rapoport,
\sl Shimuravarietaeten und Gerben, 
\rm J. reine angew. Math. 378 (1987), p. 113--220.

\Ref[Ma]
H. Matsumura,
\sl Commutative algebra,
\rm The Benjamin/Cummings Publishing Co., Inc., 1980.

\Ref[Mi1]
J. S. Milne,
\sl Canonical models of (mixed) Shimura varieties and automorphic vector bundles, 
\rm  Automorphic Forms, Shimura varieties and L-functions, Vol. I, Perspectives in Math., Vol. 10, Acad. Press 1990.

\Ref[Mi2]
J. S. Milne,
\sl The points on a Shimura variety modulo a prime of good
reduction,
\rm The Zeta function of Picard modular surfaces, Les Publications CRM, Montreal 1992, p. 153--255.

\Ref[Mi3]
J. S. Milne,
\sl Shimura varieties and motives,
\rm Proc. Symp. Pure Math. {\bf 55} (1994), Part 2, p. 447--523.

\Ref[Mi4]
J. S. Milne,
\sl Descent for Shimura varieties,
\rm Michigan Math. J. 46 (1999), p. 203--208.

\Ref[Mo]
Y. Morita,
\sl Reduction mod $\grB$ of Shimura curves,
\rm Hokkaido Math. J. {\bf 10} (1981), p. 209--238.

\Ref[MFK]
D. Mumford, J. Fogarty, and F. Kirwan,
\sl Geometric invariant theory,
\rm Ergeb. der Math. und ihrer Grenz., Vol. 34, Spinger-Verlag, third enlarged edition, 1994.

\Ref[MS]
J. S. Milne and K.-y. Shih,
\sl Conjugates of Shimura varieties,
\rm  Hodge cycles,
motives, and Shimura varieties, LNM {\bf 900}, Springer-Verlag, 1982, p. 280--356.

\Ref[Na]
M. Nagata,
\sl Imbedding of an abstract variety in a complete variety,
\rm J. Math. Kyoto Univ., {\bf 2} (1962), p. 1--10.

\Ref[N\'e] 
A. N\'eron, 
\sl Mod\`eles minimaux des vari\'et\'es ab\'eliennes, 
\rm Publ. Math. IHES {\bf 21} (1964).

\Ref[Oo]
F. Oort,
\sl Subvarieties of moduli spaces,
\rm Inv. Math. {\bf 24} (1975), p. 95--119.

\Ref[Re]
H. Reimann,
\sl The semi-simple zeta function of quaternionic Shimura varieties,
\rm LNM {\bf 1657} (1997), Springer--Verlag, Berlin, viii + 143 p.

\Ref[Sa]
I. Satake,
\sl Holomorphic imbeddings of symmetric domains into a Siegel space,
\rm  Am. J. Math. 87 (1965), p. 425--461.

\Ref[Saa]
N. Saavedra Rivano,
\sl Cat\'egories tannakiennes,
\rm LNM 265, Springer-Verlag, 1972.

\Ref[SGA1]
A. Grothendieck \'et al.,
\sl Rev\^etements \'etales et groupe fondamental,
\rm LNM {\bf 224}, 1971, Springer--Verlag.

\Ref[SGA3]
A. Grothendieck \'et al.,
\sl Sch\'emas en groupes,
\rm LNM 152-3, Vol. II-III, 1970.

\Ref[Ti1]
J. Tits,
\sl Classification of algebraic semisimple groups,
\rm Proc. Sympos. Pure Math., Vol. 9, A.M.S., Providence, R.I., 1966, p. 33--62.

\Ref[Ti2]
J. Tits,
\sl Reductive groups over local fields, 
\rm Proc. Symp. Pure Math. 33, Part 1, 1979, p. 29--69.

\Ref[Va1]
A. Vasiu,
\sl Integral canonical models of Shimura varieties of preabelian type,
\rm Asian J. Math., Vol. 3, No. 2, p. 401--518, June 1999.

\Ref[Va2] 
A. Vasiu,
\sl A purity theorem for abelian schemes,
\rm accepted for publ. in Michigan Math. J.

\Ref[Va3] 
A. Vasiu,
\sl Integral models in mixed characteristic (0,2) of Shimura varieties of PEL type,
\rm submitted for publ. on 5/16/02. 

\Ref[Va4]
A. Vasiu,
\sl The Mumford--Tate Conjecture and Shimura Varieties, Part I,
\rm submitted for publ. on 12/11/02, improved version 9/15/02. 

\Ref[Va5] 
A. Vasiu,
\sl Generalized Serre--Tate Ordinary Theory,
\rm submitted for publ. on 3/3/03. 

\Ref[Va6] 
A. Vasiu,
\sl A motivic conjecture of Milne,
\rm submitted for publ. on 8/20/03.

\Ref[Vo]
P. Vojta,
\sl Nagata's embedding theorem,
\rm preprint (appendix to planned book).

\Ref[Zi]
T. Zink,
\sl Isogenieklassen von Punkten von Shimuramannigfaltigkeiten mit Werten in einem endlichen K\"orper,
\rm Math. Nachr. {\bf 112} (1983), p. 103--124.

\Ref[Wi]
J.-P. Wintenberger,
\sl Un scindage de la filtration de Hodge pour certaines variet\'es algebriques sur les corps locaux,
\rm Ann. of Math. (2) 119 (1984), p. 511--548.

}}

\enddocument